%% file: spring-dimer-current.tex
\documentclass[12pt]{amsart}

\allowdisplaybreaks
\numberwithin{equation}{subsection}
\everymath{\displaystyle}

\newtheorem{theorem}{Theorem}[section]
\newtheorem{proposition}[theorem]{Proposition}
\newtheorem{lemma}[theorem]{Lemma}
\newtheorem{corollary}[theorem]{Corollary}
\newtheorem{remark}[theorem]{Remark}

\DeclareSymbolFont{bbold}{U}{bbold}{m}{n}
\DeclareSymbolFontAlphabet{\mathbbold}{bbold}

\DeclareFontFamily{OT1}{pzc}{}
\DeclareFontShape{OT1}{pzc}{m}{it}{<-> s * [1.10] pzcmi7t}{}
\DeclareMathAlphabet{\mathpzc}{OT1}{pzc}{m}{it}

\usepackage{adjustbox}
\usepackage{amsmath}
\usepackage{amsfonts}
\usepackage{amssymb}
\usepackage{enumitem}\setenumerate{label={\bf(\alph*)},leftmargin=0pt,listparindent=\parindent,itemsep=10pt,topsep=10pt,parsep=0pt,widest=a,fullwidth}
\usepackage{mathtools}
\usepackage{pgfplots}
\usepackage{subfig}
\usepackage{wrapfig}

\usepackage{hyperref} 

\newcommand{\an}{\mathfrak{A}}

\renewcommand{\b}{\mathbf{B}}
\newcommand{\bignorm}[1]{\big\|#1\big\|}
\DeclareMathOperator{\boot}{boot}
\newcommand{\bunderbrace}[2]{\underbrace{#1}_{\dst{#2}}}

\newcommand{\cep}{c_{\ep}}

\newcommand{\dby}[2]{\frac{d{#1}}{d{#2}}}
\newcommand{\ds}{\ ds}
\newcommand{\dst}[1]{{\displaystyle{#1}}}
\newcommand{\dt}{\ dt}
\newcommand{\dX}{\ dX}

\newcommand{\ep}{\epsilon}
\newcommand{\epbar}{\overline{\ep}}

\newcommand{\ft}{\mathfrak{F}}
\newcommand{\fb}{\mathbf{f}}

\renewcommand{\hat}[1]{\widehat{#1}}
\newcommand{\hb}{{\bf{h}}}

\newcommand{\ip}[2]{\left\langle#1,#2\right\rangle}
\DeclareMathOperator{\incr}{inc}
\newcommand{\ind}{\scalebox{1.15}{$\mathbbold{1}$}}

\renewcommand{\j}{\mathpzc{j}}
\newcommand{\jb}{\mathbf{j}}

\DeclareMathOperator{\Kappa}{K}

\renewcommand{\l}{\mathpzc{l}}
\DeclareMathOperator{\lip}{lip}
\renewcommand{\L}{\mathcal{L}}

\DeclareMathOperator{\map}{map}
\newcommand{\medsum}{\textstyle{\sum}}
\renewcommand{\mod}{\text{mod}}

\newcommand{\nanorhs}{\mathfrak{R}}
\newcommand{\nbar}{\overline{N}}
\newcommand{\norm}[1]{\left\|#1\right\|}
\newcommand{\nsyst}{\mathfrak{N}}
\newcommand{\nsystb}{\boldsymbol{\nsyst}}

\newcommand{\pb}{{\bf{p}}}
\DeclareMathOperator{\per}{per}

\newcommand{\quadword}[1]{\qquad\text{ #1 }\qquad}
\newcommand{\qstar}{q_{\star}}

\DeclareMathOperator{\sech}{sech}
\newcommand{\set}[2]{\!\left\{ #1 \ \middle| \ #2\right\}}

\newcommand{\tbar}{\overline{t}}
\newcommand{\tPi}{\tilde{\Pi}}
\DeclareMathOperator{\TT}{T}
\newcommand{\tv}{\tilde{v}}
\newcommand{\tvb}{\tilde{\mathbf{v}}}

\newcommand{\ub}{\mathbf{u}}
\newcommand{\ubar}{\overline{u}}

\newcommand{\vb}{\mathbf{v}}

\newcommand{\alphab}{\boldsymbol{\alpha}}
\newcommand{\betab}{\boldsymbol{\beta}}
\newcommand{\etab}{\boldsymbol{\eta}}
\newcommand{\nub}{\boldsymbol{\nu}}
\newcommand{\phib}{\boldsymbol{\phi}}
\newcommand{\Phib}{\boldsymbol{\Phi}}
\newcommand{\psib}{\boldsymbol{\psi}}
\newcommand{\Psib}{\boldsymbol{\Psi}}
\newcommand{\sigmab}{\boldsymbol{\sigma}}
\newcommand{\thetab}{\boldsymbol{\theta}}
\newcommand{\Thetab}{\boldsymbol{\Theta}}
\newcommand{\varphib}{\boldsymbol{\varphi}}

\renewcommand{\kappa}{\varkappa}

\renewcommand{\tilde}[1]{\widetilde{#1}}

\newcommand{\tlambda}{\tilde{\lambda}}
\newcommand{\tLambda}{\tilde{\Lambda}}
\newcommand{\tmu}{\tilde{\mu}}
\newcommand{\trho}{\tilde{\varrho}}

\newcommand{\tvarpi}{\tilde{\varpi}}
\newcommand{\txi}{\tilde{\xi}}

\newcommand{\tJ}{\tilde{J}}
\newcommand{\tL}{\tilde{L}}

\newcommand{\A}{\mathcal{A}}
\newcommand{\B}{\mathcal{B}}
\newcommand{\C}{\mathcal{C}}
\newcommand{\D}{\mathcal{D}}
\newcommand{\xtr}{\mathcal{E}}
\newcommand{\F}{\mathcal{F}}
\newcommand{\G}{\mathcal{G}}
\renewcommand{\H}{\mathcal{H}}
\newcommand{\I}{\mathcal{I}}
\newcommand{\K}{\mathcal{K}}
\renewcommand{\L}{\mathcal{L}}
\newcommand{\M}{\mathcal{M}}
\newcommand{\nl}{\mathcal{N}}
\renewcommand{\O}{\mathcal{O}}
\renewcommand{\P}{\mathcal{P}}
\newcommand{\cb}{\mathcal{Q}}
\newcommand{\rhs}{\mathcal{R}}
\newcommand{\T}{\mathcal{T}}
\newcommand{\U}{\mathcal{U}}
\newcommand{\W}{\mathcal{W}}
\newcommand{\X}{\mathcal{X}}
\newcommand{\Y}{\mathcal{Y}}

\newcommand{\R}{\mathbb{R}}
\newcommand{\N}{\mathbb{N}}
\newcommand{\Z}{\mathbb{Z}}

\makeatletter
\DeclareRobustCommand\widecheck[1]{{\mathpalette\@widecheck{#1}}}
\def\@widecheck#1#2{%
    \setbox\z@\hbox{\m@th$#1#2$}%
    \setbox\tw@\hbox{\m@th$#1%
       \widehat{%
          \vrule\@width\z@\@height\ht\z@
          \vrule\@height\z@\@width\wd\z@}$}%
    \dp\tw@-\ht\z@
    \@tempdima\ht\z@ \advance\@tempdima2\ht\tw@ \divide\@tempdima\thr@@
    \setbox\tw@\hbox{%
       \raise\@tempdima\hbox{\scalebox{1}[-1]{\lower\@tempdima\box
\tw@}}}%
    {\ooalign{\box\tw@ \cr \box\z@}}}
\makeatother

\usepackage{tikz}
\usetikzlibrary{positioning}
\usetikzlibrary{patterns}
\usetikzlibrary{shapes}
\usetikzlibrary{matrix, arrows,decorations.pathmorphing}
\usetikzlibrary{decorations.markings}
\usetikzlibrary{fadings}
\usetikzlibrary{arrows.meta,bending}
\usepackage{tikzscale}

\newcommand{\spring}[3]{
\draw[line width =.055cm] (#3,0)
--(#3+#1,#2)
--(#3+3*#1,-#2)
--(#3+5*#1,#2)
--(#3+7*#1,-#2)
--(#3+9*#1,#2)
--(#3+11*#1,-#2)
--(#3+12*#1,0)
}

\tikzstyle{springG} = [decorate,decoration={aspect=.65, segment length=.56cm, amplitude=.3cm,coil},line width = .055cm] 

\newcommand{\leftpointycorrectionfactor}[1]{
\fill (#1,0)--(.0235+#1,-0.0141)--(.0155+#1,-.75pt)--(#1,-.75pt)--cycle;
\draw[line width = .025pt] (#1,-.75pt)--(#1,0)--(.0235+#1,-.0141);
}

\newcommand{\rightpointycorrectionfactor}[1]{
\fill (#1,0)--(#1,.75pt)--(#1-.0155,.75pt)--(#1-.0235,.0141)--cycle;
\draw[line width = .025pt] (#1,.75pt)--(#1,0)--(#1-.0235,.0141);
}

\newcommand{\leftcurlycorrectionfactor}[1]{
\fill (#1,-.01) rectangle (#1+.026,-.75pt);
\draw[line width = 1pt] (#1,-.01) -- (#1,-.75pt);
}

\newcommand{\rightcurlycorrectionfactor}[1]{
\fill (#1-.01,.025) rectangle (#1+.01,-.025);
}

\title[Nanpteron-stegoton traveling waves in spring dimer FPUT lattices]{Nanopteron-stegoton traveling waves in spring dimer Fermi-Pasta-Ulam-Tsingou lattices}
\author{Timothy E. Faver}
\address{Department of Mathematics, Drexel University, 3141 Chestnut St, Philadelphia, PA 19104, tef36@drexel.edu}

\begin{document}
\keywords{FPU, FPUT, nonlinear hetergeneous lattice, dimer, solitary traveling wave, periodic traveling wave, singular perturbation, nanopteron, stegoton, composition operator}

\begin{abstract}
We study the existence of traveling waves in a spring dimer Fermi-Pasta-Ulam-Tsingou (FPUT) lattice.  
This is a one-dimensional lattice of identical particles connected by alternating nonlinear springs. 
Following the work of Faver and Wright on the mass dimer, or diatomic, lattice, we find that the lattice equations in the long wave regime are singularly perturbed and apply a method of Beale to produce nanopteron traveling waves with wave speed slightly greater than the lattice's speed of sound. 
The nanopteron wave profiles are the superposition of an exponentially decaying term (which itself is a small perturbation of a KdV $\sech^2$-type soliton) and a periodic term of very small amplitude.
Generalizing our work in the diatomic case, we allow the nonlinearity in the spring forces to have the more complicated form ``quadratic plus higher order terms.''  
This necessitates the use of composition operators to phrase the long wave problem, and these operators require delicate estimates due to the characteristic superposition of function types from Beale's ansatz.  
Unlike the diatomic case, the value of the leading order term in the traveling wave profiles alternates between particle sites, so that the spring dimer traveling waves are also ``stegotons,'' in the terminology of LeVeque and Yong.
This behavior is absent in the mass dimer and confirms the approximation results of Gaison, Moskow, Wright, and Zhang for the spring dimer.
\end{abstract}

\maketitle

\input{"section_introduction"}
\input{"section_tw_equations"}

\input{"section_periodic_solutions"}
\input{"section_nanopteron_equations"}
\input{"section_existence_properties_solutions"}
\appendix
\input{"appendix_function_spaces_master"}
\input{"appendix_existence_of_periodic_solutions"}
\input{"appendix_nanopteron_estimates_master"}

\bibliographystyle{alpha}
\bibliography{spring-dimer-bib}
\end{document}

%% file: section_introduction.tex
\section{Introduction}

\subsection{The lattice equations}
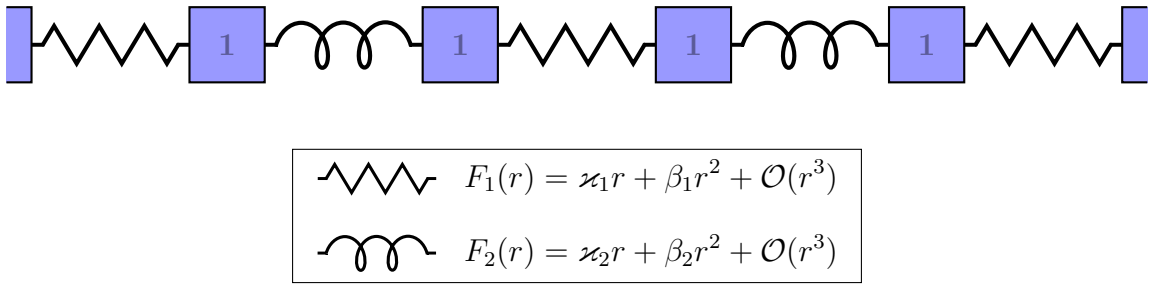
\begin{figure}
$$
\input{fig_spring_dimer}
$$

$$
\input{fig_spring_dimer_legend}
$$
\caption{The spring dimer lattice}
\label{spring dimer fig}
\end{figure}

We may assemble a one-dimensional lattice by placing infinitely many particles on a horizontal line and connecting each particle to the particles on its immediate left and right by springs. Such a construct is a Fermi-Pasta-Ulam-Tsingou (FPUT) lattice \cite{fput-original}, \cite{dauxois}. Depending on the material properties that we ascribe to the particles and springs --- that is, the masses of the particles and the forces that the springs exert when stretched --- we can vary the behavior of the lattice considerably.  In this article, we assume that the particles all have identical mass and that the springs alternate, as detailed below and sketched in Figure \ref{spring dimer fig}.

We index the particles, their masses, and the springs by integers $j \in \Z$.  The $j$th spring connects the $j$th particle (on the left) with the $(j+1)$st particle (on the right).  All particles have the same mass $m > 0$, but we vary two properties of the springs.  First, the equilibrium length of the $j$th spring is $\ell_j$, where
$$
\ell_j 
= \begin{cases}
\ell_1 > 0 , &j \text{ is odd} \\
\ell_2 > 0, &j \text{ is even.}
\end{cases}
$$
Second, the $j$th spring exerts a force $F_j(r)$ when stretched a length $r$ from its equilibrium length, where
$$
F_j(r)
=
\begin{cases}
\kappa_1r+\beta_1r^2+r^3\nbar_1(r), &j \text{ is odd} \\
\kappa_2r+\beta_2r^2+r^3\nbar_2(r), &j \text{ is even.}
\end{cases}
$$
The linear spring force coefficients satisfy
$$
0 < \kappa_2 < \kappa_1 
$$
and the quadratic coefficients $\beta_1$ and $\beta_2$ are nonzero.  Last, we assume $\nbar_1$, $\nbar_2 \in \C^{\infty}(\R)$.  
We denote the position of the $j$th particle at time $\tbar$ by $\ubar_j(\tbar)$.  Then Newton's law implies that the position functions satisfy
\begin{equation}\label{newton}
m_j\frac{d^2\ubar_j}{d\tbar^2} 
= F_{j}(\ubar_{j+1}-\ubar_j-\ell_j)-F_{j-1}(\ubar_j-\ubar_{j-1}-\ell_{j-1}).
\end{equation}
We call this species of lattice the \textit{spring dimer} lattice, as opposed to the \textit{mass dimer} or \textit{diatomic} lattice consisting of alternating particles of two different masses and only one kind of spring, which has been more frequently studied \cite{faver-wright}, \cite{hw}, \cite{betti-pelinovsky}, \cite{cb-etal}, \cite{qin}.  

Both of these lattices are spatially heterogeneous generalizations of the monatomic species treated comprehensively by Friesecke and Wattis \cite{friesecke-wattis} and Friesecke and Pego \cite{friesecke-pego1,friesecke-pego2,friesecke-pego3,friesecke-pego4}. In these articles, it is shown that the monatomic lattice, which has only one kind of particle and one kind of spring force, possesses solitary traveling waves.  The mathematical and physical model provided by a one-dimensional lattice has a rich history of theory and applications that predates these articles, and even the Fermi-Pasta-Ulam-Tsingou experiments, by centuries; see Brillouin \cite{brill} for a pre-1950s history and Pankov \cite{pankov} for a contemporary mathematical overview.  The articles \cite{kevrekidis} and \cite{cpkd} explore a plethora of modern applications of lattices.

\subsection{Nondimensionalization} 
We simplify the original lattice equations \eqref{newton} in several ways.  First, we can eliminate the spring lengths $\ell_1$ and $\ell_2$ by writing
$$
\ubar_{2j+1} 
= \breve{u}_{2j+1} +j\ell_1 + j\ell_2 
\quadword{and} 
\ubar_{2j} 
= \breve{u}_{2j} + j\ell_1+(j-1)\ell_2.
$$
Then
$$
\ubar_{2j+1}-\ubar_{2j}-\ell_1
= \breve{u}_{2j+1} - \breve{u}_{2j}
\quadword{and}
\ubar_{2j}-\ubar_{2j-1}-\ell_2
= \breve{u}_{2j} - \breve{u}_{2j-1}.
$$
We next express the system in terms of relative displacements.  We set
$$
\breve{r}_j  
:= \breve{u}_{j+1}-\breve{u}_j.
$$
Then \eqref{newton} becomes
\begin{equation}\label{no spring lengths}
\begin{cases}
m\dby{^2\breve{u}_{2j+1}}{\tbar^2}
= F_1(\breve{r}_{2j+1})-F_2(\breve{r}_{2j}) \\
\\
m\dby{^2\breve{u}_{2j}}{\tbar^2}
=  F_2(\breve{r}_{2j})-F_1(\breve{r}_{2j-1}).
\end{cases}
\end{equation}
Now we make the rescaling
\begin{equation}\label{rel disp nondimen}
\breve{u}_j(\tbar) 
= a_1u_j(a_2\tbar)
\quadword{and}
r_j := u_{j+1} - u_j
\end{equation}
where $u_j = u_j(t)$ and 
$$
a_1 := \frac{\kappa_2}{\beta_2}
\quadword{and}
a_2 :=\sqrt{\frac{\kappa_2}{m}}.
$$
After canceling some common factors and defining 
$$
\kappa := \frac{\kappa_1}{\kappa_2} > 1,
\qquad
\beta := \frac{\beta_1}{\beta_2} \ne 0,
\quadword{and}
N_j(r) := \frac{a_1^2}{\kappa_2}\nbar_j(a_1r), \ j = 1,2,
$$
we convert the system \eqref{no spring lengths} to the nondimensionalized equations
$$
\begin{cases}
\ddot{u}_{2j+1} 
= \kappa r_{2j+1} + \beta r_{2j+1}^2+r_{2j+1}^3N_1(r_{2j+1})
-  r_{2j}-r_{2j}^2 - r_{2j}^3N_{2}(r_{2j}) \\
\\
\ddot{u}_{2j} 
= r_{2j} + r_{2j}^2+r_{2j}^3N_2(r_{2j})
-\kappa r_{2j-1}-\beta r_{2j-1}^2 - r_{2j-1}^3N_1(r_{2j-1}).
\end{cases}
$$
From these we can compute the equations of motion solely in terms of relative displacement:
\begin{equation}\label{nondimen lattice eqns}
\left\{\!\!\begin{tabular}{rl}
$\ddot{r}_{2j+1}$
&= $-2\left(\kappa{r}_{2j+1}+\beta{r}_{2j+1}^2+r_{2j+1}^3N_1(r_{2j+1})\right)
+\left((r_{2j+2}+r_{2j+2}^2+r_{2j+2}^3N_2(r_{2j+2})\right)$ \\
&$+\left(r_{2j}+r_{2j}^2+r_{2j}^3N_2(r_{2j})\right)$ \\
\\
$\ddot{r}_{2j}$
&= $-2\left(r_{2j}+r_{2j}^2+r_{2j}^3N_2(r_{2j})\right) 
+ \left(\kappa{r}_{2j+1}+\beta{r}_{2j+1}^2 + r_{2j+1}^3N_1(r_{2j+1})\right)$ \\
&$+\left(\kappa{r}_{2j-1}+\beta{r}_{2j-1}^2+r_{2j-1}^3N_1(r_{2j-1})\right).$
\end{tabular}\right.
\end{equation}

\subsection{Main results} 
\begin{figure}
$$
\input{fig_nanopteron}
$$
\caption{The nanopteron ($\O(\ep^{\infty}) = $ ``small beyond all orders of $\ep$'')}
\label{nanopteron fig}
\end{figure}
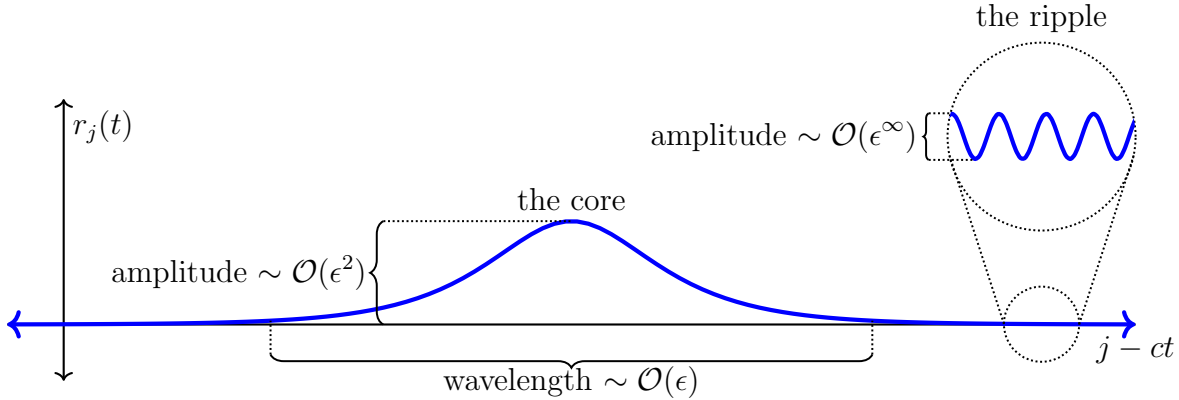

We prove a result analogous to our conclusions for the mass dimer lattice in \cite{faver-wright}: the equations \eqref{nondimen lattice eqns} for relative displacements possess \textit{nanopteron} traveling wave solutions.  Following Boyd \cite{boyd}, the nanopteron is a wave that is asymptotic to a periodic oscillation at $\pm\infty$, where the amplitude of the oscillation, or ``ripple,'' is small beyond all orders of the amplitude of the ``core'' of the wave.  We state informally our main result in Theorem \ref{main theorem} below, provide a precise version in Corollary \ref{main corollary}, and sketch the nanopteron in Figure \ref{nanopteron fig}. 

\begin{theorem}\label{main theorem}
Suppose that the spring force coefficient ratios satisfy $\kappa > 1$, $\beta \ne 0$, and $\beta \ne -\kappa^3$.  Then for wave speeds $c \ge c_{\kappa} := \sqrt{2\kappa/(1+\kappa)}$, the nondimensionalized lattice equations \eqref{nondimen lattice eqns} possess traveling wave solutions of the form $r_j = v_j^{\ep} + p_j^{\ep},$ where
\begin{itemize}[leftmargin=*]
\item $\ep = \sqrt{c^2-c_{\kappa}^2}$ is the distance between the wave speed $c$ and the ``speed of sound'' $c_{\kappa}$.
\item $v_j^{\ep}$ is a small perturbation of a $\sech^2$-type function that solves a KdV traveling wave equation whose coefficients depend on $\kappa$ and $\beta$.  This is the ``core'' of the nanopteron, and its amplitude is roughly $\ep^2$.  
\item $p_j^{\ep}$ is a periodic function whose amplitude is small beyond all orders of $\ep$ and whose frequency is $\O(1)$ in $\ep$.  This is the ``ripple'' of the nanopteron.
\end{itemize}
\end{theorem}

Our analysis follows closely our previous work on the mass dimer lattice, which in turn was based on the approach of Beale for capillary-gravity waves in \cite{beale2} and of Amick and Toland for a related model equation in \cite{amick-toland}.  In \cite{faver-wright}, the fundamental physical parameter for the mass dimer lattice was the mass ratio $w = m_1/m_2 > 1$.  The analogue to $w$ for the spring dimer problem is of course the linear spring force ratio $\kappa = \kappa_1/\kappa_2 > 1$.  To our great convenience, $\kappa$ appears in place of $w$ in several key functions and operators, chiefly in \eqref{eigenvalues}, so that we may import more or less directly many of the key strategies and estimates from our work on the mass dimer.  There are, naturally, some technical differences here, chiefly in the set-up of the long wave limit and the nanopteron equations, but these we quickly overcome.  Notably, while we could solve the mass dimer problem for all mass ratios $w > 1$, we do not solve the spring dimer problem for all linear force ratios $\kappa > 1$ and quadratic force ratios $\beta \ne 0$.  Instead, we require $\beta \ne -\kappa^3$ to preserve the validity of a KdV approximation in the long wave limit.  The solutions that we do find exhibit the stegoton-type behavior indicated in \cite{gmwz}: roughly speaking, the relative displacements for even spatial index $j$ carry an extra factor of $\kappa$ that the displacements for odd $j$ do not, and this results in a ``jagged'' appearance in graphs of approximations to the displacements.

We also explore a generalization that we raised, but did not address rigorously, in the mass dimer problem.  There, we considered spring forces of the form $F(r) = r+r^2$.  We have now permitted higher-order terms in each of the forces, i.e., $F_j(r) = \kappa_j{r} + \beta_j{r}^2 + \O(r^3)$. This is how \cite{friesecke-pego1} handled their monatomic lattice with higher-order forces from the start.  Like us, they used composition operator estimates to glean critical information about their problem.  We do the same, but, following Beale, our ansatz for the traveling wave profiles for \eqref{nondimen lattice eqns} is the superposition of an $L^2$ function and and $L^{\infty}$ function.  In general, such functions are neither $L^2$ nor $L^{\infty}$, so we have no common function space in which to work; moreover, we need certain peculiar estimates on our compositions that are unlikely to exist in standard treatments of these operators.  All this is to say that the higher-order terms introduce numerous challenges at a technical, though not quite conceptual, level.  

Last, although the underlying linear operators in these nonlinear composition operators arise explicitly from the structure of the spring dimer problem, algebraically they are identical to the composition operators that would surface in the mass dimer problem with higher-order forces.  So, our work effectively addresses the problem of higher order forces from mass dimer perspective, as well.

%% file: fig_spring_dimer.tex
\begin{tikzpicture}[thick]
\def\circleradius{.35};
\def\squareside{1};
\def\connectorlength{.15};
\def\springheight{.25};
\def\springlength{.15};
\def\connectorlength{.15};

\fill[blue,opacity=.4] (-\squareside/3,\squareside/2) rectangle (0,-\squareside/2);
\draw (-\squareside/3,\squareside/2) -- (0,\squareside/2)--(0,-\squareside/2)--(-\squareside/3,-\squareside/2);

\draw[line width = 1.5pt] (0,0) -- (\connectorlength,0);
\leftpointycorrectionfactor{\connectorlength};
\spring{\springlength}{\springheight}{\connectorlength};
\rightpointycorrectionfactor{\connectorlength+12*\springlength};
\draw[line width = 1.5pt] (\connectorlength+12*\springlength,0)--(2*\connectorlength+12*\springlength,0);

\draw[fill=blue,opacity=.4,draw opacity=1] (2*\connectorlength+12*\springlength,\squareside/2) 
rectangle (2*\connectorlength+12*\springlength+\squareside,-\squareside/2)node[pos=.5,text opacity=1]{\bf{1}};

\draw[line width = 1.5pt]  (2*\connectorlength+12*\springlength+\squareside,0)--(3*\connectorlength+12*\springlength+\squareside,0);
\leftcurlycorrectionfactor{3*\connectorlength+12*\springlength+\squareside};
\draw[springG] (3*\connectorlength+12*\springlength+\squareside,-.01)--(3*\connectorlength+24*\springlength+\squareside,.0005);
\rightcurlycorrectionfactor{3*\connectorlength+24*\springlength+\squareside};
\draw[line width = 1.5pt] (3*\connectorlength+24*\springlength+\squareside,0)--(4*\connectorlength+24*\springlength+\squareside,0);

\draw[fill=blue,opacity=.4,draw opacity=1]  (4*\connectorlength+24*\springlength+\squareside,\squareside/2)
rectangle (4*\connectorlength+24*\springlength+2*\squareside,-\squareside/2)node[pos=.5,text opacity=1]{\bf{1}};

\draw[line width = 1.5pt] (4*\connectorlength+24*\springlength+2*\squareside,0) -- (5*\connectorlength+24*\springlength+2*\squareside,0);
\leftpointycorrectionfactor{5*\connectorlength+24*\springlength+2*\squareside};
\spring{\springlength}{\springheight}{5*\connectorlength+24*\springlength+2*\squareside};
\rightpointycorrectionfactor{5*\connectorlength+36*\springlength+2*\squareside};
\draw[line width = 1.5pt] (5*\connectorlength+36*\springlength+2*\squareside,0)--(6*\connectorlength+36*\springlength+2*\squareside,0);

\draw[fill=blue,opacity=.4,draw opacity=1] (6*\connectorlength+36*\springlength+2*\squareside,\squareside/2) 
rectangle (6*\connectorlength+36*\springlength+3*\squareside,-\squareside/2)node[pos=.5,text opacity=1]{\bf{1}};

\draw[line width = 1.5pt] (6*\connectorlength+36*\springlength+3*\squareside,0)--(7*\connectorlength+36*\springlength+3*\squareside,0);
\leftcurlycorrectionfactor{7*\connectorlength+36*\springlength+3*\squareside};
\draw[springG] (7*\connectorlength+36*\springlength+3*\squareside,-.01)--(7*\connectorlength+48*\springlength+3*\squareside,.0005);
\rightcurlycorrectionfactor{7*\connectorlength+48*\springlength+3*\squareside};
\draw[line width = 1.5pt] (7*\connectorlength+48*\springlength+3*\squareside,0)--(8*\connectorlength+48*\springlength+3*\squareside,0);

\draw[fill=blue,opacity=.4,draw opacity=1] (8*\connectorlength+48*\springlength+3*\squareside,\squareside/2)
rectangle (8*\connectorlength+48*\springlength+4*\squareside,-\squareside/2)node[pos=.5,text opacity=1]{\bf{1}};

\draw[line width = 1.5pt] (8*\connectorlength+48*\springlength+4*\squareside,0) -- (9*\connectorlength+48*\springlength+4*\squareside,0);
\leftpointycorrectionfactor{9*\connectorlength+48*\springlength+4*\squareside};
\spring{\springlength}{\springheight}{9*\connectorlength+48*\springlength+4*\squareside};
\rightpointycorrectionfactor{9*\connectorlength+60*\springlength+4*\squareside};
\draw[line width = 1.5pt] (9*\connectorlength+60*\springlength+4*\squareside,0)--(10*\connectorlength+60*\springlength+4*\squareside,0);

\fill[blue,opacity=.4] (10*\connectorlength+60*\springlength+4*\squareside,\squareside/2) 
rectangle (10*\connectorlength+60*\springlength+4*\squareside+\squareside/3,-\squareside/2);

\draw (10*\connectorlength+60*\springlength+4*\squareside+\squareside/3,\squareside/2)
--(10*\connectorlength+60*\springlength+4*\squareside,\squareside/2)
--(10*\connectorlength+60*\springlength+4*\squareside,-\squareside/2)
--(10*\connectorlength+60*\springlength+4*\squareside+\squareside/3,-\squareside/2);

\end{tikzpicture}

%% file: fig_spring_dimer_legend.tex
\begin{tikzpicture}
\node[rectangle,draw,thin]{
\begin{tabular}{ll}
\adjustbox{valign=c}{\scalebox{.75}{\begin{tikzpicture}[solid]
\def\springlength{.15};
\def\connectorlength{.15};
\def\springheight{.25};
\def\totalspringlength{12*\springlength};

\draw[line width = 1.5pt] (0,0) -- (\connectorlength,0);
\leftpointycorrectionfactor{\connectorlength};
\spring{\springlength}{\springheight}{\connectorlength};
\draw[line width = 1.5pt]  (\connectorlength + 12*\springlength,0)--(2*\connectorlength+12*\springlength,0);
\rightpointycorrectionfactor{\connectorlength+\totalspringlength};
\end{tikzpicture}}}
&$F_{1}(r) = \kappa_1{r}+\beta_1{r}^2+\O(r^3)$ \\
\\
\adjustbox{valign=c}{\scalebox{.75}{\begin{tikzpicture}[solid]
\def\springlength{.15};
\def\connectorlength{.15};
\def\totalspringlength{12*\springlength};
\def\squareside{1};

\draw[line width = 1.5pt]  (2*\connectorlength+12*\springlength+\squareside,0)--(3*\connectorlength+12*\springlength+\squareside,0);
\leftcurlycorrectionfactor{3*\connectorlength+12*\springlength+\squareside};
\draw[springG] (3*\connectorlength+12*\springlength+\squareside,-.01)--(3*\connectorlength+24*\springlength+\squareside,.0005);
\rightcurlycorrectionfactor{3*\connectorlength+24*\springlength+\squareside};
\draw[line width = 1.5pt] (3*\connectorlength+24*\springlength+\squareside,0)--(4*\connectorlength+24*\springlength+\squareside,0);

\end{tikzpicture}}}
&$F_{2}(r) = \kappa_2{r}+\beta_2{r}^2+\O(r^3)$ 
\end{tabular}};
\end{tikzpicture}

%% file: fig_nanopteron.tex
\begin{tikzpicture}

\draw[thick] (-7.5,0)--(7.5,0)node[below]{$j-ct$};
\draw[thick,<->] (-6.75,-.75)--(-6.75,3)node[right=15pt,below]{$r_j(t)$};
\draw[blue,ultra thick,<->] plot[domain=-7.5:7.5, samples=100] (\x,{2.75/(exp(\x)+exp(-\x))});
\node[above] at (0,1.375){the core};


\draw [line width = .75pt,decorate,decoration={brace,amplitude=5pt}] (-2.5,0)--node[midway,left=2pt]{amplitude $\sim \O(\ep^2)$}(-2.5,1.375);
\draw[densely dotted, line width = .75pt] (-2.5,1.375)--(0,1.375);

\draw[line width =.75pt,decorate,decoration={brace,amplitude=5pt,mirror}] (-4,-.4)--node[midway,below]{wavelength $\sim\O(\ep)$}(4,-.4);
\draw[densely dotted, line width = .75pt] (-4,-.4)--(-4,.0504);
\draw[densely dotted, line width = .75pt] (4,-.4)--(4,.0504);


\draw[densely dotted,line width = .75pt] (6.25,{2/(exp(6.25)+exp(-6.25))}) circle(.5);
\draw[densely dotted, line width = .75pt] (5.75,{2/(exp(5.75)+exp(-5.75))}) -- (5.0253,2.25);
\draw[densely dotted, line width = .75pt] (6.75,{2/(exp(6.75)+exp(-6.75))}) -- (7.4747,2.25);
\draw[densely dotted, line width = .75pt] (6.25,2.5) circle(1.25);

\begin{scope}
\clip (6.25,2.5) circle(1.25);
\node at (6,2.5){\begin{tikzpicture}
\draw[ultra thick,blue] plot[domain=-3:3,samples=100,smooth] (\x,{-.3*cos(10*\x r)});
\end{tikzpicture}};
\end{scope}

\node[above] at (6.25,3.75){the ripple};


\draw [line width = .75pt,decorate,decoration={brace,amplitude=2pt}] (4.75,2.19)--node[midway,left]{amplitude $\sim\O(\ep^{\infty})$}(4.75,2.81);
\draw[line width = .75pt, densely dotted] (4.75,2.19)--(5.36,2.19);
\draw[line width = .75pt, densely dotted] (4.75,2.81)--(5.05,2.81);

\end{tikzpicture}

%% file: section_tw_equations.tex
\section{The Traveling Wave Equations}
\label{tw equations section}

\subsection{The traveling wave ansatz} 
Set
\begin{equation}\label{tw ansatz}
r(j,t) = \begin{cases}
p_1(j-c{t}), &j \text{ is odd} \\
p_2(j-c{t}), &j \text{ is even}.
\end{cases}
\end{equation}
With $S^d$ as the ``shift by $d$'' operator, i.e., 
$$
(S^df)(x) = f(x+d),
$$
we find
$$
\begin{cases}
c^2p_1'' = -2(\kappa{p}_1+\beta{p}_1^2+p_1^3N_1(p_1)) 
+S^1(p_2+p_2^2+p_2^3N_2(p_2)) 
+S^{-1}(p_2+p_2^2+p_2^3N_2(p_2)) \\
\\
c^2p_2'' = -2(p_2+p_2^2+p_2^3N_2(p_2)) 
+S^1(\kappa{p}_1+\beta{p}_1^2+p_1^3N_1(p_1)) 
+S^{-1}(\kappa{p}_1+\beta{p}_1^2+p_1^3N_1(p_1)).
\end{cases}
$$
We rewrite these equations for $p_1$ and $p_2$ in matrix-vector form.  Let 
$$
\pb := \begin{pmatrix*}
p_1 \\
p_2
\end{pmatrix*}
\quadword{and}
N(\pb)
:= \begin{pmatrix*}
N_1(p_1) \\
N_2(p_2)
\end{pmatrix*}.
$$
Then $\pb$ satisfies
\begin{equation}\label{original matrix system compressed}
c^2\pb''+L_{\kappa}\pb + L_{\beta}\pb^{.2} +L_1(\pb^{.3}.N(\pb))
= 0,
\end{equation}
where, given $\alpha \in \R$, we set
\begin{equation}\label{convenient factorization}
L_{\alpha} 
:= \begin{bmatrix*}
2\alpha &-(S^1+S^{-1}) \\
-\alpha(S^1+S^{-1}) &2
\end{bmatrix*}
=
\bunderbrace{
\begin{bmatrix*}
2 &-(S^1+S^{-1}) \\
-(S^1+S^{-1}) &2
\end{bmatrix*}}{L_1}
\bunderbrace{\begin{bmatrix*}
\alpha &0 \\
0 &1
\end{bmatrix*}}{M_{\alpha}}.
\end{equation}
The  ``dot'' notation in \eqref{original matrix system compressed} is componentwise squaring, cubing, or multiplication in the spirit of Matlab, e.g., 
$$
\pb^{.2} := \begin{pmatrix*}
p_1^2 \\
p_2^2
\end{pmatrix*}.
$$


\subsection{Diagonalization} We treat the operators $L_{\kappa}$ and $L_{\beta}$ as Fourier multipliers (cf. Appendix \ref{fourier multipliers appendix}), so that 
$$
\ft[L_{\alpha}\fb](k) = \tL_{\alpha}(k)\ft[\fb](k),
$$
where
$$
\tL_{\alpha}(k)
:= \begin{bmatrix*}
2\alpha &-2\cos(k) \\
-2\alpha\cos(k) &2
\end{bmatrix*}.
$$
We wish to diagonalize $L_{\kappa}$, and we begin by computing the eigenvalues of $\tL_{\kappa}(k)$; they are
\begin{equation}\label{eigenvalues}
\tlambda_{\pm}(k)
:= 1+\kappa \pm \trho(k),
\qquad
\qquad
\trho(k)
:= \sqrt{(1-\kappa)^2+4\kappa\cos^2(k)}.
\end{equation}
What is critical --- and felicitous --- is that the functions $\tlambda_{\pm}$ here are the same as they were in equation (2.5) in \cite{faver-wright} with the parameter $w > 1$ in the mass dimer problem replaced by $\kappa > 1$ in our spring dimer problem.  This allows us to import a tremendous amount of results from that paper with little to no changes.  

To find the corresponding eigenvectors, we first set
$$
\tv_-(k)
:= \frac{2-\tlambda_-(k)}{2\kappa\cos(k)}
\quadword{and}
\tv_+(k)
:= \frac{2\kappa-\tlambda_+(k)}{2\cos(k)}.
$$
Of course $\tv_{\pm}$ are not defined at numbers $k$ such that $\cos(k) = 0$, but we will show that $\tv_{\pm}$ have removable singularities there.  Then
$$
\tvb_+(k)
:= 
\begin{pmatrix*}
1 \\
\tv_+(k)
\end{pmatrix*}
\quadword{and}
\tvb_-(k)
:=
\begin{pmatrix*}
\tv_-(k) \\
1
\end{pmatrix*}
$$
are the eigenvectors corresponding to $\tlambda_{\pm}(k)$.  With
$$
\tJ(k) 
:= \begin{bmatrix*}
\tvb_-(k) &\tvb_+(k)
\end{bmatrix*},
\qquad
\tJ_1(k)
:= \tJ(k)^{-1},
\quadword{and}
\tLambda(k) 
:= \begin{bmatrix*}
\tlambda_-(k) &0 \\
0 &\tlambda_+(k)
\end{bmatrix*},
$$
we have 
$$
\tL_{\kappa}(k) 
= \tJ(k)\tLambda(k)\tJ_1(k).
$$
Taking $J$, $J_1$, and $\Lambda$ to be the Fourier multipliers with the symbols $\tJ$, $\tJ_1$, and $\tLambda$, we diagonalize the operator $L_{\kappa}$:
$$
L_{\kappa} = J\Lambda{J}_1.
$$
Before we exploit this diagonalization, we summarize the essential properties of these various Fourier multipliers.  
%

\begin{proposition}\label{eigenprops}
\begin{enumerate}[label={\bf(\roman*)}]
\item $\tlambda_-(0) = 0$ and $\tlambda_+(0) = 2+2\kappa$.
\item For $q > 0$, let 
$$
\Sigma_{q} = \set{z \in \mathbb{C}}{|\Im(z)| < q}
$$
and let $\overline{\Sigma}_{q}$ be its closure.  There exists $q_0 > 0$ such that $\tlambda_{\pm}$ and $\tv_{\pm}$ extend to even, $\pi$-periodic, bounded, complex-valued analytic functions on $\overline{\Sigma}_{q_0}$.
\item For all $k \in \R$,
\begin{equation}\label{CRZ bifurcation condition}
|\tlambda_{\pm}'(k)| \le 2
\quadword{and}
|\tlambda_{\pm}'(k)| \le 2c_{\kappa}|k|,
\end{equation}
where
$$
c_{\kappa} := \sqrt{\frac{1}{2}\tlambda_-''(0)}
= \sqrt{\frac{2\kappa}{1+\kappa}} > 1.
$$
\item If $c^2 > c_{\kappa}^2$, then
$$
c^2k^2 - \tlambda_-(k) > 0
$$
for all $k \in \R$.  
\item There exists $c_- \in (0,1)$ such that for all $c \ge c_-$ there is a unique $\Omega_c > 0$ satisfying
\begin{equation}
c^2\Omega_c^2-\tlambda_+(\Omega_c) = 0.
\end{equation}
This number $\Omega_c$ also satisfies
\begin{equation}\label{omegac intermediate bound}
\frac{\sqrt{2\kappa}}{c} \le \Omega_c \le \frac{\sqrt{2+2\kappa}}{c},
\end{equation}
and there is a constant $b_0 > 0$ such that 
\begin{equation}\label{omegac derivative bound}
|2c^2\Omega_c - \tlambda_+'(\Omega_c)| \ge b_0.
\end{equation}
\end{enumerate}
\end{proposition}

\begin{proof}
Since our eigenvalues $\tlambda_{\pm}$ are the same as those in the mass-dimer problem, with $w$ replaced by $\kappa$, all of these results were proved in Lemma 2.1 in \cite{faver-wright}, with the exception of (ii) for the functions $\tv_{\pm}$.  This is straightforward, as the proof of that lemma gives $q_0 > 0$ such that $\Re[(1-\kappa)^2+4\kappa\cos^2(z)] > 0$ for $|\Im(z)| \le q_0$.  We then use the principal square root to extend $\trho(k)$, and thereby $\tlambda_{\pm}(k)$, into analytic maps from $\overline{\Sigma}_{q_0}$ to $\mathbb{C}$.  

For $z \in \mathbb{C}$, set $f(z) = 2-\tlambda_-(z)$ and $z_k = (2k+1)/2\pi$.  Then $f(z_k) = 0$, so $f(z) = (z-z_k)f_k(z)$ for some analytic function $f_k$.  Likewise, since $\cos(\cdot)$ has simple zeros at $z_k$, we can write $\cos(z) = (z-z_k)g_k(z)$, where $g_k(z_k) \ne 0$.  Hence
$$
\tv_-(z) = \frac{2-\tlambda_-(z)}{2\kappa\cos(z)} = \frac{(z-z_k)f_k(z)}{2\kappa(z-z_k)g_k(z)} = \frac{f_k(z)}{2\kappa{g}_k(z)}.
$$
That is, $\tv_-$ has a removable singularity at each $z_k$, and so $\tv_-$ is analytic on $\overline{\Sigma}_{q_0}$.  The same argument shows that $\tv_+$ is analytic on $\overline{\Sigma}_{q_0}$.  
\end{proof}

Now, using the factorization $L_{\kappa} = J\Lambda{J}_1$, we diagonalize the system \eqref{original matrix system compressed}.  Set
\begin{equation}\label{diag cov}
\hb = J_1\pb.
\end{equation}
Then \eqref{original matrix system compressed} is equivalent to
\begin{equation}\label{tw post diag}
c^2\partial_x^2\hb + \Lambda\hb + J_1L_{\beta}[(J\hb)^{.2}] +J_1L_1[(J\hb)^{.3}.N(J\hb)]
= 0.
\end{equation}
Using the factorization \eqref{convenient factorization} of the operators $L_{\kappa}$ and $L_{\beta}$, we have
$$
J_1L_{\beta}
= J_1L_1M_{\beta} 
= J_1L_1M_{\kappa}M_{\beta/\kappa} \\
= J_1L_{\kappa}M_{\beta/\kappa}
= J_1(J\Lambda{J}_1)M_{\beta/\kappa}
= \Lambda{J}_1M_{\beta/\kappa}.
$$
In a similar way, we obtain
$$
J_1L_1 
= \Lambda{J}_1M_{1/\kappa}.
$$
The advantage is that the quadratic and cubic factors in \eqref{tw post diag} now each have the same prefactor of $\Lambda$, so that \eqref{tw post diag} becomes
\begin{equation}\label{tw post diag and factoring}
(c^2\partial_x^2+\Lambda)\hb + \Lambda{J}_1M_{\beta/\kappa}[(J\hb)^{.2}]+ \Lambda{J}_1M_{1/\kappa}[(J\hb)^{.3}.N(J\hb)]
= 0
\end{equation}
With
\begin{equation}\label{B defn}
B(\hb,\grave{\hb})
= \begin{pmatrix*}
B_1(\hb,\grave{\hb}) \\
B_2(\hb,\grave{\hb})
\end{pmatrix*}
:= J_1M_{\beta/\kappa}[(J\hb).(J\grave{\hb})],
\end{equation}
\begin{equation}\label{N defn}
\nl(\hb)
:= \hb.N(\hb),
\end{equation}
and
\begin{equation}\label{Q defn}
\cb(\hb,\grave{\hb},\breve{\hb})
= \begin{pmatrix*}
\cb_1(\hb,\grave{\hb},\breve{\hb}) \\
\cb_2(\hb,\grave{\hb},\breve{\hb})
\end{pmatrix*}
:= J_1M_{1/\kappa}[(J\hb).(J\grave{\hb}).\nl(J\breve{\hb})]
\end{equation}
we further compress \eqref{tw post diag and factoring} to
\begin{equation}\label{tw post diag, factoring, and labeling}
(c^2\partial_x^2+\Lambda)\hb
+ \Lambda{B}(\hb,\hb)
+\Lambda\cb(\hb,\hb,\hb)
=0.
\end{equation}

\subsection{The Friesecke-Pego cancelation}\label{fp cancelation section}
The first component of \eqref{tw post diag, factoring, and labeling} is
\begin{equation}\label{first comp spatial}
(c^2\partial_x^2+\lambda_-)h_1 + \lambda_-B_1(\hb,\hb) + \lambda_-\cb_1(\hb,\hb,\hb) =0.
\end{equation}
Applying the Fourier transform, this becomes
\begin{equation}\label{second comp}
-(c^2k^2-\tlambda_-(k))\hat{h}_1(k) + \tlambda_-(k)\ft[B_1(\hb,\hb)](k) + \tlambda_-(k)\ft[\cb_1(\hb,\hb,\hb)](k) = 0,
\end{equation}
By Proposition \ref{eigenprops}, we have $c^2k^2-\tlambda_-(k) > 0$ for all $k \ne 0$, as long as $c^2 > c_{\kappa}^2$.  We assume this lower bound on $c^2$ from now on, so that \eqref{second comp} becomes
$$
\hat{h}_1(k) + \tvarpi_c(k)\ft[B_1(\hb,\hb)](k) + \tvarpi_c(k)\ft[\cb_1(\hb,\hb,\hb)](k) = 0,
$$
where
$$
\tvarpi_c(k) 
:= -\frac{\tlambda_-(k)}{c^2k^2-\tlambda_-(k)}.
$$
It is easy to calculate that $\tvarpi_c$ has a removable singularity at $k=0$, so $\tvarpi_c$ extends to an even, $\pi$-periodic, bounded complex-valued analytic function on the strip $\overline{\Sigma}_{q_0}$ from Proposition \ref{eigenprops}.  

Let $\varpi_c$ be the Fourier multiplier with symbol $\tvarpi_c$.  Then any function $\hb = (h_1,h_2)$ that solves
$$
h_1 + \varpi_cB_1(\hb,\hb) + \varpi_c\cb_1(\hb,\hb,\hb) = 0,
$$
will solve \eqref{first comp spatial}, and so we can find solutions to the entire system \eqref{tw post diag, factoring, and labeling} by studying
$$
\H_c(\hb)
:=
\begin{bmatrix*}
1 &0 \\
0 &c^2\partial_x^2+\lambda_+
\end{bmatrix*}\hb
+
\begin{bmatrix*}
\varpi_c &0 \\
0 &\lambda_+
\end{bmatrix*}
B(\hb,\hb)
+
\begin{bmatrix*}
\varpi_c &0 \\
0 &\lambda_+
\end{bmatrix*}
\cb(\hb,\hb,\hb)
=0.
$$


\begin{lemma}\label{even even symmetry lemma}
If $h_1$ and $h_2$ are both even and $\hb = (h_1,h_2)$, then the components of $\H_c(\hb)$ are also both even.
\end{lemma}

\begin{proof}
We use the fact that if the symbol of a Fourier multiplier $\mu$ is even and $f$ is an even function, then $\mu{f}$ is even. Observe that the Fourier multipliers in the definition of $\H_c$ --- which are $c^2\partial_x^2+\lambda_+$, $\lambda_+$, $\varpi_c$, $J_1$, and $J$ --- all have even symbols.  Moreover, multiplication and composition of even functions of course preserves evenness.  Together with the convenient structure of $\H_c$, this proves the lemma.
\end{proof}


\subsection{The long wave scaling}\label{lw scaling section}
This is our final change of variables.  We set 
\begin{equation}\label{lw scaling}
\hb(x) = \ep^2\thetab(\ep{x}),
\end{equation}
where $\thetab(X) = (\theta_1(X),\theta_2(X))$, and we take the wave speed $c$ to satisfy
$$
c^2 = \cep^2 := c_{\kappa}^2 + \ep^2.
$$
That is, we intend to solve 
\begin{equation}\label{Thetab ep problem}
\H_{\cep}(\ep^2\thetab(\ep\cdot)) = 0
\end{equation}
for $\thetab$ with $\ep$ small. 


Let $\varpi^{\ep}$ be the Fourier multiplier with symbol
\begin{equation}\label{tvarpi ep defn}
\tilde{\varpi^{\ep}}(K) 
:=  \ep^2\tvarpi_{\cep}(\ep{k})
= -\frac{\ep^2\tlambda_-(\ep{k})}{\cep^2(\ep{k})^2-\tlambda_-(\ep{k})}.
\end{equation}
Per the convention outlined in Appendix \ref{fourier multipliers appendix}, for any other Fourier multiplier $\mu$ in the definition of $\H_c$, let $\mu^{\ep}$ have the symbol $\tilde{\mu^{\ep}}(k) = \tmu(\ep{k})$, where of course $\tmu$ is the symbol of $\mu$.  Let
\begin{equation}\label{B ep defn}
B^{\ep}(\thetab,\grave{\thetab}) = 
\begin{pmatrix*}
B_1^{\ep}(\thetab,\grave{\thetab}) \\
B_2^{\ep}(\thetab,\grave{\thetab})
\end{pmatrix*}
:= J_1^{\ep}M_{\beta/\kappa}[(J^{\ep}\thetab).(J^{\ep}\grave{\thetab})]
\end{equation}
and
\begin{equation}\label{Q ep defn}
\cb^{\ep}(\thetab,\grave{\thetab},\breve{\thetab}) 
=
\begin{pmatrix*}
\cb_1^{\ep}(\thetab,\grave{\thetab},\breve{\thetab}) \\
\cb_2^{\ep}(\thetab,\grave{\thetab},\breve{\thetab})
\end{pmatrix*}
:= J_1^{\ep}M_{\beta/\kappa}\left((J^{\ep}\thetab).(J^{\ep}\grave{\thetab}).\nl(\ep^2J^{\ep}\breve{\thetab})\right).
\end{equation}
Then by the scaling properties of Fourier multipliers, our problem \eqref{Thetab ep problem} is equivalent to 
\begin{equation}\label{lw eqns}
\Thetab_{\ep}(\thetab)
:= \D_1^{\ep}\thetab + \D_2^{\ep}B^{\ep}(\thetab,\thetab) + \D_2^{\ep}\cb^{\ep}(\thetab,\thetab,\thetab)
= 0,
\end{equation}
where
$$
\D_1^{\ep}
:= \begin{bmatrix*}
1 &0 \\
0 &(c_{\kappa}^2+\ep^2)\ep^2\partial_X^2 + \lambda_+^{\ep}
\end{bmatrix*}
\quadword{and}
\D_2^{\ep}
:= \begin{bmatrix*}
\varpi^{\ep}&0 \\
0 &\ep^2\lambda_+^{\ep}
\end{bmatrix*}.
$$
Note that because the small parameter $\ep^2$ multiplies the second derivative operator $\partial_X^2$ in $\D_1^{\ep}$, our problem \eqref{Thetab ep problem} is singularly perturbed just as in the mass dimer problem.  This was not a feature of the monatomic problem in \cite{friesecke-pego1}, as the long wave problem there involved only an equation analogous to \eqref{first comp spatial}; the subsequent factoring and cancelation in the monatomic equation, which proceeded like ours in Section \ref{fp cancelation section}, removed that singularity.  We also remark that the nonlinear operator $\Thetab_{\ep}$ --- or, more precisely, the three terms in the definition of $\Thetab_{\ep}$ in \eqref{lw eqns} --- inherits the ``even $\times$ even'' symmetry of $\H_{\cep}$ from Lemma \ref{even even symmetry lemma}.

\subsection{The formal long wave limit} In this section we formally define what $\Thetab_0$ should be by assigning meaning to $\D_1^0$, $\D_2^0$, $B^0$, and $\cb^0$ in a natural way.  First, Proposition \ref{eigenprops} gives
$$
\tlambda_+(0) = 2+2\kappa
\quadword{and} 
\tlambda_-(0) = 0,
$$
so that 
$$
\tJ(0) 
= \begin{bmatrix*}[r]
1/\kappa &1 \\
1 &-1
\end{bmatrix*}
\quadword{and}
\tJ_1(0) 
= \frac{\kappa}{\kappa+1}\begin{bmatrix*}[r]
1 &1 \\
1 &-1/\kappa
\end{bmatrix*}.
$$
This motivates the definition of the Fourier multipliers $J^0$ and $J_1^0$ as 
\begin{equation}\label{J0}
J^0 
:= \begin{bmatrix*}[r]
1/\kappa &1 \\
1 &-1
\end{bmatrix*}
\quadword{and}
J_1^0 
:= \frac{\kappa}{\kappa+1}\begin{bmatrix*}[r]
1 &1 \\
1 &-1/\kappa
\end{bmatrix*}.
\end{equation}
and from this and \eqref{B ep defn} we have
\begin{equation}\label{B0}
B^0(\thetab,\grave{\thetab}):=
\frac{\kappa}{\kappa+1}\begin{pmatrix*}
\beta/\kappa(\theta_1/\kappa+\theta_2)(\grave{\theta}_1/\kappa+\grave{\theta}_2)+(\theta_1-\theta_2)(\grave{\theta}_1-\grave{\theta}_2) \\
\beta/\kappa(\theta_1/\kappa+\theta_2)(\grave{\theta}_1/\kappa+\grave{\theta}_2)-(\theta_1-\theta_2)(\grave{\theta}_1-\grave{\theta}_2)/\kappa
\end{pmatrix*}.
\end{equation}

We could define $\cb^0$ in the same way using $J^0$ and $J_1^0$, but it is straightforward to see that $\cb^0$ will be identically zero thanks to the extra factor of $\ep^2$ that $\cb^{\ep}$ carries within $\nl$, per \eqref{Q ep defn} and \eqref{N defn}. This factor of $\ep^2$ will resurface frequently in the depths of the estimates to come.

Using the same calculation that we did in \cite{faver-wright}, we compute the Taylor expansion of $\tlambda_+$ and find the natural definition of $\tvarpi^0$ to be
$$
\tvarpi^0(k) = -\frac{c_{\kappa}^2}{1+\alpha_{\kappa}k^2},
$$
where
$$
\alpha_{\kappa} = \frac{c_{\kappa}^2}{3}\frac{1-\kappa+\kappa^2}{(1+\kappa)^2}.
$$
So, we set
$$
\varpi^0 := -c_{\kappa}^2(1-\alpha_{\kappa}\partial_X^2)^{-1}.
$$

All together, this implies (formally) that 
\begin{equation}\label{Theta0}
\Thetab_0(\thetab)
= \begin{bmatrix*}
1 &0 \\
0 &2+2\kappa
\end{bmatrix*} 
\thetab
+ \begin{bmatrix*}
\varpi^0 &0 \\
0 &0
\end{bmatrix*}
B^0(\thetab,\thetab)
\end{equation}
Now we consider the problem $\Thetab_0(\thetab) = 0$ with $\thetab = (\theta_1,\theta_2)$. We find from the second component of \eqref{Theta0} that
$$
(2+2\kappa)\theta_2 = 0,
$$
so $\theta_2 = 0$.  Then the first component reduces to
$$
0 = 
\theta_1 + \varpi^0\left(\frac{\kappa}{\kappa+1}\left(\frac{\beta}{\kappa^3}\theta_1^2+\theta_1^2\right)\right) 
=
\theta_1 - c_{\kappa}^2(1-\alpha_{\kappa}\partial_X^2)^{-1}\left(\frac{\kappa}{\kappa+1}\left(\frac{\beta}{\kappa^3}+1\right)\theta_1^2\right).
$$
Applying $1-\alpha_{\kappa}\partial_X^2$ to both sides, we get
\begin{equation}\label{KdV ode}
\alpha_{\kappa}\theta_2'' - \theta_2 + c_{\kappa}^2\frac{\kappa}{\kappa+1}\left(\frac{\beta}{\kappa^3}+1\right)\theta_2^2 = 0.
\end{equation}
This is a rescaling of the ordinary differential equation that gives the $\sech^2$-type traveling wave profiles for the KdV equation, provided that 
\begin{equation}\label{beta kappa rel}
\frac{\beta}{\kappa^3} +1 \ne 0
\end{equation}
to keep the nonlinear term $\theta_2^2$ present. If we require $\beta$ and $\kappa$ to satisfy \eqref{beta kappa rel}, then the solution to \eqref{KdV ode} is
\begin{equation}\label{sigma soliton defn}
\theta_2(X) 
= \sigma(X) 
:= \frac{3}{2c_{\kappa}^2}\left(\frac{\kappa+1}{\kappa}\right)\left(\frac{\beta}{\kappa^3}+1\right)^{-1}\sech^2\left(\frac{X}{2\sqrt{\alpha_w}}\right).
\end{equation}


%% file: section_periodic_solutions.tex
\section{Periodic Solutions}
\label{periodic solutions section}

From \eqref{lw eqns}, the linearization of $\Thetab_{\ep}$ at $\thetab = 0$ is the operator
$$
\D_1^{\ep}
=
\begin{bmatrix*}
1 &0 \\
0 &\cep^2\ep^2\partial_X^2+\lambda_+^{\ep}
\end{bmatrix*}.
$$
Using the scaling properties of Fourier multipliers, we see that if $\D_1^{\ep}\thetab = 0$ for some $2\pi$-periodic function $\thetab = (\theta_1,\theta_2)$, then $\theta_2 = 0$ and 
$$
\big(\cep^2(\ep{k})^2-\tlambda_+(\ep{k})\big)\hat{\theta}_2(k) 
= 0
$$
for all $k$.  From Theorem \ref{eigenprops}, let $\Omega_{\cep}$ be the unique positive number such that 
\begin{equation}\label{Omega ep zero defn}
\cep^2\Omega_{\cep}^2-\tlambda_+(\Omega_{\cep}) 
= 0,
\end{equation}
and take 
\begin{equation}\label{omega ep defn}
\omega_{\ep} 
:= \frac{\Omega_{\cep}}{\ep}
\end{equation}
so that
\begin{equation}\label{omega ep zero prop}
\cep^2(\ep\omega_{\ep})^2 - \tlambda_+(\ep\omega_{\ep}) = 0.
\end{equation}
It follows that $\D_1^{\ep}[\cos(\omega_{\ep}\cdot)\jb] = 0$, and, if we restrict $\D_1^{\ep}$ to functions $\thetab$ with $\theta_1$ and $\theta_2$ both even, then, per Lemma \ref{even even symmetry lemma}, one can show that the kernel of $\D_1^{\ep}$ is in fact spanned by $\cos(\omega_{\ep}\cdot)\jb$.  As in \cite{faver-wright}, the search for periodic solutions of \eqref{lw eqns} then fits naturally into the set-up of ``bifurcation from a simple eigenvalue'' from Crandall and Rabinowitz \cite{crandall-rabinowitz} and Zeidler \cite{zeidler}, and the critical bifurcation parameter is $\omega_{\ep}$. Our result, proved in Appendix \ref{periodic solutions appendix}, is the following theorem.

\begin{theorem}\label{periodic main theorem}
There exist $\ep_{\per}$, $a_{\per} > 0$ such that for all $\ep \in (0,\ep_{\per})$, there are maps
\begin{align*}
[-a_{\per},a_{\per}] &\to \R \colon a \mapsto \omega_{\ep}^a \\
[-a_{\per},a_{\per}] &\to \C_{\per}^{\infty} \cap \{even \ functions\} \colon a \mapsto \psi_{\ep,1}^a \\
[-a_{\per},a_{\per}] &\to \C_{\per}^{\infty} \cap \{even \ functions\} \colon a \mapsto \psi_{\ep,2}^a
\end{align*}
such that the following hold. 

\begin{enumerate}[label={\bf(\roman*)}]
\item If
$$
\nub := \cos(\cdot)\jb,
\qquad
\psib_{\ep}^a := \begin{pmatrix*}
\psi_{\ep,1}^a \\
\psi_{\ep,2}^a
\end{pmatrix*},
\quadword{and}
\varphib_{\ep}^a(X) := \nub(\omega_{\ep}^aX) + \psib_{\ep}^a(\omega_{\ep}^aX),
$$
then $\thetab := a\varphib_{\ep}^a$ solves \eqref{lw eqns} for all $|a| \le a_{\per}$ and $\ep \in (0,\ep_{\per})$.

\item The frequency $\omega_{\ep}^0$ satisfies $\omega_{\ep}^0 = \omega_{\ep}$ as defined in \eqref{omega ep defn} above.  We say that $\omega_{\ep} = \O(1/\ep)$ in the sense that there are constants $C_1$, $C_2 > 0$ such that 
$$
\frac{C_1}{\ep} < \omega_{\ep} < \frac{C_2}{\ep}
$$
for all $\ep \in (0,\ep_{\per})$.

\item $\psi_{\ep,1}^0 = \psi_{\ep,2}^0 = 0$.

\item For all $r \ge 0$, there is $C_r > 0$ such that
$$
|\ep\omega_{\ep}^a| + \norm{\psib_{\ep}^a}_{\C_{\per}^r \times \C_{\per}^r} \le C_r
$$
and
$$
|\omega_{\ep}^a - \omega_{\ep}^{\grave{a}}| + \norm{\psib_{\ep}^a - \psib_{\ep}^{\grave{a}}}_{\C_{\per}^r \times \C_{\per}^r} \le C_r|a-\grave{a}|
$$
for all $|a|$, $|\grave{a}| \le a_{\per}$ and $0 < \ep < \ep_{\per}$.
\end{enumerate}
\end{theorem}


%% file: section_nanopteron_equations.tex
\section{The Nanopteron Equations}
\label{nanopteron equations section}

\subsection{Beale's ansatz} We recall from Appendix \ref{sobolev spaces appendix} that $E_q^1$ is the space of even, exponentially decaying functions in $H^1$:
$$
E_q^1 
:= \set{f \in H^1}{f \text{ is even and} \cosh(q\cdot)f \in H^1}
$$
and that 
$$
\norm{f}_{1,q} := \norm{\cosh(q\cdot)f}_{H^1}.
$$
We return to our main problem $\Thetab_{\ep}(\thetab) = 0$ from \eqref{lw eqns} and make Beale's ansatz:
\begin{equation}\label{beale}
\thetab 
= \an_{\ep}(\etab,a) 
:= \sigmab + a\varphib_{\ep}^a + \etab,
\end{equation}
where 
\begin{itemize}[leftmargin=*]
\item $\etab = (\eta_1,\eta_2) \in E_q^1 \times E_q^1$;
\item $a \in \R$;
\item $\varphib_{\ep}^a$ is periodic and $a\varphib_{\ep}^a$ satisfies $\Thetab_{\ep}(a\varphib_{\ep}^a) = 0$, per Theorem \ref{periodic main theorem};
\item $\sigmab := (\sigma,0)$, where $\sigma$ solves the KdV profile equation \eqref{KdV ode}.
\end{itemize}
Beale's ansatz introduces three unknowns into our problem: the amplitude $a$ of the periodic ripple and the decaying terms $\eta_1$ and $\eta_2$.  

We find that $\Thetab_{\ep}(\an_{\ep}(\etab,a)) = 0$ is componentwise equivalent to
\begin{equation}\label{beale first pass}
\begin{cases}
\eta_1 = -\sum_{k=1}^5 \left(\j_{k1}^{\ep}(\etab,a) + \j_{k2}^{\ep}(\etab,a)\right) - \j_6^{\ep}(\etab,a) =: \nanorhs_1^{\ep}(\etab,a) \\
(\cep^2\ep^2\partial_X^2 + \lambda_+^{\ep})\eta_2 = -\sum_{k=1}^5 \left(\l_{k1}^{\ep}(\etab,a) + \l_{k2}^{\ep}(\etab,a)\right) - \l_6^{\ep}(\etab,a) =: \nanorhs_2^{\ep}(\etab,a).
\end{cases}
\end{equation}
where we have used the bilinearity of $B^{\ep}$ and of $\cb^{\ep}(\cdot,\cdot,\an_{\ep}(\etab,a))$ in its first two arguments to break up the terms as 
\begin{align*}
\j_{11}^{\ep}(\etab,a) &:= \sigma + \varpi^{\ep}B_1^{\ep}(\sigmab,\sigmab) &\j_{12}^{\ep}(\etab,a) &:= \varpi^{\ep}\cb_2^{\ep}(\sigmab,\sigmab,\an_{\ep}(\etab,a)) \\
\j_{21}^{\ep}(\etab,a) &:= 2\varpi^{\ep}B_1^{\ep}(\sigmab,\etab) &\j_{22}^{\ep}(\etab,a) &:=\varpi^{\ep}\cb_1^{\ep}(\sigmab,\etab,\an_{\ep}(\etab,a)) \\
\j_{31}^{\ep}(\etab,a) &:= 2a\varpi^{\ep}B_1^{\ep}(\sigmab,\varphib_{\ep}^a) &\j_{32}^{\ep}(\etab,a) &:= 2a\varpi^{\ep}\cb_1^{\ep}(\sigmab,\varphib_{\ep}^a,\an_{\ep}(\etab,a)) \\
\j_{41}^{\ep}(\etab,a) &:= 2a\varpi^{\ep}B_1(\etab,\varphib_{\ep}^a) &\j_{42}^{\ep}(\etab,a) &=2a\varpi^{\ep}\cb_1^{\ep}(\etab,\varphib_{\ep}^a,\an_{\ep}(\etab,a)) \\
\j_{51}^{\ep}(\etab,a) &:= 2\varpi^{\ep}B_1^{\ep}(\etab,\etab) &\j_{52}^{\ep}(\etab,a) &:= 2\varpi^{\ep}\cb_1^{\ep}(\etab,\etab,\an_{\ep}(\etab,a)) \\
\\
\l_{11}^{\ep}(\etab,a) &:= \lambda_+^{\ep}B_2^{\ep}(\sigmab,\sigmab) &\l_{12}^{\ep}(\etab,a) &:= \lambda_+^{\ep}\cb_2^{\ep}(\sigmab,\sigmab,\an_{\ep}(\etab,a)) \\
\l_{21}^{\ep}(\etab,a) &:= 2\lambda_+^{\ep}B_2^{\ep}(\sigmab,\etab) &\l_{22}^{\ep}(\etab,a) &:= 2\lambda_+^{\ep}\cb_2^{\ep}(\sigmab,\etab,\an_{\ep}(\etab,a)) \\
\l_{31}^{\ep}(\etab,a) &:= 2a\lambda_+^{\ep}B_2^{\ep}(\sigmab,\varphib_{\ep}^a) & \l_{32}^{\ep}(\etab,a) &:= 2a\lambda_+^{\ep}\cb_2^{\ep}(\sigmab,\varphib_{\ep}^a,\an_{\ep}(\etab,a)) \\
\l_{41}^{\ep}(\etab,a) &:= 2a\lambda_+^{\ep}B_2^{\ep}(\etab,\varphib_{\ep}^a) &\l_{42}^{\ep}(\etab,a) &:= 2a\lambda_+^{\ep}\cb_2^{\ep}(\etab,\varphib_{\ep}^a,\an_{\ep}(\etab,a)) \\
\l_{51}^{\ep}(\etab,a) &:= \lambda_+^{\ep}B_2^{\ep}(\etab,\etab) &\l_{52}^{\ep}(\etab,a) &:= \lambda_+^{\ep}\cb_2^{\ep}(\etab,\etab,\an_{\ep}(\etab,a)).
\end{align*}
\begin{align*}
\j_{6}^{\ep}(\etab,a) &:= a^2\varpi^{\ep}\left[\cb_1^{\ep}(\varphib_{\ep}^a,\varphib_{\ep}^a,\an_{\ep}(\etab,a))-\cb_1^{\ep}(\varphib_{\ep}^a,\varphib_{\ep}^a,a\varphib_{\ep}^a)\right] \\
\l_{6}^{\ep}(\etab,a) &:= a^2\lambda_+^{\ep}\left[\cb_2^{\ep}(\varphib_{\ep}^a,\varphib_{\ep}^a,\an_{\ep}(\etab,a))-\cb_2^{\ep}(\varphib_{\ep}^a,\varphib_{\ep}^a,a\varphib_{\ep}^a)\right].
\end{align*}

\subsection{Adjustments to the nanopteron equations}
We need to modify this system in several ways before it becomes amenable to our quantitative contraction mapping argument.  First, as it stands, the term $\j_{21}^{\ep}$ is $\O(1)$ in $\ep$, which is inadequate for our intended methods.  So, we add the term $2\varpi^0B_1^0(\sigmab,\etab)$ to both sides of the equation for $\eta_1$.  
Expanding $B_1^0(\sigmab,\etab)$ from its definition in \eqref{B0}, we have
\begin{equation}\label{death star defns}
2\varpi^0B_1^0(\sigmab,\etab)
= \bunderbrace{\frac{2\kappa}{\kappa+1}\left(\frac{\beta}{\kappa^3}+1\right)\varpi^0(\sigma\eta_1)}{-\K_1\eta_1}
+ \bunderbrace{\frac{2\kappa}{\kappa+1}\left(\frac{\beta}{\kappa^2}-1\right)\varpi^0(\sigma\eta_2)}{\K_2\eta_2}.
\end{equation}
Then subtracting $\K_2\eta_2$ from both sides, we find
\begin{equation}\label{A defn}
\bunderbrace{\eta_1-\K_2\eta_1}{\A\eta_1} = \nanorhs_1^{\ep,\mod}(\etab,a)- \K_2\eta_2,
\end{equation}
where 
$$
\nanorhs_1^{\ep,\mod}(\etab,a)
:= -\sum_{k=1}^5 \left(\j_{k1}^{\ep,\mod}(\etab,a)+\j_{k2}^{\ep}(\etab,a)\right)-\j_6^{\ep}(\etab,a)
$$
and
$$
\j_{1k}^{\ep,\mod}(\etab,a)
:= \begin{cases}
2\varpi^{\ep}B_1^{\ep}(\sigmab,\etab)-2\varpi^0B_1^0(\sigmab,\etab), k = 2 \\
\\
\j_{1k}^{\ep}(\etab,a), k \ne 2.
\end{cases}
$$
The operator $\A$ given in \eqref{A defn}
is invertible on $E_q^1$ by Proposition \ref{faver wright nanopteron} for $q$ sufficiently small, and so we may solve for $\eta_1$:
$$
\eta_1 = \A^{-1}\nanorhs_1^{\ep,\mod}(\etab,a)-\A^{-1}\K_2\eta_2.
$$
The term $\A^{-1}\K_2\eta_2$ is still $\O(1)$ in $\ep$, which will ruin our contraction estimates.  However, once we establish our fixed point equation for $\eta_2$, we will rewrite the system yet again in a manner that eliminates this difficulty.  

%
%
In the equation for $\eta_2$ in \eqref{beale first pass}, which is
\begin{equation}\label{T ep defn}
\bunderbrace{(\cep^2\ep^2\partial_X^2 + \lambda_+^{\ep})\eta_2}{\T_{\ep}\eta_2} = \nanorhs_2^{\ep}(\etab,a),
\end{equation}
the operator $\T_{\ep}$ is not invertible since \eqref{omega ep zero prop} implies
$$
\hat{\T_{\ep}f}(\pm\omega_{\ep}) 
= 0
$$
for any function $f$.  That is, $\T_{\ep}$ is not surjective.  Equivalently, if $g$ is an even function in the range of $\T_{\ep}$, then  
\begin{equation}\label{iota ep defn}
\bunderbrace{\int_{\R} g(X)\cos(\omega_{\ep}X)\dX}{\iota_{\ep}[g]}
= \hat{g}(\pm\omega_{\ep})
= 0.
\end{equation}
So, if $\Thetab_{\ep}(\an_{\ep}(\etab,a))=0$, then $(\etab,a)$ must meet the ``solvability condition''
$$
\iota_{\ep}[\nanorhs_2^{\ep}(\etab,a)] 
= 0.
$$
Since we cannot merely invert $\T_{\ep}$ to solve for $\eta_2$ above, we instead follow the route established and motivated in \cite{faver-wright}.   Let
\begin{equation}\label{upsilon ep defn}
\begin{aligned}
\nub_{\ep} 
&:= \varphib_{\ep}^0 := \cos(\omega_{\ep}\cdot)\jb \\
\\
\chi_{\ep} 
&:= \lambda_+^{\ep}J_1^{\ep}[(J_2^0\sigmab).(J_2^{\ep}\nub_{\ep})]\cdot\jb \\
\\
\upsilon_{\ep} 
&:= \iota_{\ep}[\chi_{\ep}].
\end{aligned}
\end{equation}
Subsequent estimates, which we detail in Proposition \ref{faver wright nanopteron}, reveal that $\upsilon_{\ep}$ is bounded away from zero and that $\T_{\ep}$ is invertible on $E_q^1 \cap \ker(\iota_{\ep})$ for $q$ sufficiently small, and so we may define, for any $f \in E_q^1$, 
\begin{equation}\label{Pep defn orig}
\P_{\ep}f 
:= \T_{\ep}^{-1}\left(f-\frac{1}{\upsilon_{\ep}}\iota_{\ep}[f]\chi_{\ep}\right).
\end{equation}
Then with
$$
\l_{31}^{\ep,\mod}(\etab,a) 
:= \begin{cases}
-2a\chi_{\ep} + \l_{31}^{\ep}(\etab,a), k = 3 \\
\\
\l_{k1}^{\ep}(\etab,a), \ k \ne 3, 
\end{cases} 
$$
and
$$
\nanorhs_2^{\ep,\mod}(\etab,a) 
:= \sum_{k= 1}^5 \l_{k1}^{\ep,\mod}(\etab,a)+\sum_{k=1}^5 \l_{k2}^{\ep}(\etab,a) + \l_6^{\ep}(\etab,a),
$$
the pair of equations
$$
\begin{cases}
\T_{\ep}\eta_2
= \nanorhs_2^{\ep}(\etab,a) \\
\\
\iota_{\ep}[\nanorhs_2^{\ep}(\etab,a)]
=0
\end{cases}
$$
is equivalent to
$$
\begin{cases}
\eta_2
= \ep^2\P_{\ep}\nanorhs_2^{\ep,\mod}(\etab,a) =: \nsyst_2^{\ep}(\etab,a) \\
\\
a = \frac{1}{2\upsilon_{\ep}}\iota_{\ep}\big[\nanorhs_2^{\ep,\mod}(\etab,a)\big] =: \nsyst_3^{\ep}(\etab,a).
\end{cases}
$$
That is, our problem $\Thetab_{\ep}(\an_{\ep}(\etab,a)) = 0$ is equivalent to the fixed point problem
\begin{equation}\label{original fp}
\begin{cases}
\eta_1 = \A^{-1}\nanorhs_1^{\ep,\mod}(\etab,a)-\A^{-1}\K_2\eta_2 \\
\\
\eta_2 = \nsyst_2^{\ep}(\etab,a) \\
\\
a = \nsyst_3^{\ep}(\etab,a).
\end{cases}
\end{equation}

Now we can eliminate the difficulty with the term $\A^{-1}\K_2\eta_2$.  A pair $(\etab,a) \in E_q^1 \times E_q^1 \times \R$ is a fixed point solution to \eqref{original fp} if and only if $(\etab,a)$ solves
\begin{equation}\label{new fp}
\begin{cases}
\eta_1 = \A^{-1}\nanorhs_1^{\ep,\mod}(\etab,a)-\A^{-1}\K_2\nsyst_2^{\ep}(\etab,a) =: \nsyst_1^{\ep}(\etab,a)\\
\\
\eta_2 = \nsyst_2^{\ep}(\etab,a) \\
\\
a = \nsyst_3^{\ep}(\etab,a).
\end{cases}
\end{equation}
The term $\A^{-1}\K_2\nsyst_2^{\ep}(\etab,a)$ in the revised equation for $\eta_1$ turns out to have the ``right'' estimates in $\ep$ for our contraction mapping argument below.  We conclude that $\Thetab_{\ep}(\an_{\ep}(\etab,a)) = 0$ if and only if 
\begin{equation}\label{ultimate fp problem}
(\etab,a) = \left(\nsyst_1^{\ep}(\etab,a),\nsyst_2^{\ep}(\etab,a),\nsyst_3^{\ep}(\etab,a)\right) =: \nsystb^{\ep}(\etab,a),
\end{equation}
and we will solve this fixed point problem in the following section.

%% file: section_existence_properties_solutions.tex
\section{Existence and Properties of Solutions}
\label{main theorem section}

We model our existence proof on the approach of \cite{hw} for the small mass ratio, which in turn is a refinement of the contraction mapping proof in \cite{faver-wright}.  Let $q_{\star} > 0$ be as in Appendix \ref{nanopteron estimates from FW}.  For $r \ge 0$, set
$$
\X^r 
:=
\begin{cases}
E_{q_{\star}/2}^1 \times E_{q_{\star}/2}^1 \times \R, &r = 0 \\
\\
E_{q_{\star}}^r \times E_{q_{\star}}^r \times \R, &r > 1
\end{cases}
$$
and for $r$, $\ep$, $\tau > 0$, let
$$
\U_{\ep,\tau}^r
:= \set{(\etab,a) \in \X^r}{\norm{\etab}_{r,q_{\star}} \le \tau\ep, \ |a| \le \tau\ep^r}.
$$

\subsection{Existence and uniqueness of nanopteron solutions}

We base our contraction mapping argument on the following collection of estimates, which are proved in Appendix \ref{all the pretty estimates appendix}.

\begin{proposition}\label{main workhorse estimates}
There exists $\ep_{\star} > 0$ with the following properties.

\begin{enumerate}[label={\bf(\roman*)}]
\item There exists $\tau_{\star} > 0$ such that if $\ep \in (0,\ep_{\star})$, then

\begin{itemize}
\item (Mapping)
$$
(\etab,a) \in \U_{\ep_{\star},\tau_{\star}}^1
\Longrightarrow \nsystb^{\ep}(\etab,a) \in \U_{\ep,\tau_{\star}}^1
$$

\item (Lipschitz) 
$$
(\etab,a),(\grave{\etab},\grave{a}) \in \U_{\ep,\tau_{\star}}^1
\Longrightarrow \norm{\nsystb^{\ep}(\etab,a)-\nsystb^{\ep}(\grave{\etab},\grave{a})}_{\X^0} 
\le \frac{1}{2}\norm{(\etab,a)-(\grave{\etab},\grave{a})}_{\X^0}.
$$
\end{itemize}

\item (Bootstrapping) For all $r \ge 1$ and $\tau \in (0,a_{\per})$, there exists $\overline{\tau} = \overline{\tau}(\tau,r)$ such that if $\ep \in (0,\ep_{\star})$, then
$$
(\etab,a) \in \U_{\ep,\tau}^r \Longrightarrow \nsystb^{\ep}(\etab,a) \in \U_{\ep,\overline{\tau}}^{r+1}.
$$
\end{enumerate}
\end{proposition}

These estimates are essentially the same as the ones achieved in Lemma 8.1 of \cite{hw}, and a proof identical to that of their principal result, Theorem 8.2, produces the following solution to our ultimate fixed point problem \eqref{ultimate fp problem}.

\begin{theorem}\label{main theorem}
Let $\ep \in (0,\ep_{\star})$.  There exists a unique pair $(\etab_{\ep},a_{\ep}) \in \U_{\ep,\tau_{\star}}^1$ such that $\nsystb^{\ep}(\etab_{\ep},a_{\ep}) = (\etab_{\ep},a_{\ep})$.  This solution $(\etab_{\ep},a_{\ep})$ has the following additional properties:
\begin{enumerate}[label={\bf(\roman*)}]
\item $\etab_{\ep} \in \cap_{r=1}^{\infty} E_{\qstar}^r \times E_{\qstar}^r$;
\item For all $r \ge 0$, there is $C_r > 0$ such that 
$$
\norm{\etab_{\ep}}_{r,\qstar} \le C_r\ep \quadword{and} |a_{\ep}| \le C_r\ep^r
$$
for all $\ep \in (0,\ep_{\star})$.
\end{enumerate}
\end{theorem}

%
%
%

\subsection{Stegotons}  We now translate the fixed point solutions of Theorem \ref{main theorem} and Beale's ansatz into the language of relative displacements for our lattice problem \ref{nondimen lattice eqns}. 

\begin{corollary}\label{main corollary}
Let $\kappa > 1$ and $\beta \ne 0$ satisfy $\beta \ne -\kappa^3$.  There exist $\ep_{\star}$, $\qstar > 0$ such that for all $\ep \in (0,\ep_{\star})$, there is a solution for the relative displacements $r_j(t)$ of the nondimensionalized lattice equations \eqref{nondimen lattice eqns} in the form
$$
r_j(t)
= 
\bunderbrace{\kappa^{((-1)^j+1)/2}\frac{3\kappa(\kappa+1)}{c_{\kappa}^2(\beta+\kappa^3)}\ep^2\sech^2\left(\frac{\ep(j-\cep{t})}{2\sqrt{\alpha_{\kappa}}}\right)}{\text{\scalebox{.8}{{Note the extra factor of $\kappa$ for $j$ even.}}}}
+ v_j^{\ep}(\ep(j-\cep{t}))
+ p_j^{\ep}(j-\cep{t}),
$$
where 
\begin{enumerate}[label={\bf(\roman*)}]
\item $v_1^{\ep}$, $v_2^{\ep} \in \cap_{r=1}^{\infty} H_{\qstar}^r$.
\item $p_1^{\ep}$, $p_2^{\ep} \in \cap_{r=1}^{\infty} W^{r,\infty}$.
\item For each $r \ge 0$ there is a constant $C_r >0$ such that 
$$
\norm{v_j^{\ep}}_{H_{\qstar}^r} \le C_r\ep^3 \quadword{and} \norm{p_j^{\ep}}_{W^{r,\infty}} \le C_r\ep^{r}
$$
for all $\ep \in (0,\ep_{\star})$ and $j=1$, 2.
\item $p_1^{\ep}$ and $p_2^{\ep}$ are periodic with period $P_{\ep}$, and there is a constant $C >0$ such that $|P_{\ep}| \le C$ for all $\ep \in (0,\ep_{\star})$.
\end{enumerate}
\end{corollary}

\begin{proof}
The traveling wave ansatz \eqref{tw ansatz} with wave speed $c = \cep = \sqrt{c_{\kappa}^2+\ep^2}$ gave
$$
r_j(t) = \begin{cases}
p_1(j-\cep{t}), &j \text{ is odd} \\
p_2(j-\cep{t}), &j \text{ is even},
\end{cases}
$$
and the change of variables \eqref{diag cov} and the long wave scaling \eqref{lw scaling} converted $\pb = (p_1,p_2)$ into
$$
\pb(x) 
= (J\hb)(x)
= \ep^2(J\thetab(\ep\cdot))(x)
= \ep^2(J^{\ep}\thetab)(\ep{x}).
$$
With $\thetab = \sigmab + a_{\ep}\varphib_{\ep}^{a_{\ep}}+\etab_{\ep}$ from Beale's ansatz \eqref{beale} and Theorem \ref{main theorem}, we find
$$
\pb(x) = \ep^2(J^{\ep}\sigmab)(\ep{x})
+ \ep^2a_{\ep}(J^{\ep}\varphib_{\ep}^{a_{\ep}})(\ep{x})
+ \ep^2(J^{\ep}\etab_{\ep})(\ep{x}).
$$

We now want to isolate what will be the lowest order term in $\ep$.  With $J^0$ defined in \eqref{J0}, we have
$$
\pb(x) = \ep^2(J^0\sigmab)(\ep{x})
+ \ep^2a_{\ep}(J^{\ep}\varphib_{\ep}^{a_{\ep}})(\ep{x})
+ \ep^2\left((J^{\ep}\etab_{\ep})(\ep{x}) + ((J^{\ep}-J^0)\sigmab)(\ep{x})\right),
$$
where
\begin{equation}\label{rel disp periodic est}
\ep^2\norm{a_{\ep}(J^{\ep}\varphib_{\ep}^{a_{\ep}})(\ep\cdot)}_{W^{r,\infty}}
\le \ep^2(C_r\ep^r)\ep^r\norm{J^{\ep}\varphib_{\ep}^{a_{\ep}}}_{W^{r,\infty}} 
\le C_r\ep^{r+2}
\end{equation}
and
\begin{equation}\label{rel disp decay est}
\ep^2\norm{J^{\ep}\etab_{\ep} + (J^{\ep}-J^0)\sigmab}_{r,q_{\star}} \le C_r\ep^3.
\end{equation}
by the estimates in Theorem \ref{main theorem} and \eqref{Jep0 estimate}.  We abbreviate
$$
\begin{pmatrix*}
p_1^{\ep}(x) \\
p_2^{\ep}(x)
\end{pmatrix*}
:= \ep^2a_{\ep}(J^{\ep}\varphib_{\ep}^{a_{\ep}})(\ep{x})
\quadword{and}
\begin{pmatrix*}
v_1^{\ep}(x) \\
v_2^{\ep}(x)
\end{pmatrix*}
:=
\ep^2\left((J^{\ep}\etab_{\ep})({x}) + ((J^{\ep}-J^0)\sigmab)({x})\right).
$$
We get the estimates for $p_j$ and $v_j$ from \eqref{rel disp periodic est} and \eqref{rel disp decay est}.  For the period of $p_j$, observe that by Theorem \ref{periodic main theorem}
$$
(J^{\ep}\varphib_{\ep}^{a_{\ep}})(\ep{x})
= (J^{\ep}\nub(\omega_{\ep}^{a_{\ep}}\cdot))(\ep{x})
+ (J^{\ep}\psib_{\ep}^{a_{\ep}}(\omega_{\ep}^{a_{\ep}}\cdot))(\ep{x})
= (J^{\ep}\nub)(\ep\omega_{\ep}^{a_{\ep}}x)
+ (J^{\ep}\psib_{\ep}^{a_{\ep}})(\ep\omega_{\ep}^{a_{\ep}}x),
$$
where $\ep\omega_{\ep}^{a_{\ep}}$ is uniformly bounded in $\ep$ and $\nub$ and $\psib_{\ep}^{a_{\ep}}$ are $2\pi$-periodic.

Finally, 
$$
J^0\sigmab 
= \begin{bmatrix*}[r]
1/\kappa &1 \\
1 &-1
\end{bmatrix*}
\begin{pmatrix*}
\sigma \\
0
\end{pmatrix*}
= \begin{pmatrix*} 
\sigma/\kappa \\
\sigma
\end{pmatrix*},
$$
and so for $j$ odd we have
$$
r_j(t)
= \ep^2\frac{1}{\kappa}\sigma(\ep(j-\cep{t}))
+ v_1^{\ep}(\ep(j-\cep{t})) + p_1^{\ep}(j-\cep{t}),
$$
while for $j$ even,
\[
r_j(t) 
= \ep^2\sigma(\ep(j-\cep{t}))
+ v_2^{\ep}(\ep(j-\cep{t})) + p_2^{\ep}(j-\cep{t}). \qedhere
\]
\end{proof}

\begin{figure}
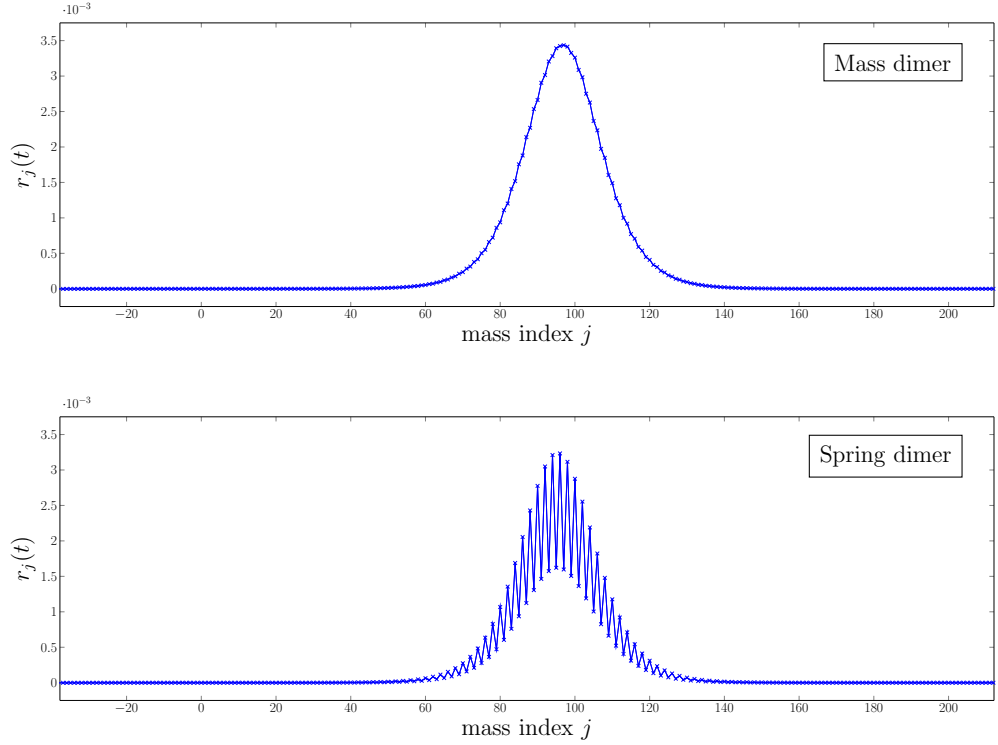

{\centering
\begin{tabular}{c}
\input{fig_mass_no_stegoton} \\
\\
\input{fig_spring_stegoton}
\end{tabular}
}
\caption{Relative displacements for mass and spring dimers}
\label{fig - rel disp}
\end{figure}

\begin{remark}
We note that the leading order terms in $r_j(t)$ differ by a factor of $\kappa$ depending on whether $j$ is even or odd. This is a feature not present in the relative displacements for the mass dimers (cf. Corollary 6.4 in \cite{faver-wright}), and it leads to the spiky, jagged graphs evocatively called ``stegotons'' in \cite{leveque1,leveque2}.  This behavior for lattices with alternating springs was observed in \cite{gmwz} and using the same simulations from that paper, we plot in Figure \ref{fig - rel disp} (approximate) solutions to the lattice equations \eqref{nondimen lattice eqns}.  There is clearly a marked difference between the two graphs, even though we have chosen the essential parameters $w$ and $\kappa$ to be equal to 2 in each case. 
\end{remark}

Moreover, the estimates in Corollary \ref{main corollary} are consistent with the results of \cite{gmwz}, which establishes 
\begin{equation}\label{gmwz approx}
r_j(t) = \frac{\ep^2}{\Kappa_j}\left(U_-(\ep(j-c_{\kappa}t),\ep^3t) + U_+(\ep(j+c_{\kappa}t),\ep^3t)\right) + \O(\ep^{5/2}),
\end{equation}
where 
$$
\Kappa_j = \begin{cases}
1, &j \text{ is even} \\
\kappa, &j \text{ is odd}
\end{cases}
$$
and $U_{\pm} = U_{\pm}(X,T)$ solve the KdV equations
\begin{equation}\label{KdV GMWZ}
\mp\frac{1}{c_{\kappa}} \partial_T[U_{\pm}]
+ \frac{1}{6}\left(\frac{1-\kappa+\kappa^2}{(1+\kappa)^2}\right)\partial_X^3[U_{\pm}]
+ \frac{\kappa}{\kappa+1}\left(\frac{\beta}{\kappa^3}+1\right)U_{\pm}\partial_X[U_{\pm}]
= 0.
\end{equation}
Making the traveling wave ansatz
$$
U_{\pm}(X,T) = V\left(X\pm\frac{T}{2c_{\kappa}}\right)
$$
rearranging, and integrating, \eqref{KdV GMWZ} becomes
$$
\alpha_{\kappa}V'' - V + c_{\kappa}^2\frac{\kappa}{\kappa+1}\left(\frac{\beta}{\kappa^3}+1\right)V^2 
= 0,
$$
which is precisely the ordinary differential equation \eqref{KdV ode} that we derived in our study of the formal long wave limit at $\ep = 0$. Hence $V = \sigma$ as defined in \eqref{sigma soliton defn} and
$$
U_{\pm}(X,T) = \sigma\left(X\pm\frac{T}{2c_{\kappa}}\right),
$$
and so in the approximation \eqref{gmwz approx} for $r_j(t)$ we have
$$
U_{-}(\ep(j- c_{\kappa}t),\ep^3t)
= \sigma\left(\ep(j - c_{\kappa}t) - \frac{\ep^3t}{2c_{\kappa}}\right) 
= \sigma\left(\ep\left(j -\left(c_{\kappa} + \frac{1}{2c_{\kappa}}\ep^2\right)t\right)\right).
$$
Since
$$
\cep = \sqrt{c_{\kappa}^2+\ep^2} = c_{\kappa} + \frac{1}{2c_{\kappa}}\ep^2 + \O(\ep^4)
$$
and $\sigma' \in L^{\infty}$, we have (for $|t| \le T_0$)
$$
\sigma\left(\ep\left(j - \left(c_{\kappa} + \frac{1}{2c_{\kappa}}\ep^2\right)t\right)\right)
= \sigma(\ep(j- \cep{t})) + \O(\ep^5).
$$
Last, since $\sigma \in L^{\infty}$, we have
$$
\ep^2U_+(\ep(j+c_{\kappa}t),\ep^3t)
= \ep^2\sigma\left(\ep(j+c_{\kappa}t)+\frac{\ep^3t}{2c_{\kappa}}\right)
= \O(\ep^2).
$$
All together, \eqref{gmwz approx} becomes
$$
r_j(t)
= \frac{\ep^2}{\Kappa_j}\sigma(\ep(j-\cep{t})) + \O(\ep^2),
$$
which agrees with the expression for $r_j(t)$ in Corollary \ref{main corollary} and the subsequent estimates.

%% file: appendix_function_spaces_master.tex
\section{Function Spaces} 
\label{function spaces appendix}

\input{appendix_calculus}
\input{appendix_fourier_transform}
\input{appendix_sobolev_spaces}

\input{appendix_fourier_multipliers}
\input{appendix_function_estimates}

%% file: appendix_calculus.tex
\subsection{Calculus}

\subsubsection{Leibniz's rule} We will often use Leibniz's rule for an arbitrary derivative of a product: 
$$
\partial_X^r[fg] 
= \sum_{k=0}^r {r \choose k} \partial_X^k[f]\partial_X^{r-k}[g].
$$

\subsubsection{Fa\'{a} di Bruno's formula} 
We employ the convenient expression of Fa\'{a} di Bruno's formula for the chain rule found in \cite{mortini}.  For $k$, $r \in \N$ with $k \le r$, let
$$
\Sigma_k^r 
= \set{\sigmab \in \N^k}{\sigma_1 \ge \sigma_2 \ge \cdots \ge \sigma_k \ge 1, \ |\sigmab| = r},
$$
where
$$
|\sigmab| 
:= \sum_{j=1}^k \sigma_j.
$$

\begin{remark}\label{FdB remark}
\begin{enumerate}[label = {\bf(\roman*)}]
\item It is apparent from the definition that $\Sigma_1^r = \{r\}$.
\item If $\sigmab \in \Sigma_k^r$ with $2 \le k \le r$, then $\sigma_j < r$ for all $j$.
\end{enumerate}
\end{remark}

\begin{theorem}[Fa\'{a} di Bruno]
Let $N,f \in \C^r(\R)$.  Then
$$
\partial_X^r[N(f)] 
= \sum_{k=1}^r \partial_X^k[N](f)\sum_{\sigmab \in \Sigma_k^r} C_{\sigmab}\prod_{j=1}^k \partial_X^{\sigma_j}[f],
$$
where the $C_{\sigmab}$ are positive constants that depend on $r$ and $k$ but are independent of $f$ and $X$.  
\end{theorem}

%% file: appendix_fourier_transform.tex
\subsection{The Fourier transform}\label{fourier transform appendix}

\subsubsection{The periodic Fourier transform}
For an integrable, $2P$-periodic function $f$ defined on $\R$, we set
$$
\ft[f](k)
= \hat{f}(k)
=: \frac{1}{\sqrt{2P}}\int_{-P}^{P} f(X)e^{ik\pi{X}/P} \dX.
$$

\subsubsection{The Fourier transform on $\R$}
For a function $f \in L^2(\R)$, we set
$$
\ft[f](k) 
= \hat{f}(k) 
:= \frac{1}{\sqrt{2\pi}}\int_{\R} f(X)e^{-ikX} \dX
\quadword{and}
\ft^{-1}[f](k) 
= \widecheck{f}(k) 
:= \frac{1}{\sqrt{2\pi}}\int_{\R} f(X)e^{ikX} \dX.
$$

%% file: appendix_sobolev_spaces.tex
\subsection{Sobolev Spaces}\label{sobolev spaces appendix}

\input{appendix_sobolev_spaces_defns}

%% file: appendix_sobolev_spaces_defns.tex
\subsubsection{Sobolev spaces on $\R$}\label{basic definitions and norms for sobolev spaces}

For an integer $r \ge 1$ we let $H^r$ be the usual Sobolev space of functions in $L^2(\R)$ whose $j$th weak partial derivative, $\partial_X^j[f]$, exists and is in $L^2$ for $j = 0,\ldots, r$.  Our preferred norm for $H^r$ is 
\begin{equation}\label{preferred R sobolev norm}
\norm{f}_{H^r} := \norm{f}_{L^2} + \norm{\partial_X^r[f]}_{L^2}.
\end{equation}

\subsubsection{Sobolev spaces of exponentially decaying functions} For $r$, $q \ge 0$ we set
$$
H_q^r 
:= \set{f \in H^r}{\cosh(q\cdot)f \in H^r} = \set{f \in H^r}{\cosh^q(\cdot)f \in H^r}
$$
and, for integer $r$, we use the norm
\begin{equation}\label{Hrq norm defn}
\norm{f}_{r,q} 
:= \norm{\cosh(q\cdot)f}_{L^2} + \norm{\cosh(q\cdot)\partial_X^r[f]}_{L^2}.
\end{equation}
This norm is equivalent to the norm $f \mapsto \norm{\cosh(q\cdot)f}_{H^r}$, with $\norm{\cdot}_{H^r}$ defined as in \eqref{preferred R sobolev norm}, and also the norm $f \mapsto \norm{\cosh^q(\cdot)f}_{L^2} + \norm{\cosh^q(\cdot)\partial_X^r[f]}_{L^2}$.  

We will use the following bevy of elementary estimates on functions in $H_q^r$ frequently. 

\begin{proposition}\label{all the hrq basic estimates prop}
Let $r$, $q > 0$.  
\begin{enumerate}[label={\bf(\roman*)}]
\item $\norm{f}_{L^{\infty}} \le C_r\norm{f}_{H^r}$
\item $\norm{fg}_{H^r} \le C_r\norm{f}_{W^{r,\infty}}\norm{g}_{H^r}$
\item $\norm{\partial_X^k[\cosh(q\cdot)f]}_{L^{\infty}} \le C_{r-k}\norm{f}_{r,q}$

\item $\norm{f}_{H^r} \le C_{r,q}\norm{f}_{r,q}$

\item $\norm{fg}_{r,q} \le C_{r,q}\norm{f}_{r,q}\norm{g}_{r,q}$

\item $\norm{fg}_{r,q} \le C_r\norm{g}_{W^{r,\infty}}\norm{f}_{r,q}$

\item $\norm{f}_{r,q} \le C_{r,q'-q}\norm{f}_{r,q'}, \ q \le q'$
\end{enumerate}
\end{proposition}

\subsubsection{Periodic Sobolev spaces} \label{periodic sobolev spaces appendix}
For an integrable, $2\pi$-periodic function $f$ on $\R$, we set
$$
\norm{f}_{L_{\per}^2} := \left(\sum_{k \in \Z} |\hat{f}(k)|^2\right)^{1/2}
$$
and, for integer $r$,
\begin{equation}\label{Hrper norm}
\norm{f}_{H_{\per}^r} := \norm{f}_{L_{\per}^2} + \norm{\partial_X^r[f]}_{L_{\per}^2}.
\end{equation}
We will also use the inner product
$$
\ip{f}{g}_{H_{\per}^r} := \sum_{k \in \Z} (1+k^2)^r\hat{f}(k)\overline{\hat{g}(k)}.
$$
We define
$$
H_{\per}^r = \set{f \in L_{\per}^2}{\norm{f}_{H_{\per}^2} < \infty}.
$$
For $r = 0$ we let $H_{\per}^0 := L_{\per}^2$.  The norm \eqref{Hrper norm} is of course equivalent to the more familiar norm 
$$
f \mapsto \left(\sum_{k \in \Z} (1+k^2)^r|\hat{f}(k)|^2\right)^{1/2}
$$
on $H_{\per}^r$; see \cite{kress} for details.  We prefer \eqref{Hrper norm} for its similarity to the highly convenient norms \eqref{preferred R sobolev norm} and \eqref{Hrq norm defn}.

%

We will also need three familiar estimates:
\begin{itemize}[leftmargin=*]
\item The Sobolev product estimate
$$
\norm{fg}_{H_{\per}^r} \le C_r\norm{f}_{H_{\per}^r}\norm{g}_{H_{\per}^r};
$$
\item 
The Sobolev embedding estimate
$$
\norm{f}_{\C_{\per}^{r-1}} \le C_r\norm{f}_{H_{\per}^r}
$$
for $r \ge 1$, where
$$
\norm{f}_{\C_{\per}^k} := \sum_{j=0}^k \norm{\partial_X^j[f]}_{L_{\per}^{\infty}};
$$
\item And the Fourier transform estimate
$$
|\hat{f}(k)| \le \sqrt{2\pi}\norm{f}_{L^{\infty}} \le C_r\norm{f}_{H_{\per}^r}.
$$
\end{itemize}

%% file: appendix_fourier_multipliers.tex
\subsection{Fourier multipliers}
\label{fourier multipliers appendix}

Let $\tmu \in L^{\infty}(\R)$.  For $f \in L^2(\R)$, we define the Fourier multiplier operator $\mu$ with symbol $\tmu$ by
\begin{equation}\label{fm for L2}
\mu{f} := \ft^{-1}[\mu\hat{f}].
\end{equation}
That is, $\hat{\mu{f}}(k) = \tmu(k)\hat{f}(k)$, $k \in \R$.

Similarly, for an integrable, $2P$-periodic $f$ on $\R$, we define $\mu{f}$ to be the function whose Fourier coefficients are $\tmu(k\pi/P)\hat{f}(k)$, which is to say,
\begin{equation}\label{fm for periodic}
\mu{f}(X) := \sum_{k \in \Z} e^{ik\pi{X}/P}\tmu\left(\frac{k\pi}{P}\right)\hat{f}(k).
\end{equation}

In both cases, we have the following important scaling property of Fourier multipliers: if $\omega \in \R$ and $f$ is a function, then $\mu$ acts on the scaled function $f(\omega\cdot)$ by 
$$
\mu[f(\omega\cdot)] = (\mu^{\omega}f)(\omega \cdot),
$$
where $\mu^{\omega}$ is the Fourier multiplier with symbol $\tilde{\mu^{\omega}}(k) := \tmu(\omega{k})$.  

Finally, if $f \in L^2(\R)$ and $g$ is integrable and $2P$-periodic, then we set 
$$
\mu(f+g) := \mu{f} + \mu{g},
$$
where $\mu{f}$ is defined per \eqref{fm for L2} and $\mu{g}$ by \eqref{fm for periodic}.

%% file: appendix_function_estimates.tex
\subsection{Estimates on function spaces}\label{function estimates appendix}
In this appendix we prove a number of estimates that directly facilitate our treatment of the higher-order terms in the spring forces.  All of these estimates treat composition operators: given a map $N$ and functions $f$ and $\grave{f}$, we need to control the norms of $N(f)$ and $N(f)-N(\grave{f})$ in a suitable function space.  If our functions were always in $H^r$, this would be a well-understood problem \cite{brezis-mironescu}; to bound $N(f)$, we could rely on, for example, the proof of estimate (2.4) in \cite{moser} or Proposition 3.9 in \cite{taylor}, both of which state, roughly, that if $N \in \C^{\infty}(\R)$ and $N(0) = 0$, then
$$
\norm{N(f)}_{H^r} \le C(1+\norm{f}_{H^r}),
$$
where $C$ depends on $\norm{f}_{L^{\infty}}$. 

However, our functions will always be the superposition of a function $f \in H_q^r$ with a function $\varphi \in W^{r,\infty}$.  Since $\partial_X^r[f]$ need not be bounded and $\varphi$ need not be square-integrable, the sum $f+\varphi$ belongs to neither of the spaces $H_q^r$ nor $W^{r,\infty}$, so we cannot appeal to the existing literature for our estimates.  Moreover, the $W^{r,\infty}$-norm of $\varphi$ will depend delicately on the small parameter $\ep$, and we want careful, uniform estimates in $\ep$.  So, we develop our composition operator estimates from scratch.  Our proofs are mostly straightforward, but thorough, applications of Fa\`{a} di Bruno's rule to express precisely the derivatives of certain compositions.  

\subsubsection{A mapping estimate}

\begin{proposition}\label{mapping prop}
Let $r \in \N$, $N \in \C^r(\R)$, and $C_{\star}$, $q > 0$.  There exists an increasing map $\M \colon (0,\infty) \to (0,\infty)$ with the following property. Let $a$, $\ep \in (0,1)$.  Suppose $f$, $g \in H_q^r$ and $\varphi \in W^{r,\infty}$ with $\norm{\partial_X^j[\varphi]}_{L^{\infty}} \le C_{\star}\ep^{-j}$ for $j = 0, \ldots,r$.  Then
\begin{equation}\label{mapping}
\norm{fN(\ep(g+a\varphi))}_{r,q} \le \M[\norm{g}_{r,q}](1+|a|\ep^{1-r})\norm{f}_{r,q}.
\end{equation}
The function $\M$ depends on $N$, $r$, $q$, and $C_{\star}$, but not on $f$, $\ep$, or $\varphi$.  
\end{proposition}

\begin{remark}\label{bounding remark}
Each of the estimates that we prove in this appendix will include a rather opaque nonnegative factor $\M$ that depends, by its (convoluted) construction, on the norms of one or more functions involved in the estimate.  So, we can think of $\M$ as a map from $\R_+^n$ to $\R_+$ for some $n \ge 1$, where $\R_+ = [0,\infty)$.  It turns out that $\M$ will always have the property that
\begin{equation}\label{M is locally bounded}
\sup_{\substack{\vb \in \R_+^n \\ |\vb| \le r}} \M[\vb] 
< \infty
\end{equation}
for any $r > 0$, where $|\vb| = \sum_{k=1}^n v_k$.  Then we may define a new map $\M_{\incr} \colon \R_+^n \to \R_+$ by
$$
\M_{\incr}[\ub] 
:= \sup_{\substack{\vb \in \R_+^n \\ |\vb| \le |\ub|}} \M[\vb].
$$
This map $\M_{\incr}$ enjoys three properties:
\begin{itemize}[leftmargin=*]
\item $\M_{\incr}$ is ``radially increasing'' in the sense that if $|\ub| \le |\grave{\ub}|$, then $\M_{\incr}[\ub] \le \M_{\incr}[\grave{\ub}]$.  In particular, if $n =1$, then $\M$ is increasing in the usual sense.
\item As a consequence of the first property, $\M_{\incr}$ is locally bounded in the sense that
$$
\sup_{\ub \in \R_+^n \\ |\ub| \le r} \M_{\incr}[\ub] < \infty
$$
for any $r > 0$.
\item It is obvious from the definition of $\M_{\incr}$ that $\M[\ub] \le \M_{\incr}[\ub]$.
\end{itemize}
The first two properties of $\M_{\incr}$ above will be very useful for the nanopteron estimates. The third property allows us to replace the function $\M$ that we construct in a given proof by its relative $\M_{\incr}$, and we will do so without further comment. So, for example, when we prove Proposition \ref{mapping prop}, our proof will only demonstrate the property \eqref{M is locally bounded}.  
\end{remark}

\begin{remark}
In the proof of Proposition \ref{mapping prop} and the following proofs we will have a great many constants that depend more or less innocuously on different parameters. A constant $C_r$ depends only on $r$; a constant $C_{r,q}$ depends only on $r$ and $q$; and a constant $C_{\star,r,q}$ depends only on $C_{\star}$, $r$, and $q$.  The value of these constants may change from line to line, but their dependence remains rests firmly and solely on their subscripts.  We will maintain this convention throughout the rest of this appendix and Appendices \ref{periodic solutions appendix} and \ref{all the pretty estimates appendix}.
\end{remark}

Now we are ready to prove Proposition \ref{mapping prop}.

\begin{proof} 
We will construct three functions $\M_0$, $\M_1$, and $\M_2$ in this proof as we progress toward the ultimate $\M$ that satisfies \eqref{mapping}.  To simplify notation, we will write $\M_0[g] = \M[\norm{g}_{r,q}]$, etc.  Since
\begin{equation}\label{first mapping equation}
\norm{fN(\ep(g+\varphi))}_{r,q} 
= \norm{\cosh(q\cdot)fN(\ep(g+a\varphi))}_{L^2} 
+ \norm{\cosh(q\cdot)\partial_X^r[fN(\ep(g+a\varphi))]}_{L^2}.
\end{equation}
it suffices to verify that each $L^2$-norm in the sum above has an upper bound of the form given on the right side of \eqref{mapping}.  

First, let
$$
\I(g) 
= \set{Y \in \R}{|Y| \le C_{r,q}\norm{g}_{r,q} + C_{\star}}.
$$
Since
\begin{align*}
\ep|g(X)+a\varphi(X)| 
&\le \norm{g}_{L^{\infty}} + |a|\norm{\varphi}_{L^{\infty}} \\
\\
&\le C_{r,q}\norm{g}_{r,q} + C_{\star}\ep^0  \\
\\
&\le C_{r,q}\norm{g}_{r,q}+C_{\star},
\end{align*}
we have
$$
\norm{\partial_X^{k}[N](\ep(g+a\varphi))}_{L^{\infty}} 
\le \norm{N}_{W^{r,\infty}(\I(g))}.  
$$
Then the first term in \eqref{first mapping equation} is easy to handle:
$$
\norm{\cosh(q\cdot)fN(\ep(g+a\varphi))}_{L^2} 
\le \norm{N(\ep(g+a\varphi))}_{L^{\infty}}\norm{\cosh(q\cdot)f}_{L^2}
\le \norm{N}_{W^{r,\infty}(\I(g))}\norm{\cosh(q\cdot)f}_{r,q}.
$$

Next, Leibniz' rule reduces the study of the second term in \eqref{first mapping equation} to estimating terms of the form
$$
\norm{\cosh(q\cdot)\partial_X^m[f]\partial_X^n[N(\ep(g+a\varphi))]}_{L^2},
$$
where $0 \le m,n \le r$ and $m+n \le r$.  It turns out that estimating the $m=0, n = r$ term is both the most complicated and the most instructive term, so we do it first.  That is, we estimate
$$
\norm{\cosh(q\cdot)f\partial_X^r[N(\ep(g+a\varphi))]}_{L^2}.
$$

We use Fa\'{a} di Bruno's rule to expand
$$
\partial_X^r[N(\ep(g+a\varphi))] 
= \sum_{k=1}^r \partial_X^k[N](\ep(g+a\varphi))\sum_{\sigmab \in \Sigma_k^r} C_{\sigmab}\prod_{j=1}^k \partial_X^{\sigma_j}[\ep(g+a\varphi)].
$$

A first pass then reduces our estimate to
\begin{equation}\label{first pass}
\norm{\cosh(q\cdot)f\partial_X^r[N(\ep(g+a\varphi))] }_{L^2}
\le C_r\norm{N}_{W^{r,\infty}(\I(g))}\sum_{k=1}^r\sum_{\sigmab \in \Sigma_k^r}\norm{\cosh(q\cdot)f\prod_{j=1}^k \ep\partial_X^{\sigma_j}[g+a\varphi]}_{L^2}.
\end{equation}

When $k=1$, we have $\Sigma_1^r = \{r\}$ by Remark \ref{FdB remark}, and
\begin{align*}
\norm{\cosh(q\cdot)f\ep\partial_X^r[g+a\varphi]}_{L^2} 
&\le \norm{\cosh(q\cdot)f}_{L^{\infty}}\norm{\partial_X^r[g]}_{L^2} + \norm{\cosh(q\cdot)f}_{L^2}\ep|a|\norm{\partial_X^r[\varphi]}_{L^{\infty}} \\
\\
&\le C_{r,q}\norm{f}_{r,q}\norm{g}_{r,q} + \norm{f}_{r,q}|a|C_{\star}\ep^{1-r} \\
\\
&\le \M_0[g]\left(1+|a|\ep^{1-r}\right)\norm{f}_{r,q},
\end{align*}
where we have set 
$$
\M_0[g] = \max\left\{C_{r,q}\norm{g}_{r,q}, C_{\star}\right\}.
$$

When $2 \le k \le r$, all of the factors in the product in \eqref{first pass} will be $L^{\infty}$ because the order of each derivative $\partial_X^{\sigma_j}[g]$ will be at most $r-1$.  Then
\begin{equation}\label{Eq3}
\begin{aligned}
\norm{\prod_{j=1}^k \ep\partial_X^{\sigma_j}[g+a\varphi]}_{L^{\infty}}
&\le \prod_{j=1}^k \left(C_{r,q}\norm{g}_{r,q} + \ep|a|\norm{\partial_X^{\sigma_j}[\varphi]}_{L^{\infty}}\right) \\
\\
&\le C_{\star,r,q}\prod_{j=1}^k \left(\norm{g}_{r,q} + |a|\ep^{1-\sigma_j}\right). 
\end{aligned}
\end{equation}
We use Lemma \ref{product of sums lemma} to write this product as 
\begin{equation}\label{Eq4}
\begin{aligned}
\prod_{j=1}^k \left(\norm{g}_{r,q} + |a|\ep^{1-\sigma_j}\right) 
&= \sum_{\substack{\alphab,\betab \in \{0,1\}^k \\ \alpha_j + \beta_j = 1 \ \forall{j}}} \prod_{j=1}^k \norm{g}_{r,q}^{\alpha_j}\left(|a|\ep^{1-\sigma_j}\right)^{\beta_j} \\
\\
&= \norm{g}_{r,q}^k + \sum_{\substack{\alphab,\betab \in \{0,1\}^k \\ \alpha_j + \beta_j = 1 \ \forall{j} \\ \betab \ne \mathbf{0}}}\prod_{j=1}^k \norm{g}_{r,q}^{\alpha_j}\left(|a|\ep^{1-\sigma_j}\right)^{\beta_j}.
\end{aligned}
\end{equation}
We focus on the second term.  First, because $0 \le \medsum_{j=1}^k \alpha_j \le k$, we have
\begin{equation}\label{Eq1}
\prod_{j=1}^k \norm{g}_{r,q}^{\alpha_j} 
= \norm{g}_{r,q}^{\sum_{j=1}^k\alpha_j}
\le \sum_{j=0}^k \norm{g}_{r,q}^j 
=:\M_1[g].
\end{equation}
Next, because $|a| < 1$ and $\ep < 1$, we have
$$
\prod_{j=1}^k \left(|a|\ep^{1-\sigma_j}\right)^{\beta_j} 
\le |a|\ep\prod_{j=1}^k \ep^{-\sigma_j\beta_j}
= |a|\ep\ep^{-\sum_{j=1}^k\sigma_j\beta_j}.
$$
Since $0 \le \sigma_j\beta_j \le \sigma_j$, this in turn becomes
\begin{equation}\label{Eq2}
|a|\ep\ep^{-\sum_{j=1}^k\sigma_j\beta_j}
\le |a|\ep\ep^{-\sum_{j=1}^k \sigma_j}
= |a|\ep\ep^{|\sigmab|}
= |a|\ep^{1-r}.
\end{equation}
Here we have used the stipulation $|\sigmab| = r$ from Remark \ref{FdB remark}.  

All together, \eqref{Eq1} and \eqref{Eq2} imply
$$
\sum_{\substack{\alphab,\betab \in \{0,1\}^k \\ \alpha_j + \beta_j = 1 \ \forall{j} \\ \betab \ne \mathbf{0}}}\prod_{j=1}^k \norm{g}_{r,q}^{\alpha_j}\left(|a|\ep^{1-\sigma_j}\right)^{\beta_j}
\le 
\M_1[g]\sum_{\substack{\alphab,\betab \in \{0,1\}^k \\ \alpha_j + \beta_j = 1 \ \forall{j} \\ \betab \ne \mathbf{0}}} |a|\ep^{1-r}
\le 
r^2\M_1[g]|a|\ep^{1-r},
$$
and so \eqref{Eq3} and \eqref{Eq4} imply
\begin{equation}\label{product bound for mapping}
\norm{\prod_{j=1}^k\ep\partial_X^{\sigma_j}[g+a\varphi]}_{L^{\infty}} 
\le \M_1[g] + r^2\M_1[g]|a|\ep^{1-r} 
\le \M_2[g]\left(1+|a|\ep^{1-r}\right),
\end{equation}
where $\M_2[g] = r^2\M_1[g]$.
Returning to \eqref{first pass}, we find
\begin{align*}
\norm{\cosh(q\cdot)f\partial_X^r[N(\ep(g+a\varphi))]}_{L^2}
&\le C_r\norm{N}_{W^{r,\infty}(\I(g))}\M_0[g]\left(1+|a|\ep^{1-r}\right)\norm{f}_{r,q} \\
\\
&+C_r\norm{N}_{W^{r,\infty}(\I(g))}\sum_{k=2}^r\sum_{\sigmab \in \Sigma_k^r} \norm{\cosh(q\cdot)f}_{L^2}\M_2[g]\left(1+|a|\ep^{1-r}\right) \\
\\
&\le \bunderbrace{C_r\norm{N}_{W^{r,\infty}(\I(g))}\M_2[g]}{\M[g]}\left(1+|a|\ep^{1-r}\right)\norm{f}_{r,q}. 
\end{align*}
This estimate has the same form as the right side of \eqref{mapping}, so we have completed our work on $\norm{\cosh(q\cdot)f\partial_X^r[N(\ep(g+a\varphi))]}_{L^2}$.
To treat the other terms
$$
\norm{\cosh(q\cdot)\partial_X^m[f]\partial_X^n[N(\ep(g+a\varphi))]}_{L^2}
$$
where $0 \le m,n \le r, m+n \le r,$ and $(m,n) \ne (0,r)$, note that these strictures on $m$ and $n$ imply $n \le r-1$.  So, when we expand $\partial_X^n[N(\ep(g+a\varphi))]$ with Fa\`{a} di Bruno's rule, all derivatives $\partial_X^j[g+\varphi]$ will be $L^{\infty}$, while of course $\cosh(q\cdot)\partial_X^m[f]$ will be $L^2$.  Then we proceed exactly as we did above in the long treatment of the case $2 \le k \le r$.
\end{proof}

\subsubsection{A Lipschitz estimate in $H_q^r$}\label{Hrq lipschitz section}

\begin{proposition}\label{general hrq lipschitz prop}
Let  $r \in \N$, $N \in \C^{r+1}(\R)$, and $C_{\star},q > 0$.  There exists a radially increasing map $\M \colon \R_+^2 \to \R_+$ with the following property.  Let $a$, $\ep \in (0,1)$.  Suppose $f$, $\grave{f} \in H_q^r$ and $\varphi \in W^{r,\infty}$ with $\norm{\partial_X^j[\varphi]}_{L^{\infty}} \le C_{\star}\ep^{-j}$ for $j = 0, \ldots,r$.  Then
\begin{equation}\label{general hrq lipschitz}
\bignorm{N(\ep(f+a\varphi))-N(\ep(\grave{f}+a\varphi))}_{r,q} 
\le \M[\norm{f}_{r,q},\bignorm{\grave{f}}_{r,q}]\left(1+|a|\ep^{1-r}\right)\bignorm{f-\grave{f}}_{r,q}.
\end{equation}
The function $\M$ depends on $N$, $r$, $q$, and $C_{\star}$ but not on $f$, $\grave{f}$, $\ep$, or $\varphi$.
\end{proposition}

\begin{proof}
Following our notation in the proof of Proposition \ref{mapping prop}, we will build $\M$ out of several functions $\M_0,\ldots,\M_5$, and suppress norms within these functions, writing, for example,
$$
\M_0[\norm{f}_{r,q},\bignorm{\grave{f}}_{r,q}] 
= \M_0[f,\grave{f}] 
\quadword{and}
\M_1[\norm{f}_{r,q}] = \M_1[f].
$$

We have
\begin{align*}
N(\ep(f+a\varphi))-N(\ep(\grave{f}+a\varphi))  
&= \bunderbrace{\cosh(q\cdot)\left(N(\ep(f+a\varphi))-N(\ep(\grave{f}+a\varphi))\right)}{\Delta_1} \\
\\
&+ \bunderbrace{\cosh(q\cdot)\partial_X^r\left[N(\ep(f+a\varphi))-N(\ep(\grave{f}+a\varphi))\right]}{\Delta_2}.
\end{align*}
Since
$$
|\ep(f(X)+a\varphi(X))| 
\le \norm{f}_{L^{\infty}} + \norm{\varphi}_{L^{\infty}} 
\le C_{r,q}\norm{f}_{r,q} + C_{\star},
$$
for any $0 \le k \le r$ we have
\begin{equation} \label{N wrinfty est}
\norm{\partial_X^k[N](\ep(f + a\varphi))}_{L^{\infty}} 
\le \norm{N}_{W^{r+1,\infty}(\I(f,\grave{f}))}
=: \M_0[f,\grave{f}], 
\end{equation}
where
$$
\I(f,\grave{f}) 
= \set{Y \in \R}{|Y| \le C_{r,q}\norm{f}_{r,q} + C_{r,q}\bignorm{\grave{f}}_{r,q} + 2C_{\star}}.
$$
If we replace $f$ with $\grave{f}$ in \eqref{N wrinfty est}, the same estimate still holds. We will use this estimate frequently throughout the rest of the proof, starting with a bound on $\Delta_1^2$:
\begin{equation}\label{Delta1 estimate}
\begin{aligned}
\norm{\Delta_1}_{L^2}^2 &=
\int_{\R} \cosh^2(qX)\left|N(\ep(f(X)+a\varphi(X)))-N(\ep(\grave{f}(X)+a\varphi(X)))\right|^2 \dX \\
\\
&\le \norm{\partial_X[N]}_{L^{\infty}(\I(f,\grave{f}))}^2\int_{\R} \cosh^2(qX)|f(X)-\grave{f}(X)|^2 \dX  \\
\\
&\le \norm{N}_{W^{r+1,\infty}(\I(f,\grave{f}))}^2\bignorm{f-\grave{f}}_{r,q}^2.
\end{aligned}
\end{equation}
So, $\norm{\Delta_1}_{L^2}$ has the bound
\begin{equation}\label{Delta1 bound}
\norm{\Delta_1}_{L^2} \le \M_0[f,\grave{f}]\bignorm{f-\grave{f}}_{r,q}.
\end{equation}

To estimate $\Delta_2$, we first rewrite
\begin{align*}
\Delta_2
&= \sum_{k=1}^r \partial_X^k[N](\ep({f}+a\varphi))\sum_{\sigmab \in \Sigma_k^r} C_{\sigmab}\prod_{j=1}^k \partial_X^{\sigma_j}[\ep(f+a\varphi)] \\
\\
&- \sum_{k=1}^r \partial_X^k[N](\ep(\grave{f}+a\varphi))\sum_{\sigmab \in \Sigma_k^r} C_{\sigmab}\prod_{j=1}^k \partial_X^{\sigma_j}[\ep(\grave{f}+a\varphi)] \\
\\
&= \sum_{k=1}^r \bunderbrace{\left(\partial_X^k[N](\ep(f+a\varphi))-\partial_X^k[N](\ep(\grave{f}+a\varphi))\right)}{\Delta_3(k)}\sum_{\sigmab \in \Sigma_k^r} \bunderbrace{\prod_{j=1}^k \partial_X^{\sigma_j}[\ep(f+a\varphi)]}{\Pi(\sigmab)} \\
\\
&+ \sum_{k=1}^r \partial_X^k[N](\ep(\grave{f}+a\varphi))\sum_{\sigmab \in \Sigma_k^r} \bunderbrace{\left(\prod_{j=1}^k \partial_X^{\sigma_j}[\ep(f+a\varphi)] - \prod_{j=1}^k\partial_X^{\sigma_j}[\ep(\grave{f}+a\varphi)]\right)}{\Delta_4(\sigmab)}.
\end{align*}

So, now we need to estimate
\begin{equation}\label{Delta2 first pass}
\norm{\Delta_2}_{L^2}
\le C_r\sum_{k=1}^r \sum_{\sigmab \in \Sigma_k^r} \bunderbrace{\norm{\cosh(q\cdot)\Delta_3(k)\Pi(\sigmab)}_{L^2}}{\TT_1(\sigmab,k)} + \bunderbrace{\norm{\cosh(q\cdot)\partial_X^k[N](\ep(f+a\varphi))\Delta_4(\sigmab)}_{L^2}}{\TT_2(\sigmab,k)}.
\end{equation}
When $k=1$, $\sigmab = \{r\}$, and $\TT_(\{r\},1)$ reduces to
\begin{align*}
\norm{\cosh(q\cdot)\Delta_3(1)\Pi_1(r)} 
&= \bignorm{\cosh(q\cdot)\left(\partial_X[N](\ep(f+a\varphi))-\partial_X[N](\ep(f+a\varphi))\right)\partial_X^r[\ep(f+a\varphi)]}_{L^2} \\
\\
&\le \bignorm{\left(\partial_X[N](\ep(f+a\varphi))-\partial_X[N](\ep(f+a\varphi))\right)\cosh(q\cdot)\partial_X^r[f]}_{L^2} \\
\\
&+ \bignorm{\cosh(q\cdot)\left(\partial_X[N](\ep(f+a\varphi))-\partial_X[N](\ep(f+a\varphi))\right)\partial_X^r[\varphi]}_{L^2} \\
\\
&\le \bignorm{\partial_X[N](\ep(f+a\varphi))-\partial_X[N](\ep(f+a\varphi))}_{L^{\infty}}\norm{\cosh(q\cdot)\partial_X^r[f]}_{L^2} \\
\\
&+ \bignorm{\cosh(q\cdot)\left(\partial_X[N](\ep(f+a\varphi))-\partial_X[N](\ep(f+a\varphi))\right)}_{L^2}\norm{\partial_X^r[\varphi]}_{L^{\infty}} \\
\\
&\le \bignorm{N}_{W^{r+1,\infty}(\I(f,\grave{f}))}\bignorm{f-\grave{f}}_{r,q}\norm{f}_{r,q} + \norm{N}_{W^{r+1,\infty}(\I(f,\grave{f}))}\bignorm{f-\grave{f}}_{r,q}C_{\star}\ep^{-r} \\
\\
&\le \M_0[f,\grave{f}](1+|a|\ep^{1-r})\bignorm{f-\grave{f}}_{r,q}.
\end{align*}
Similarly, $\TT_2(\{r\},1)$ is bounded by
\begin{align*}
\bignorm{\cosh(q\cdot)\partial_X[N](\ep(f+a\varphi))\Delta_4(r)}_{L^2} 
&= \bignorm{\cosh(q\cdot)\partial_X[N](\ep(f+a\varphi))\left(\partial_X^r[\ep(f+a\varphi)]-\partial_X^r[\ep(f+a\varphi)]\right)}_{L^2} \\
\\
&\le \bignorm{\partial_X[N](\ep(f+a\varphi))}_{L^{\infty}}\bignorm{\cosh(q\cdot)\partial_X[f-\grave{f}]}_{L^2} \\
\\
&\le \M_0[f,\grave{f}](1+|a|\ep^{1-r})\bignorm{f-\grave{f}}_{r,q}.
\end{align*}

Now let $2 \le k \le r$. Since all of the derivatives $\partial_X^{\sigma_j}[f]$ are at most order $r-1$ when $\sigmab \in \Sigma_k^r$, they are $L^{\infty}$.  So we estimate
$$
\TT_1(\sigmab,k)
\le \bignorm{\cosh(q\cdot)\left(\partial_X^k[N](\ep(f+a\varphi))-\partial_X^k[N](\ep(f+a\varphi))\right)}_{L^2}\prod_{j=1}^k \norm{\partial_X^{\sigma_j}[\ep(f+a\varphi)]}_{L^{\infty}}
$$
We bound the $L^2$ factor above using \eqref{N wrinfty est}:
$$
\bignorm{\cosh(q\cdot)\left(\partial_X^k[N](\ep(f+a\varphi))-\partial_X^k[N](\ep(f+a\varphi))\right)}_{L^2}
\le \M_0[f,\grave{f}]\bignorm{f-\grave{f}}_{r,q}.
$$
And we bound the product using \textit{exactly} the same reasoning that led to the estimate \eqref{product bound for mapping} in the proof of Proposition \ref{mapping}.  That is, we find
$$
\prod_{j=1}^k \norm{\partial_X^{\sigma_j}[\ep(f+a\varphi)]}_{L^{\infty}} 
\le \M_1[f]\left(1+|a|\ep^{1-r}\right),
$$
where $\M_1 \colon \R_+ \to \R_+$ is a continuous function that depends on $r$, $q$, and $C_{\star}$, but not on $a$, $\ep$, $\varphi$, or $f$.  Thus
\begin{equation}\label{T1 bound}
\TT_1(\sigmab,k) 
\le \M_0[f,\grave{f}]\M_1[f]\left(1+|a|\ep^{1-r}\right)\bignorm{f-\grave{f}}_{r,q}
\end{equation}
for $2 \le k \le r$ and $\sigmab \in \Sigma_k^r$.

Next, a first glance at $\TT_2(\sigmab,k)$ shows
$$
\TT_2(\sigmab,k)
\le \M_0[f,\grave{f}]\norm{\Delta_4(\sigmab)}_{L^2}.
$$
We use Lemma \ref{lipschitz product lemma} to rewrite
\begin{equation}\label{T2eq1}
\begin{aligned}
\Delta_4(\sigmab) 
&= \prod_{j=1}^k \partial_X^{\sigma_j}[\ep(f+a\varphi)] - \prod_{j=1}^k\partial_X^{\sigma_j}[\ep(f+a\varphi)] \\
\\
&= \sum_{j=1}^k \bunderbrace{\left(\partial_X^{\sigma_j}[\ep(f+a\varphi)]-\partial_X^{\sigma_j}[\ep(f+a\varphi)]\right)}{\Delta_5(\sigma_j)}
\bunderbrace{\prod_{\ell=1}^{j-1}\partial_X^{\sigma_{\ell}}[\ep(f+a\varphi)]}{\Pi_{j-1}(\sigmab)}
\bunderbrace{\prod_{\ell=j+1}^k \partial_X^{\sigma_{\ell}}[\ep(f+a\varphi)]}{\Pi^{j+1}(\sigmab)}.
\end{aligned}
\end{equation}
We easily estimate
$$
\norm{\cosh(q\cdot)\Delta_5(\sigma_j)}_{L^2} 
\le \M_0[f,\grave{f}]\bignorm{f-\grave{f}}_{r,q}.
$$

The products $\Pi_{j-1}(\sigmab)$ and $\Pi^{j+1}(\sigmab)$ require a little more care.  Since the derivatives $\partial_X^{\sigma_j}[f]$ and $\partial_X^j[\grave{f}]$ are still order at most $r-1$, these products are $L^{\infty}$.  More precisely, we can carefully replicate the steps that led to the estimate \eqref{product bound for mapping} to find
$$
\norm{\Pi_{j-1}(\sigmab)}_{L^{\infty}}
\le \M_2[f]\left(1+|a|\ep^{1-\sum_{\ell=1}^{j-1}\sigma_{\ell}}\right)
$$
and
$$
\norm{\Pi^{j+1}(\sigmab)}_{L^{\infty}}
\le \M_3[\grave{f}]\left(1+|a|\ep^{1-\sum_{\ell = j+1}^k \sigma_j}\right).
$$
Multiplying these estimates and using, as always, the assumptions $0 < \ep < 1$ and $|a| < 1$ and the relation $r = |\sigmab| = \medsum_{j=1}^k \sigma_j$, we find
\begin{equation}\label{T2eq2}
\norm{\Pi_{j-1}(\sigmab)\Pi^{j+1}(\sigmab)}_{L^{\infty}}
\le \M_4[f,\grave{f}]\left(1+|a|\ep^{1-r}\right).
\end{equation}
Then \eqref{T2eq1} and \eqref{T2eq2} together yield
\begin{equation}\label{T2 bound}
\TT_2(\sigmab,k)
\le r\M_0[f,\grave{f}]^2\M_3[f,\grave{f}]\left(1+|a|\ep^{1-r}\right)\bignorm{f-\grave{f}}_{r,q}.
\end{equation}

Now that we have bounded both $\TT_1(\sigmab,k)$ and $\TT_2(\sigmab_k,k)$ for $k = 2,\ldots,r$ and $\sigmab \in \Sigma_k^r$, we can use \eqref{T1 bound} and \eqref{T2 bound} to our estimate \eqref{Delta2 first pass} for $\Delta_2$ and conclude
$$
\norm{\Delta_2}_{L^2} 
\le \M_5[f,\grave{f}]\left(1+|a|\ep^{1-r}\right)\bignorm{f-\grave{f}}_{r,q}.
$$
for some continuous map $\M_5 \colon \R_+^2 \to \R_+$ that is independent of $\ep$, $|a|$, and $\varphi$ (but dependent on $r$ and $q$).  This, together with \eqref{Delta1 bound}, gives \eqref{general hrq lipschitz}.
\end{proof}

By taking $\varphi = 0$ and $\ep = a = 1/2$ in Proposition \ref{general hrq lipschitz prop}, we obtain the following simpler estimate.

\begin{corollary}
Let $N \in \C^{r+1}(\R)$.  There exists a radially increasing map $\M \colon \R_+^4 \to \R_+$ such that 
$$
\bignorm{N(f)-N(\grave{f})}_{r,q} \le \M[\norm{f}_{r,q},\bignorm{\grave{f}}_{r,q},\norm{\varphi}_{W^{r,\infty}},\norm{\grave{\varphi}}_{W^{r,\infty}}]\bignorm{f-\grave{f}}_{r,q}
$$
for all $f,\grave{f} \in H_q^r$.
\end{corollary}

\subsubsection{A Lipschitz estimate in $H_q^1$}\label{H1q lipschitz section}

\begin{proposition}\label{easy peasy lipschitz}
Let $N \in \C^1(\R)$.  There exists a radially increasing map $\R_+^4 \to \R_+$ such that 
$$
\bignorm{(f-\grave{f})(g+\varphi)N(h+\grave{\varphi})}_{1,q}
\le \M[\norm{g}_{1,q},\norm{h}_{1,q},\norm{\varphi}_{W^{1,\infty}},\norm{\grave{\varphi}}_{W^{1,\infty}}]\bignorm{f-\grave{f}}_{1,q}
$$
for all $f$, $\grave{f}$, $g$, $h\in H_q^1$ and $\varphi$, $\grave{\varphi} \in W^{1,\infty}$.
\end{proposition}

\begin{proof}
This is a straightforward calculation that requires only one pass with the chain rule, so we omit the details.
\end{proof}

\subsubsection{Lipschitz estimates in $W^{1,\infty}$} \label{W1infty lipschitz section}
The estimates in this section are much simpler than the preceding $H_q^r$ mapping and Lipschitz estimates because we work only with $r=1$ and so do not need to keep careful track of powers of $\ep$. However, we do need a ``decay borrowing'' product estimate, which we take from Lemma A.2 in \cite{faver-wright}: for $r \ge 0$ and $q > 0$, we have

\begin{equation}\label{decay borrowing}
\norm{fg}_{r,q/2} 
\le C_r\norm{f}_{r,q}\norm{\sech\left(\frac{q}{2}\cdot\right)g}_{W^{r,\infty}}.
\end{equation}

\begin{proposition}\label{lipschitz on phi}
Let $N \in \C^2(\R)$.  There exists a radially increasing map $\M \colon \R_+^3 \to \R_+$ such that 
\begin{equation}\label{lipschitz on phi est}
\norm{f(N(g+\varphi)-N(g+\grave{\varphi}))}_{1,q/2}
\le \M[g,\varphi,\grave{\varphi}]\norm{\sech\left(\frac{q}{2}\cdot\right)(\varphi-\grave{\varphi})}_{W^{1,\infty}}\norm{f}_{1,q}
\end{equation}
for all $f$, $g \in H_q^1$ and $\varphi$, $\grave{\varphi} \in W^{1,\infty}$.
\end{proposition}

\begin{proof}
We have
\begin{align*}
\norm{f(N(g+\varphi)-N(g+\grave{\varphi}))}_{1,q/2} &= \bunderbrace{\norm{\cosh\left(\frac{q}{2}\cdot\right)f(N(g+\varphi)-N(g+\grave{\varphi}))}_{L^2}}{\Delta_1} \\
\\
&+ \bunderbrace{\norm{\cosh\left(\frac{q}{2}\cdot\right)\partial_X[f(N(g+\varphi)-N(g+\grave{\varphi}))]}_{L^2}}{\Delta_2}.
\end{align*}

The first estimate for $\Delta_1^2$ is similar to that for $\Delta_1^2$ in the proof of Proposition \ref{general hrq lipschitz prop}, i.e., a direct calculation with the integral yields
$$
\Delta_1
\le \norm{N}_{W^{2,\infty}(\I(g,\varphi,\grave{\varphi}))}\norm{\cosh\left(\frac{q}{2}\cdot\right)f(\varphi-\grave{\varphi})}_{L^2}.
$$
Next, we use the decay-borrowing estimate \eqref{decay borrowing}:
\begin{align*}
\Delta_1 &\le \norm{N}_{W^{2,\infty}(\I(g,\varphi,\grave{\varphi}))}\norm{f(\varphi-\grave{\varphi})}_{0,q/2} \\
\\
&\le C\norm{N}_{W^{2,\infty}(\I(g,\varphi,\grave{\varphi}))}\norm{f}_{0,q}\norm{\sech\left(\left|q-\frac{q}{2}\right|\cdot\right)(\varphi-\grave{\varphi})}_{W^{0,\infty}}\\
\\
&\le C\norm{N}_{W^{2,\infty}(\I(g,\varphi,\grave{\varphi}))}\norm{f}_{1,q}\norm{\sech\left(\frac{q}{2}\cdot\right)(\varphi-\grave{\varphi})}_{W^{1,\infty}}.
\end{align*}
This is exactly the kind of estimate we want for $\Delta_1$, and so we move on to $\Delta_2$.  Since we only ever take one derivative in this proof, we write $N' = \partial_X[N]$, etc. We have
%
%
%
\begin{align*}
\Delta_2 &\le \bunderbrace{\norm{\cosh\left(\frac{q}{2}\cdot\right)f'(N(g+\varphi)-N(g+\grave{\varphi}))}_{L^2}}{\Delta_3} \\
\\
&+\bunderbrace{\norm{\cosh\left(\frac{q}{2}\cdot\right)f(N'(g+\varphi)(g'+\varphi')-N'(g+\grave{\varphi})(g'+\grave{\varphi}'))}_{L^2}}{\Delta_4}.
\end{align*}
We get an estimate on $\Delta_3$ of the form in \eqref{lipschitz on phi est} by working directly with the integral, and we have

\begin{align*}
\Delta_4 &\le \bunderbrace{\norm{\cosh\left(\frac{q}{2}\cdot\right)fN'(g+\varphi)(\varphi'-\grave{\varphi}')}_{L^2}}{\Delta_{5}} 
+ \bunderbrace{\norm{\cosh\left(\frac{q}{2}\cdot\right)f(N'(g+\varphi)-N'(g+\grave{\varphi}))g'}_{L^2}}{\Delta_{6}} \\
\\
&+ \bunderbrace{\norm{\cosh\left(\frac{q}{2}\cdot\right)f(N'(g+\varphi)-N'(g+\grave{\varphi}))\grave{\varphi}'}_{L^2}}{\Delta_{7}}.
\end{align*}

For $\Delta_5$, we have
\begin{align*}
\Delta_{5} &= \norm{fN'(g+\varphi)(\varphi'-\grave{\varphi}')}_{0,q/2} \\
\\
&\le \norm{N'(g+\varphi)}_{L^{\infty}}\norm{f(\varphi'-\grave{\varphi}')}_{0,q/2} \\
\\
&\le C\norm{N}_{W^{2,\infty}(\I(g,\varphi,\grave{\varphi}))}\norm{f}_{0,q}\norm{\sech\left(\frac{q}{2}\cdot\right)(\varphi'-\grave{\varphi}')}_{W^{1,\infty}}.
\end{align*}

For $\Delta_6$, we have
\begin{align*}
\Delta_{6} &\le \norm{\cosh(q\cdot)f}_{L^{\infty}}\norm{(N'(g+\varphi)-N'(g+\grave{\varphi}))g'}_{L^2} \\
\\
&= \norm{\cosh(q\cdot)f}_{L^{\infty}}\norm{(N'(g+\varphi)-N'(g+\grave{\varphi}))g'}_{0,0} \\
\\
&\le C_q\norm{f}_{1,q}\norm{g'}_{0,q/2}\norm{\sech\left(\frac{q}{2}\cdot\right)(N'(g+\varphi)-N'(g+\grave{\varphi}))}_{W^{0,\infty}} \\
\\
&\le C_q\norm{f}_{1,q}\norm{g}_{1,q}\norm{N}_{W^{2,\infty}(\I(g,\varphi,\grave{\varphi}))}\norm{\sech\left(\frac{q}{2}\cdot\right)(\varphi-\grave{\varphi})}_{W^{1,\infty}}.
\end{align*}
For the last Lipschitz estimate on $N'$, we just work pointwise to get the estimate in $W^{0,\infty} = L^{\infty}$.

Finally, for $\Delta_7$, we use the same techniques as above to produce
\begin{align*}
\Delta_7 &\le \norm{\grave{\varphi}'}_{L^{\infty}}\norm{f(N'(g+\varphi)-N'(g+\grave{\varphi}))}_{0,q/2} \\
\\
&\le \norm{\grave{\varphi}}_{W^{1,\infty}}\norm{f}_{0,q}\norm{\sech\left(\frac{q}{2}\cdot\right)(N'(g+\varphi)-N'(g+\grave{\varphi}))}_{W^{0,\infty}} \\
\\
&\le \norm{\grave{\varphi}}_{W^{1,\infty}}\norm{f}_{1,q}\norm{N}_{W^{2,\infty}(\I(g,\varphi,\grave{\varphi}))}\norm{\sech\left(\frac{q}{2}\cdot\right)(\varphi-\grave{\varphi})}_{W^{1,\infty}}.
\end{align*}
Combining the estimates on $\Delta_1$ through $\Delta_7$, we have \eqref{lipschitz on phi est}.
\end{proof}

\begin{proposition}\label{use the ftc, luke} 
Let $N \in \C^3(\R)$.  There exists a radially increasing map $\M \colon \R_+^3 \to \R_+$ such that 
$$
\norm{(N(f+\varphi)-N(\varphi)) - (N(f+\grave{\varphi})-N(\grave{\varphi}))}_{1,q/2}
\le \M[f,\varphi,\grave{\varphi}]\norm{\sech\left(\frac{q}{2}\cdot\right)(\varphi-\grave{\varphi})}_{W^{1,\infty}}\norm{f}_{1,q}.
$$
\end{proposition}

\begin{proof}
We use the fundamental theorem of calculus twice to rewrite 
$$
(N(f+\varphi)-N(\varphi))-(N(f+\grave{\varphi})-N(\grave{\varphi}))
= f(\varphi-\grave{\varphi})\int_0^1\int_0^1 N''(tf+\grave{\varphi}+s(\varphi-\grave{\varphi})) \ds \dt.
$$
The necessary estimates are then similar to those in Proposition \ref{lipschitz on phi}.  Of these estimates, arguably the most complicated involves controlling the $L^2$-norm of
\begin{equation}\label{arguably hard}
\cosh\left(\frac{q}{2}\cdot\right)f(\varphi-\grave{\varphi})\partial_X\left[\int_0^1\int_0^1 N''(tf+\grave{\varphi}+s(\varphi-\grave{\varphi})) \ds \dt\right].
\end{equation}
After differentiating under the integral in \eqref{arguably hard}, the terms that we need to bound are by now routine; we remove a number of factors in the $L^{\infty}$-norm and then use decay borrowing on the rest.  For example, one of these terms is
$$
\Delta_{\text{hard}} =\cosh\left(\frac{q}{2}\cdot\right)f(\varphi-\grave{\varphi})(\varphi'-\grave{\varphi}')\bunderbrace{\int_0^1\int_0^1 |N'''(tf+\grave{\varphi}+s(\varphi-\grave{\varphi}))| \ds \dt}{\I},
$$
and we have
\begin{align*}
\norm{\Delta_{\text{hard}}}_{L^2} &\le \norm{\varphi-\grave{\varphi}}_{L^{\infty}}\norm{\I}_{L^{\infty}}\norm{\cosh\left(\frac{q}{2}\cdot\right)f(\varphi'-\grave{\varphi}')}_{L^2} \\
\\
&\le \left(\norm{\varphi}_{L^{\infty}}+\norm{\grave{\varphi}}_{L^{\infty}}\right)\norm{\I}_{L^{\infty}}\norm{f}_{1,q}\norm{\sech\left(\frac{q}{2}\cdot\right)(\varphi'-\grave{\varphi}')}_{W^{0,\infty}}. 
\end{align*}
We omit the other details, as they are by now routine.
%
%
%
\end{proof}

\subsubsection{Estimates in $H_{\per}^r$}  All of our estimates so far have involved the space $H_q^r$.  That is, we have proved estimates for functions that are square-integrable on all of $\R$.  However, none of our proofs relied in an essential way on the domain of integration being $\R$, and so we can replace $\R$ with $[0,2\pi]$ and find that our proofs are still valid for the space $H_{\per}^r$.  Specifically, by taking $q = a = 0$, $\varphi = 0$, and $\ep = 1/2$, we can rerun the proofs of Propositions \ref{mapping prop} and \ref{general hrq lipschitz prop} in $H_{\per}^r$ to obtain the following.

\begin{proposition}\label{Hrper mapping}
Let $r \ge 1$ and $N \in \C^r([0,2\pi])$.  There exists an increasing map $\M \colon \R_+ \to \R_+$ such that 
$$
\norm{N(f)}_{H_{\per}^r} 
\le \M[\norm{f}_{H_{\per}^r}]
$$
for all $f \in H_{\per}^r$.
\end{proposition}

\begin{proposition}\label{Hrper lipschitz}
Let $r \ge 1$ and $N \in \C^{r+1}([0,2\pi])$.  There exists a radially increasing map $\M \colon \R_+^2 \to \R_+$ such that 
$$
\bignorm{N(f)-N(\grave{f})}_{H_{\per}^r} 
\le \M[\norm{f}_{H_{\per}^r}, \bignorm{\grave{f}}_{H_{\per}^r}]\bignorm{f-\grave{f}}_{H_{\per}^r}.
$$
\end{proposition}

\subsubsection{Auxiliary identities for sums and products}

\begin{lemma}\label{product of sums lemma}
Let $\{A_j\}_{j=1}^r,\{B_j\}_{j=1}^r \subseteq \R$.  Then
$$
\prod_{j=1}^r (A_j + B_j) = \sum_{\substack{\alphab,\betab \in \{0,1\}^r \\ \alpha_j + \beta_j = 1, \ j=1,\ldots,r}} \prod_{i=1}^r A_i^{\alpha_i}B_i^{\beta_i}
$$
\end{lemma}

\begin{proof}
We induct on $r$.  If $r=1$, then
$$
\sum_{\substack{\alpha_1,\beta_1 \in \{0,1\} \\ \alpha_1+\beta_1 = 1}} A_1^{\alpha_1}B_1^{\alpha_1}
= A_1^1B_1^0 + A_1^0B_1^1 
= A_1 + B_1.
$$
Suppose the formula holds for some $r \ge 1$.  Then
\begin{align*}
\prod_{j=1}^{r+1} (A_j+B_j) 
&= (A_{r+1}+B_{r+1})\prod_{j=1}^r (A_j+B_j)^r \\
\\
&= (A_{r+1}+B_{r+1})\sum_{\substack{\alphab,\betab \in \{0,1\}^r \\ \alpha_j + \beta_j = 1, \ j=1,\ldots,r}} \prod_{i=1}^r A_i^{\alpha_i}B_i^{\beta_i} \\
\\
&= \sum_{\substack{\alphab,\betab \in \{0,1\}^r \\ \alpha_j + \beta_j = 1, \ j=1,\ldots,r}} \prod_{i=1}^r A_i^{\alpha_i}B_i^{\beta_i}A_{r+1} 
+ \sum_{\substack{\alphab,\betab \in \{0,1\}^r \\ \alpha_j + \beta_j = 1, \ j=1,\ldots,r}} \prod_{i=1}^r A_i^{\alpha_i}B_i^{\beta_i}B_{r+1} \\
\\
&= \sum_{\substack{\alphab,\betab \in \{0,1\}^{r+1} \\ \alpha_j + \beta_j = 1, \ j = 1,\ldots,r+1}} \prod_{i=1}^{r+1} A_i^{\alpha_i}B_i^{\beta_i}. \qedhere
\end{align*}
\end{proof}

\begin{lemma}\label{lipschitz product lemma}
Let $z_1,\ldots,z_n,\grave{z}_1,\ldots,\grave{z}_n \in \mathbb{C}$.  Then 
$$
\prod_{k=1}^n z_k-\prod_{k=1}^n \grave{z}_k = \sum_{k=1}^n (z_k-\grave{z}_k)\prod_{j=1}^{k-1} \grave{z}_j\prod_{j=k+1}^n z_j.
$$
\end{lemma}

\begin{proof}
When $n=1$ it is obvious, so assume it holds zor some $n$ and consider the $n+1$ case:
\begin{align*}
\prod_{k=1}^{n+1} z_k-\prod_{k=1}^{n+1} \grave{z}_k &= z_{n+1}\prod_{k=1}^nz_k - \grave{z}_{n+1}\prod_{k=1}^n \grave{z}_k \\
\\
&= z_{n+1}\prod_{k=1}^n z_k - z_{n+1}\prod_{k=1}^n \grave{z}_k + z_{n+1}\prod_{k=1}^n\grave{z}_k - \grave{z}_{k+1}\prod_{k=1}^n \grave{z}_k \\
\\
&= z_{n+1}\left(\prod_{k=1}^n z_k - \prod_{k=1}^n \grave{z}_k\right) + (z_{n+1}-\grave{z}_{n+1})\prod_{k=1}^n \grave{z}_k \\
\\
&= z_{n+1}\sum_{k=1}^n (z_k-\grave{z}_k)\left(\prod_{j=1}^{k-1}\grave{z}_j\right)\left(\prod_{j=k+1}^n z_j\right) + (z_{n+1}-\grave{z}_{n+1})\prod_{k=1}^n \grave{z}_k \\
\\
&= \sum_{k=1}^n (z_k-\grave{z}_k)\left(\prod_{j=1}^{k-1}\grave{z}_j\right)\left(\prod_{j=k+1}^{n+1}z_j\right) + (z_{n+1}-\grave{z}_{n+1})\prod_{k=1}^n \grave{z}_k \\
\\
&= \sum_{k=1}^{n+1} (z_k-\grave{z}_k)\left(\prod_{j=1}^{k-1}\grave{z}_k\right)\left(\prod_{j=k+1}^{n+1}z_j\right). \qedhere
\end{align*}
\end{proof}

%% file: appendix_existence_of_periodic_solutions.tex
\section{Existence of Periodic Solutions}
\label{periodic solutions appendix}

\subsection{Conversion to a fixed-point problem}\label{fixed point conversion appendix}
We begin by introducing a periodic profile and a frequency scaling: let $\thetab(X) = \phib(\omega{X})$, where $\phib = \phib(Y)$ is $2\pi$-periodic and $\omega \in \R$.  Then the problem of \eqref{lw eqns}, $\Thetab_{\ep}(\phib(\omega\cdot)) = 0$, converts to 
\begin{multline}\label{matrix periodic system}
\Phib_{\ep}(\phib,\omega) := 
\begin{bmatrix*}
1 &0 \\
0 &\ep^2\omega^2(c_{\kappa}^2+\ep^2)\partial_Y^2 + \lambda_+^{\ep\omega}
\end{bmatrix*}
\phib
+\begin{bmatrix*}
\varpi^{\ep,\omega} &0 \\
0 &\ep^2\lambda_+^{\ep\omega}
\end{bmatrix*}B^{\ep\omega}(\phib,\phib) 
\\
+ \begin{bmatrix*}
\varpi^{\ep,\omega} &0 \\
0 &\ep^2\lambda_+^{\ep\omega}
\end{bmatrix*}\cb^{\ep\omega}(\phib,\phib,\phib)
=0,
\end{multline}
where $\varpi^{\ep,\omega}$ has the symbol
$$
\tilde{\varpi^{\ep,\omega}}(k) 
= \tvarpi^{\ep,\omega}(k) 
:= \tvarpi^{\ep}(\omega{k})
$$
with $\tvarpi^{\ep}$ defined in \eqref{tvarpi ep defn} and the other Fourier multipliers are defined per Section \ref{lw scaling section} and the scaling properties in Appendix \ref{fourier multipliers appendix}.

At this point, we could run a standard ``bifurcation from a simple eigenvalue'' argument on $\Phib_{\ep}$ by noting that $\Phib_{\ep}(0,\omega) = 0$ for all $\omega$ and calculating, with some labor, that the linearization of $\Phib_{\ep}$ at $(0,\omega_{\ep})$ is Fredholm with index 1.  In particular, we noted back in Section \ref{periodic solutions section} that $D_{\phib}\Phib_{\ep}(0,\omega_{\ep})\cos(\cdot)\jb = 0$.  Then thanks to our restriction to ``even $\times$ even'' functions and the Friesecke-Pego cancelation, the kernel of $D_{\phib}\Phib_{\ep}(0,\omega_{\ep})$ is in fact spanned by $\cos(\cdot)\jb$.  Verification of the bifurcation condition reduces, after much computation, to the inequality \eqref{CRZ bifurcation condition}.

However, since we seek uniform estimates in $\ep$, which the Crandall-Rabinowitz-Zeidler approach will not ostensibly provide, we instead follow the proof of Crandall and Rabinowitz and rewrite the problem \eqref{matrix periodic system} in fixed point form.  Our presentation is much the same as in \cite{faver-wright} with various technical adaptations to handle the higher order terms from $\cb^{\ep}$.  

With $\nub = \cos(\cdot)\jb$ as in Theorem \ref{periodic main theorem},  let
$$
\mathcal{Z} = \set{\psib \in E_{\per}^2 \times E_{\per}^2}{\hat{\psi}_1(\pm1) = 0} = \{\nub\}^{\perp},
$$
where
$$
E_{\per}^r = \set{f \in H_{\per}^r}{f \text{ is even}}.
$$
We recall the definitions and properties of periodic Sobolev spaces from Appendix \ref{periodic sobolev spaces appendix}.  In other words, $\mathcal{Z} = \{\nub\}^{\perp}$, the orthogonal complement of $\{\nub\}$ in $E_{\per}^2 \times E_{\per}^2$.

Then with the ansatz
$$
\phib = a\nub+a\psib, \ \psib \in \mathcal{Z}
\quadword{and} \omega = \omega_{\ep}+t,
$$
the system \eqref{matrix periodic system} becomes
\begin{multline}\label{matrix periodic system ta}
\begin{bmatrix*}
1 &0 \\
0 &\ep^2(\omega_{\ep}+t)^2\cep^2\partial_Y^2 + \lambda_+^{\ep(\omega_{\ep}+t)} 
\end{bmatrix*}
(\nub+\psib)
+ 
a\begin{bmatrix*}
\varpi_{\cep}^{\ep,\omega_{\ep}+t}  &0 \\
0 & \ep^2\lambda_+^{\ep(\omega_{\ep}+t)}
\end{bmatrix*}
\B^{\ep}(\psib,t)
\\
+
a
\begin{bmatrix*}
\varpi_{\cep}^{\ep,\omega_{\ep}+t}  &0 \\
0 & \ep^2\lambda_+^{\ep(\omega_{\ep}+t)}
\end{bmatrix*}
\xtr^{\ep}(\psib,t,a)  = 0.
\end{multline}
Here we have abbreviated
$$
\B^{\ep}(\psib,t) =
\begin{pmatrix*}
\B_1^{\ep}(\psib,t) \\
\B_2^{\ep}(\psib,t)
\end{pmatrix*}
:= B^{\ep(\omega_{\ep}+t)}(\nub+\psib,\nub+\psib)
$$
and
$$
\xtr^{\ep}(\psib,t,a) 
=
\begin{pmatrix*}
\xtr_1^{\ep}(\psib,t,a) \\
\xtr_2^{\ep}(\psib,t,a)
\end{pmatrix*}
:= a\ep^2J_1^{\ep(\omega_{\ep}+t)}M_{1/\kappa}\left[\left(J^{\ep(\omega_{\ep}+t)}(\nub+\psib)\right)^{.3}.N(a\ep^2J^{\ep(\omega_{\ep}+t)}(\nub+\psib))\right].
$$

Let $\Pi_1$ be the multiplier with symbol 
$$
\tPi_1(k) := \delta_{|k|,1} = \begin{cases}
1, &|k| = 1 \\
0, &|k| \ne 1 
\end{cases}
$$
and let $\Pi_2 := \ind-\Pi_1$. Then
\begin{equation}\label{Pi1 Pi2 props}
\Pi_1\cos(\cdot) = \cos(\cdot),
\qquad
\Pi_2\cos(\cdot) = 0,
\qquad
\Pi_1\psi = 0,
\quadword{and}
\Pi_2\psi 
\end{equation}
for any $\psi \in E_{\per}^2$ with $\hat{\psi}(1) = 0$.

Let $\xi_c$ be the multiplier with symbol 
$$
\txi_c(k) := -c^2k^2+\tlambda_+(k)
$$
and let 
$$
\xi^{\ep,t} := \xi_{\cep}^{\ep(\omega_{\ep}+t)}
= \ep^2(\omega_{\ep}+t)^2\cep^2\partial_Y^2+\lambda_+^{\ep(\omega_{\ep}+t)},
$$
so  $\xi^{\ep,t}$ has the symbol 
$$
\tilde{\xi^{\ep,t}}(k) = \txi_{\cep}(\ep(\omega_{\ep}+t)k).
$$

After we apply $\Pi_1$ and $\Pi_2$ to the first component of \eqref{matrix periodic system ta} and use \eqref{Pi1 Pi2 props}, we see that \eqref{matrix periodic system ta} is equivalent to the three equations
\begin{equation}\label{prefp1}
\psi_1 + a\varpi^{\ep,\omega_{\ep}+t}\left(\B_1^{\ep}(\psib,t)+\xtr_1^{\ep}(\psib,t,a)\right) = 0.
\end{equation}
\begin{equation}\label{prefp2}
\xi^{\ep,t}\psi_2+a\ep^2\Pi_2\lambda_+^{\ep(\omega_{\ep}+t)}\left(\B_2^{\ep}(\psib,t)+\xtr_2^{\ep}(\psib,t,a)\right) = 0,
\end{equation}
and
\begin{equation}\label{prefp3}
xi^{\ep,t}\cos(\cdot)
+a\ep^2\Pi_1\lambda_+^{\ep(\omega_{\ep}+t)}\left(\B_2^{\ep}(\psib,t)+\xtr_2^{\ep}(\psib,t,a)\right) = 0.
\end{equation}
The first equation, \eqref{prefp1}, immediately converts to the fixed-point form
\begin{equation}\label{fp1}\psi_2 
= -a\varpi^{\ep,\omega_{\ep}+t}\left(\B_1^{\ep}(\psib,t)+\xtr_1^{\ep}(\psib,t,a)\right) \\
=:\Psi_1^{\ep}(\psib,t,a).
\end{equation}
The invertibility of $\xi^{\ep,t}$ on the range of $\Pi_2$, as detailed in Proposition \ref{master periodic lemma}, means that \eqref{prefp2} is equivalent to
\begin{equation}\label{fp2}
\psi_2
=-a\ep^2\left(\xi^{\ep,t}\right)^{-1}\Pi_2\lambda_+^{\ep(\omega_{\ep}+t)}\left(\B_2^{\ep}(\psib,t)+\xtr_2^{\ep}(\psib,t,a)\right) 
=:\Psi_2^{\ep}(\psib,t,a).
\end{equation}
Last, if $\psib$ is even, then \eqref{prefp3} holds if and only if the Fourier transform of its left side evaluated at $k=1$ is equal to zero.  We use \eqref{xi symbol} below to write
$$
\ft[\xi^{\ep,t}\cos(\cdot)](1) = \frac{\txi_{\cep}(\ep\omega_{\ep}+\ep{t})}{2} = \frac{(\ep{t})\Upsilon_{\ep}}{2}+ \frac{(\ep{t})^2\rhs_{\ep}(\ep{t})}{2}.
$$
Since $\Upsilon_{\ep}$ is bounded away from zero, we conclude that \eqref{prefp1} is equivalent to
\begin{equation}\label{fp3}
t = -\frac{\ep}{\Upsilon_{\ep}}\rhs_{\ep}(\ep{t})t^2-\frac{2\ep{a}}{\Upsilon_{\ep}}\ft\left[\lambda_+^{\ep(\omega_{\ep}+t)}\left(\B_2^{\ep}(\psib,t)+\xtr_2^{\ep}(\psib,t,a)\right)\right]\!(1)
=:\Psi_3^{\ep}(\psib,t,a).
\end{equation}

Now we are ready to pose our fixed point problem.  Let 
$$
\W^r 
= \left(E_{\per}^r \times E_{\per}^r\right) \cap \mathcal{Z}
\quadword{and}
\norm{\psib}_r 
:= \norm{\psib}_{\W^r} 
= \norm{\psi_1}_{H_{\per}^r}+\norm{\psi_2}_{H_{\per}^r}
$$
We will find an interval $[-a_{\per},a_{\per}] \subseteq \R$ and maps 
$$
[-a_{\per},a_{\per}] \to \W^2 \colon a \mapsto \psib_{\ep}^a = (\psi_{\ep,1}^a,\psi_{\ep,2}^a)
\quadword{and}
[-a_{\per},a_{\per}] \to \R \colon t \mapsto t_{\ep}^a
$$
such that 
$$
\begin{pmatrix*}
\psi_{1,\ep}^a \\
\psi_{2,\ep}^a \\
t_{\ep}^a
\end{pmatrix*}
=
\begin{pmatrix*}
\Psi_1^{\ep}(\psib_{\ep}^a,t_{\ep}^a,a) \\
\Psi_2^{\ep}(\psib_{\ep}^a,t_{\ep}^a,a) \\
\Psi_3^{\ep}(\psib_{\ep}^a,t_{\ep}^a,a)
\end{pmatrix*}
=:\Psib^{\ep}(\psib_{\ep}^a,t_{\ep}^a,a)
$$
for all $\ep$ in the interval $(0,\ep_{\per})$, where $\ep_{\per}$ comes from Proposition \ref{master periodic lemma}.  Once we show the existence of these maps, we will prove the additional properties and estimates in Theorem \ref{periodic main theorem}.

\subsection{The fixed point lemma}
We will use the following lemma to solve the fixed point problem and obtain the various estimates and smoothness properties.  

\begin{lemma}\label{fp}
Let $\X$ be a Banach space and and let $F_{\ep} \colon \X \times \R \to \X, \ 0 < \ep < \ep_0$, be a family of maps with the following properties: there exist continuous maps $\M_{\map} \colon \R_+ \to \R_+$, $\M_{\lip} \colon \R_+^2 \to \R_+$, and $\M_{\max} \colon \R_+ \to \R_+$ and a constant $a_1 > 0$ such that if $|a|$, $|\grave{a}| \le a_1$, then 
\begin{equation}\label{fp lemma: bound}
\sup_{0 < \ep < \ep_0} \norm{F_{\ep}(x,a)} \le \M_{\map}[\norm{x}]\left(|a| + \norm{x}^2\right),
\end{equation}
\begin{equation}\label{fp lemma: lipschitz vector}
\sup_{0 < \ep < \ep_0} \norm{F_{\ep}(x,a)-F_{\ep}(\grave{x},a)} \le \M_{\lip}[\norm{x},\norm{\grave{x}}]\left(|a|+\norm{x}+\norm{\grave{x}}\right)\norm{x-\grave{x}},
\end{equation}
and
\begin{equation}\label{fp lemma: lipschitz a}
\sup_{0 < \ep < \ep_0} \norm{F_{\ep}(x,a)-F_{\ep}(x,\grave{a})} \le \M_{\max}[\norm{x}]|a-\grave{a}|
\end{equation}
for any $x,\grave{x} \in \X$.
Then there are constants $a_0$, $r_0$ such that if $|a| \le a_0$ and $0 < \ep < \ep_0$, there exists a unique $x_{\ep}^a \in \mathfrak{B}(r_0) := \set{x \in \X}{\norm{x} \le r_0}$ such that 
\begin{equation}\label{fp periodic equation}
\norm{x_{\ep}^a} \le r_0 \quadword{and} F_{\ep}(x_{\ep}^a,a) = x_{\ep}^a.
\end{equation}
Moreover, there is a constant $C> 0$ such that if $|a|$, $|\grave{a}| \le a_0$, then 
\begin{equation}\label{fp periodic lipschitz}
\sup_{0 < \ep < \ep_0} \norm{x_{\ep}^a - x_{\ep}^{\grave{a}}} \le C|a-\grave{a}|.
\end{equation}
\end{lemma}

\begin{proof}
Let
$$
M = \max\left\{\max_{0 \le t \le 1} \M_{\map}[t], \max_{0 \le t,\grave{t} \le 1} \M_{\lip}[t,\grave{t}], \max_{0 \le t \le 1} \M_{\max}[t]\right\}.
$$
Set 
$$
r_0 = \min\left\{\frac{1}{6M},1\right\} 
\quadword{and}
a_0 = \min\left\{\frac{r_0}{2M},a_1,\frac{1}{6M}\right\}.
$$
Observe that 
$$
r_0 \le \frac{1}{2M} \Longrightarrow r_0^2 \le \frac{r_0}{2M} \Longrightarrow Mr_0^2 \le \frac{r_0}{2}.
$$
Then for $x \in \mathfrak{B}(r_0)$ and $|a| \le a_0$, we have
\begin{equation}\label{fp lemma pf: map}
\begin{aligned}
\norm{F_{\ep}(x,a)} &\le \M_{\map}[\norm{x}]\left(|a|+\norm{x}^2\right) \\
\\
&\le M(a_0+r_0^2) \\
\\
&\le M\left(\frac{r_0}{2M}\right) + Mr_0^2 \\
\\
&\le \frac{r_0}{2} + \frac{r_0}{2} \\
\\
&= r_0.
\end{aligned}
\end{equation}

Next, for $x$, $\grave{x} \in \mathfrak{B}(r_0)$ and $|a| \le a_0$, we compute
\begin{equation}\label{fp lemma pf: lipschitz}
\begin{aligned}
\norm{F_{\ep}(x,a)-F_{\ep}(\grave{x},\grave{a})}
&\le \M_{\lip}[\norm{x},\norm{\grave{x}}]\left(|a|+\norm{x}+\norm{\grave{x}}\right)\norm{x-\grave{x}} \\
\\
&\le M(|a|+2r_0)\norm{x-\grave{x}} \\
\\
&\le M\left(\frac{1}{6M}+\frac{2}{6M}\right)\norm{x-\grave{x}} \\
\\
&= \frac{1}{2}\norm{x-\grave{x}}.
\end{aligned}
\end{equation}
Together, \eqref{fp lemma pf: map} and \eqref{fp lemma pf: lipschitz} imply that each map $F_{\ep}(\cdot,a) \colon \mathfrak{B}(r_0) \to \mathfrak{B}(r_0)$ is a contraction on $\mathfrak{B}(r_0)$, which means there exists a unique $x_{\ep}^a \in \mathfrak{B}(r_0)$ such that $F_{\ep}(x_{\ep}^a,a) = x_{\ep}^a$.  

Finally, we have
\begin{align*}
\norm{x_{\ep}^a-x_{\ep}^{\grave{a}}}
&= \norm{F_{\ep}(x_{\ep}^a,a)-F_{\ep}(x_{\ep}^{\grave{a}},\grave{a})} \\
\\
&\le \norm{F_{\ep}(x_{\ep}^a,a)-F_{\ep}(x_{\ep}^a,\grave{a})} + \norm{F_{\ep}(x_{\ep}^a,\grave{a})-F_{\ep}(x_{\ep}^{\grave{a}},\grave{a})} \\
\\
&\le \M_{\max}[\norm{x_{\ep}^a}]|a-\grave{a}| + \M_{\lip}[\norm{x_{\ep}^a},\norm{x_{\ep}^{\grave{a}}}]|\grave{a}|\norm{x_{\ep}^a-x_{\ep}^{\grave{a}}} \\
\\
&\le M|a-\grave{a}| + Ma_0\norm{x_{\ep}^a-x_{\ep}^{\grave{a}}}.
\end{align*}
Since 
$$
Ma_0 \le M\left(\frac{1}{2M}\right) = \frac{1}{2},
$$
we can rearrange this last inequality to
$$
\frac{1}{2}\norm{x_{\ep}^a-x_{\ep}^{\grave{a}}} \le M|a-\grave{a}|,
$$
and thus
$$
\norm{x_{\ep}^a-x_{\ep}^{\grave{a}}} \le 2M|a-\grave{a}|. 
$$
This proves \eqref{fp periodic lipschitz}.
\end{proof}

\subsection{Solution of the fixed-point problem}
It is convenient to introduce some new notation.  Let
$$
\L_1^{\ep}(t) 
:= \begin{bmatrix*}
 \varpi^{\ep,\omega_{\ep}+t} &0 \\
0 & \ep^2\left(\xi^{\ep,t}\right)^{-1}\Pi_2
\end{bmatrix*},
$$
$$
\L_2^{\ep}(t) 
:= -\begin{bmatrix*}
1&0 \\
0 &\lambda_+^{\ep(\omega_{\ep}+t)}
\end{bmatrix*},
$$
$$
\L_3^{\ep}(t) :=
J^{\ep(\omega_{\ep}+t)},
$$
$$
\L_4^{\ep}(t) := J_1^{\ep(\omega_{\ep}+t)},
$$
$$
\G_{\ep}(\psib,t,a) :=\L_2^{\ep}(t)\left(\B^{\ep}(\psib,t)+\xtr^{\ep}(\psib,t,a)\right),
$$
and
$$
\F_{\ep}(\psib,t,a) 
:= a\L_1^{\ep}(t)\G_{\ep}(\psib,t,a)
= \begin{pmatrix*}
\Psi_1^{\ep}(\psib,t;a) \\
\Psi_2^{\ep}(\psib,t;a)
\end{pmatrix*}.
$$
We obtain the estimates \eqref{fp lemma: bound}, \eqref{fp lemma: lipschitz vector}, and \eqref{fp lemma: lipschitz a} for the function $\Psib^{\ep}$ by showing they hold for the function $\F_{\ep}$ just defined and then, separately, for $\Psi_3^{\ep}$, for which the map $\G_{\ep}$ will be useful. In our application of Lemma \ref{fp} we will take the Banach space to be $\X = \W^2$.  However, we will prove various estimates in the spaces $\W^r$ for the sake of the subsequent bootstrap arguments.  

\begin{remark}
In this appendix we will denote by $\b(\X,\Y)$ the space of bounded linear operators between the Banach spaces $\X$ and $\Y$.
\end{remark}

\subsection{Estimates for $\G_{\ep}$} 

\begin{proposition}\label{periodic workhorse for G}
Let $\ep_{\per} >0$ be as in Proposition \ref{master periodic lemma}.  We have the following estimates for all $r \ge 2$.

\begin{enumerate}[label={\bf(\roman*)}]
\item
There exists an increasing function $\M_{\map,r}^{\G} \colon \R_+ \to \R_+$ such that 
\begin{equation}\label{g bound}
\sup_{\substack{0 < \ep < \ep_{\per} \\ |t| \le 1}} \norm{\G_{\ep}(\psib,t,a)}_r 
\le \M_{\map,r}^{\G}[\norm{\psib}_r]
\end{equation}
for any $\psib \in \W^r$ and $|a| \le 1$.

\item 
There exists a radially increasing function $\M_{\lip,r}^{\G} \colon \R_+^2 \to \R_+$ such that 
\begin{equation}\label{g lip}
\sup_{0 < \ep < \ep_{\per}} \norm{\G_{\ep}(\psib,t,a)-\G_{\ep}(\grave{\psib},\grave{t},a)}_{r-1} 
\le \M_{\lip,r}^{\G}[\norm{\psib}_r,\bignorm{\grave{\psib}}_r]\left(\bignorm{\psib-\grave{\psib}}_{r} + |t-\grave{t}|\right)
\end{equation}
for any $\psib$, $\grave{\psib} \in \W^r$, $|t|$, $|\grave{t}| \le 1$, and $|a| \le 1$.  
\end{enumerate}
\end{proposition}

\begin{proof}
Throughout this proof, we assume $\ep$, $|t|$, $|\grave{t}|$, $|a| \le 1$ and $\psib \in \W^r$.
\begin{enumerate}[label={\bf(\roman*)}]
\item A first pass using \eqref{L2 L3 L4 mapping estimate} shows
\begin{align*}
\norm{\G_{\ep}(\psib,t,a)}_r
&\le |a|\norm{\L_2^{\ep}(t)}_{\b(\W^{r},\W^{r})}\left(\norm{\B^{\ep}(\psib,t)}_{r}+\norm{\xtr^{\ep}(\psib,t,a)}_{r}\right) \\
\\
&\le |a|C_{\map}^2\left(\norm{\B^{\ep}(\psib,t)}_{r}+\norm{\xtr^{\ep}(\psib,t,a)}_{r}\right).
\end{align*}
Since $r \ge 2$, we can use the Sobolev embedding to estimate the products in $\B^{\ep}$ and $\xtr^{\ep}$:
\begin{align*}
\norm{\B^{\ep}(\psib,t)}_{r}
&= \norm{B^{\ep(\omega_{\ep}+t)}(\nub+\psib,\nub+\psib)}_{r} \\
\\
&=\norm{\L_4^{\ep(\omega_{\ep}+t)}M_{\beta/\kappa}[\L_3^{\ep}(t)(\nub+\psib)]^{.2}}_{r} \\
\\
&\le C_{\map}\frac{\beta}{\kappa}\norm{\L_3^{\ep}(t)(\nub+\psib)}_{r}^2 \\
\\
&\le C_{\map}^3\frac{\beta}{\kappa}\left(\norm{\nub}_{r}^2 + 2\norm{\nub}_{r}\norm{\psib}_{r}+\norm{\psib}_{r}^2\right) \\
\\
&\le C_{\map}^3C_r\frac{\beta}{\kappa}\left(1 + \norm{\psib}_{r}+\norm{\psib}_{r}^2\right)
\end{align*}
and, similarly,
\begin{align*}
\norm{\xtr^{\ep}(\psib,t,a)}_{r} 
&= \norm{\L_4^{\ep}(t)M_{1/\kappa}[(\L_3^{\ep}(t)(\nub+\psib))^{.3}.N(a\ep^2\L_3^{\ep}(t)(\nub+\psib))]}_r \\
\\
&\le C_{\map}^2\frac{1}{\kappa}\norm{\nub+\psib}_{r}^{3}\norm{N(a\ep^2\L_3^{\ep}(t)(\nub+\psib))}_{r}.
\end{align*}
We apply Proposition \ref{Hrper mapping} to estimate
$$
\norm{N(a\ep^2\L_3^{\ep}(t)(\nub+\psib))}_{r}
\le \M_{r}[|a|\ep^2\norm{\L_3^{\ep}(t)(\nub+\psib)}_{r}]
$$
for some increasing function $\M_{r} \colon \R_+ \to \R_+$.  Since 
$$
\sup_{\substack{|a|, |t| \le 1 \\ 0 < \ep < \ep_{\per}}} |a|\ep^2\norm{\L_3^{\ep}(t)(\nub+\psib)}_{r-1} \le C_r\left(1+\norm{\psib}_{r}\right)
$$
and $\M_{r}$ is increasing, we have
$$
\norm{N(a\ep^2\L_3^{\ep}(t)(\nub+\psib))}_{r} \le \M_{r}[\norm{\psib}_{r}]. 
$$
This bound, together with the estimate on $\norm{\B^{\ep}(\psib,t)}_{r}$ above, produces \eqref{g bound}.

\item We have
\begin{align*}
\bignorm{\G_{\ep}(\psib,t,a)-\G_{\ep}(\grave{\psib},\grave{t},a)}_{r-1}
&= \bunderbrace{\norm{(\L_2^{\ep}(t)-\L_2^{\ep}(\grave{t}))(\B^{\ep}(\psib,t)+\xtr^{\ep}(\psib,t,a))}_{r-1}}{\Delta_1} \\
\\
&+ \bunderbrace{\bignorm{\L_2^{\ep}(\grave{t})(\B^{\ep}(\psib,t)-\B^{\ep}(\grave{\psib},\grave{t}))}_{r-1}}{\Delta_2} \\
\\
&+ \bunderbrace{\bignorm{\L_2^{\ep}(\grave{t})(\xtr^{\ep}(\psib,t,a)-\xtr^{\ep}(\grave{\psib},\grave{t},a))}_{r-1}}{\Delta_3}.
\end{align*}
It is straightforward to estimate $\Delta_1$ using \eqref{L2 L3 L4 Lipschitz estimate} and Proposition \ref{Hrper mapping}:
\begin{align*}
\Delta_1 &\le \norm{\L_2^{\ep}(t)-\L_2^{\ep}(\grave{t})}_{\b(\W^{r},\W^{r-1})}\norm{\B^{\ep}(\psib,t)+\xtr^{\ep}(\psib,t,a)}_{r} \\
\\
&\le C_{\lip}|t-\grave{t}|\norm{\B^{\ep}(\psib,t)+\xtr^{\ep}(\psib,t,a)}_r \\
\\
&\le C_{\lip}|t-\grave{t}|\M_r[\norm{\psib}_r].
\end{align*}

We handle $\Delta_2$ and $\Delta_3$ in essentially the same way; we treated $\Delta_2$ in \cite{faver-wright}, so we provide some more detail for $\Delta_3$ here:
\begin{align*}
\Delta_3
&\le \norm{\L_2^{\ep}(\grave{t})}_{\b(\W^{r-1},\W^{r-1})}\bignorm{\xtr^{\ep}(\psib,t,a)-\xtr^{\ep}(\grave{\psib},\grave{t},a)}_{r-1} \\
\\
&\le C_{\map} \bunderbrace{\bignorm{\xtr^{\ep}(\psib,t,a)-\xtr^{\ep}(\grave{\psib},\grave{t},a)}_{r-1}}{\Delta_3},
\end{align*}
where
\begin{multline*}
\xtr^{\ep}(\psib,t,a)-\xtr^{\ep}(\grave{\psib},\grave{t},a)
=\L_4^{\ep}(t)M_{1/\kappa}[(\L_3^{\ep}(t)(\nub+\psib))^{.3}.N(a\ep^2\L_3^{\ep}(t)(\nub+\psib))] \\
\\
- \L_4^{\ep}(\grave{t})M_{1/\kappa}[(\L_3^{\ep}(\grave{t})(\nub+\grave{\psib}))^{.3}.N(a\ep^2\L_3^{\ep}(\grave{t})(\nub+\grave{\psib}))].
\end{multline*}
Then, adding a number of zeroes, we can bound $\Delta_4$ in the natural way, and the only term the likes of which we have not seen before will be
$$
\bignorm{\L_4^{\ep}(\grave{t})M_{1/\kappa}[(\L_3^{\ep}(\grave{t})(\nub+\grave{\psib}))^{.3}.(N(a\ep^2\L_3^{\ep}(t)(\nub+\psib))-N(a\ep^2\L_3^{\ep}(\grave{t})(\nub+\grave{\psib})))]}_{r-1}.
$$
After factoring out the operators and using the Sobolev inequality for products, we invoke Proposition \ref{Hrper lipschitz} to bound
\begin{multline*}
\bignorm{N(a\ep^2\L_3(t)(\nub+\psib))-N(a\ep^2\L_3^{\ep}(\grave{t})(\nub+\grave{\psib}))}_{r-1} \\
\\
\le \M_{r-1}[\norm{a\ep^2\L_3^{\ep}(t)(\nub+\psib)}_{r-1},\bignorm{a\ep^2\L_3^{\ep}(\grave{t})(\nub+\grave{\psib})}_{r-1}]\bignorm{\L_3(t)(\nub+\psib)-\L_3^{\ep}(\grave{t})(\nub+\grave{\psib})}_{r-1}.
\end{multline*}
Since $\M_{r-1}$ is radially increasing, the uniform bound on $\L_3^{\ep}$ and the triangle inequality allow us to bound this by 
\begin{multline*}
\M_{r-1}[\norm{\psib}_r,\bignorm{\grave{\psib}}_r]\left(\norm{(\L_3^{\ep}(t)-\L_3^{\ep}(\grave{t}))\nub}_{r-1}+\norm{(\L_3^{\ep}(t)-\L_3^{\ep}(\grave{t}))\psib}_{r-1}+\bignorm{\L_3^{\ep}(\grave{t})(\psib-\grave{\psib})}_{r-1}\right),
\end{multline*}
and we estimate these terms easily enough using \eqref{L2 L3 L4 Lipschitz estimate} to achieve \eqref{g lip}. \qedhere
\end{enumerate}
\end{proof}

\subsubsection{Estimates for $\F_{\ep}$} We will obtain the estimates \eqref{fp lemma: bound}, \eqref{fp lemma: lipschitz vector}, and \eqref{fp lemma: lipschitz a} for $\F_{\ep}$ directly from the following proposition.

\begin{proposition}\label{periodic workhorse for F}
Let $\ep_{\per} >0$ be as in Proposition \ref{master periodic lemma}. The following estimates hold for all $r \ge 2$.

\begin{enumerate}[label={\bf(\roman*)}]
\item
There exists an increasing function $\M_{\map,r}^{\F} \colon \R_+ \to \R_+$ such that 
\begin{equation}\label{f bound}
\sup_{\substack{0 < \ep < \ep_{\per} \\ |t| \le 1}} \norm{\F_{\ep}(\psib,t,a)}_r 
\le |a|\M_{\map,r}^{\F}[\norm{\psib}_{r-1}]
\end{equation}
for any $\psib \in \W^r$ and $|a| \le 1$.

\item 
There exists a radially increasing function $\M_{\lip,r}^{\F} \colon \R_+^2 \to \R_+$ such that 
\begin{equation}\label{f lip}
\sup_{0 < \ep < \ep_{\per}} \bignorm{\F_{\ep}(\psib,t,a)-\F_{\ep}(\grave{\psib},\grave{t},a)}_r 
\le |a|\M_{\lip,r}^{\F}[\norm{\psib}_r,\bignorm{\grave{\psib}}_r]\left(\bignorm{\psib-\grave{\psib}}_{r-1} + |t-\grave{t}|\right)
\end{equation}
for any $\psib$, $\grave{\psib} \in \W^r$, $|t|$, $|\grave{t}| \le 1$, and $|a| \le 1$.  

\item 
There exists a continuous function $\M_{\max,r}^{\F} \colon \R_+ \to \R_+$ such that 
\begin{equation}\label{f a lip}
\sup_{0 < \ep < \ep_{\per}} \norm{\F_{\ep}(\psib,t,a)-\F_{\ep}(\psib,t,\grave{a})}_r
\le \M_{\max,r}^{\G}[\norm{\psib}_r]|a-\grave{a}|.
\end{equation}
\end{enumerate}
\end{proposition}

\begin{proof}
\begin{enumerate}[label={\bf(\roman*)}]
\item 
We use the smoothing property of $\L_1^{\ep}(t)$, \eqref{L1 mapping estimate}, and \eqref{g bound} to find
\begin{align*}
\norm{\F_{\ep}(\psib,t,a)}_r
&= |a|\norm{\L_1^{\ep}(t)\G_{\ep}(\psib,t,a)}_r \\
\\
&\le |a|\norm{\L_1^{\ep}(t)}_{\b(\W^{r-1},\W^r)}\norm{\G_{\ep}(\psib,t,a)}_{r-1} \\
\\
&\le C_{\map}|a|\M_{\map,r-1}^{\G}[\norm{\psib}_{r-1}].
\end{align*}

\item 
First, we have
\begin{align*}
\bignorm{\F_{\ep}(\psib,t,a)-\F_{\ep}(\grave{\psib},\grave{t},a)}_r
&= |a|\norm{\L_1^{\ep}(t)\G_{\ep}(\psib,t,a)-\L_1^{\ep}(\grave{t})\G_{\ep}(\grave{\psib},\grave{t},a)}_r \\
\\
&\le |a|\bunderbrace{\bignorm{(\L_1^{\ep}(t)-\L_1^{\ep}(\grave{t}))\G_{\ep}(\psib,t,a)}_r}{\Delta_1}
+ |a|\bunderbrace{\bignorm{\L_1^{\ep}(\grave{t})(\G_{\ep}(\psib,t,a)-\G_{\ep}(\grave{\psib},\grave{t},a))}_r}{\Delta_2}.
\end{align*}
We apply \eqref{L1 lipschitz estimate} to $\Delta_1$ and then \eqref{g bound} to find
\begin{align*}
\Delta_1 
&\le |a|\norm{\L_1^{\ep}(t)-\L_1^{\ep}(\grave{t})}_{\b(\W^r,\W^r)}\norm{\G_{\ep}(\psib,t,a)}_r \\
\\
&\le |a|C_{\lip}|t-\grave{t}|\M_{\map,r}^{\G}[\norm{\psib}_r].
\end{align*}
For $\Delta_2$, we need the smoothing property of $\L_1^{\ep}(\grave{t})$:
\begin{align*}
\Delta_2 &\le 
|a|\norm{\L_1^{\ep}(\grave{t})}_{\b(\W^{r-2},\W^r)}\bignorm{\G_{\ep}(\psib,t,a)-\G_{\ep}(\grave{\psib},\grave{t},a)}_{r-2} \\
\\
&\le C_{\map}|a|\M_{\lip,r}^{\G}[\norm{\psib}_r,\bignorm{\grave{\psib}}_r]\left(\bignorm{\psib-\grave{\psib}}_{r-1} + |t-\grave{t}|\right). 
\end{align*}

\item The necessary estimates are straightfoward and similar enough to the proof of \eqref{f lip} that we omit them; they rely fundamentally on the mapping and Lipschitz estimates \eqref{g bound} and \eqref{g lip} for $\G_{\ep}$ and on the Lipschitz composition estimate in Proposition \ref{Hrper lipschitz}.  
\qedhere
\end{enumerate}
\end{proof}

\subsubsection{Estimates for $\Psi_3^{\ep}$}  \label{estimates for Psi3ep appendix}
In this section, we only need estimates for $r=2$.  Fix $\ep \in (0,\ep_{\per})$.  

\begin{itemize}[leftmargin=*]
\item {\bf{Proof of \eqref{fp lemma: bound} for $\Psi_3^{\ep}$.}}
The triangle inequality gives
$$
|\Psi_3^{\ep}(\psib,t,a)|
\le \frac{1}{\Upsilon_{\ep}}|\rhs_{\ep}(\ep{t})|t^2 + \frac{2|a|}{\Upsilon_{\ep}}\left|\ft\left[\lambda_+^{\ep(\omega_{\ep}+t)}(\B_2^{\ep}(\psib,t)+\xtr_2^{\ep}(\psib,t,a))\right](1)\right|
$$
The bound \eqref{upsilon bound} on $\Upsilon_{\ep}$ and \eqref{txi mapping estimate} bound the first term above by $t^2/b_0$, and elementary properties of the Fourier transform, the Sobolev embedding, and \eqref{g bound} imply
\begin{align*}
\left|\ft\left[\lambda_+^{\ep(\omega_{\ep}+t)}(\B_2^{\ep}(\psib,t)+\xtr_2^{\ep}(\psib,t,a))\right]\!(1)\right|
&\le \norm{\lambda_+^{\ep(\omega_{\ep}+t)}(\B_2^{\ep}(\psib,t)+\xtr_2^{\ep}(\psib,t,a))}_{L^{\infty}} \\
\\
&\le C\norm{\lambda_+^{\ep(\omega_{\ep}+t)}(\B_2^{\ep}(\psib,t)+\xtr_2^{\ep}(\psib,t,a))}_2 \\
\\
&\le C\norm{\L_2^{\ep}(t)(\B^{\ep}(\psib,t)+\xtr^{\ep}(\psib,t,a))}_2 \\
\\
&= C\norm{\G_{\ep}(\psib,t,a)}_2 \\
\\
&\le C\M_{\map,2}^{\G}[\norm{\psib}_2].
\end{align*}
All together, we have 
$$
|\Psi_3^{\ep}(\psib,t,a)|
\le C_{\map}t^2+2C|a|\M_{\map,2}^{\G}[\norm{\psib}_2]. 
$$

\item {\bf{Proof of \eqref{fp lemma: lipschitz vector} for $\Psi_3^{\ep}$.}} We use the triangle inequality to estimate
\begin{align*}
|\Psi_3^{\ep}(\psib,t,a)-\Psi_3^{\ep}(\grave{\psib},\grave{t},a)|
&\le Ct^2|\rhs_{\ep}(\ep{t})-\rhs_{\ep}(\ep\grave{t})| + |\rhs_{\ep}(\ep\grave{t})||t^2-\grave{t}^2| \\
\\
&+ C|a|\norm{\G_{\ep}(\psib,t,a)-\G_{\ep}(\grave{\psib},\grave{t},a)}_1.
\end{align*}
Here we have used the Sobolev embedding inequality to bound the Fourier transform terms by their $\W^1$-norm instead of the $\W^2$-norm as we did above with the mapping estimate.  This allows us to use \eqref{g lip} to bound
$$
\norm{\G_{\ep}(\psib,t,a)-\G_{\ep}(\grave{\psib},\grave{t},a)}_1
\le \M_{\lip,2}^{\G}[\norm{\psib}_2,\bignorm{\grave{\psib}}_2]\left(\bignorm{\psib-\grave{\psib}}_2+|t-\grave{t}|\right).
$$
We estimate the first term with \eqref{txi lipschitz estimate}:
$$
|\rhs_{\ep}(\ep{t})-\rhs_{\ep}(\ep\grave{t})|
\le C_{\lip}|t-\grave{t}|
$$
and the second by the difference of squares and \eqref{txi mapping estimate}.

\item {\bf{Proof of \eqref{fp lemma: lipschitz a} for $\Psi_3^{\ep}$.}} The proof is the same as that of \eqref{f a lip} above.
\end{itemize}

\subsubsection{Proof of the remainder of Theorem \ref{periodic main theorem}} Taking $r=2$ in Proposition \ref{periodic workhorse for F} and using the estimates for $\Psi_3^{\ep}$ in Appendix  \ref{estimates for Psi3ep appendix}, we find that the map $\Psib_{\ep}$ satisfies the estimates of Lemma \ref{fp} on the space $\W^2 = \left(H_{\per}^2 \times H_{\per}^2\right) \cap \mathcal{Z}$, there exist $r_0$, $\ep_{\per}$, $a_{\per} \in (0,1)$ such that for all $0 < \ep < \ep_{\per}$ and $|a| \le a_{\per}$, there is a unique $(\psib_{\ep}^a,t_{\ep}^a) = (\psi_{1,\ep}^a,\psi_{2,\ep}^a,t_{\ep}^a) \in \W^2 \times \R$ satisfying
$$
\Psib_{\ep}(\psib_{\ep}^a,t_{\ep}^a,a) = (\psib_{\ep}^a,t_{\ep}^a)
$$
and
$$
\norm{(\psib_{\ep}^a,t_{\ep}^a)}_{\W^2 \times \R} = \norm{\psi_{\ep,1}^a}_{H_{\per}^2} + \norm{\psi_{\ep,2}^a}_{H_{\per}^2} + |t_{\ep}^a| \le r_0.
$$
We also have a constant $C_2 > 0$ such that for $|a|$, $|\grave{a}| \le a_{\per}$ and $0 < \ep < \ep_{\per}$, the Lipschitz estimate
\begin{equation}\label{periodic lipschitz conclusion}
\norm{(\psib_{\ep}^a,t_{\ep}^a)-(\psib_{\ep}^{\grave{a}},t_{\ep}^{\grave{a}})}_{\W^r \times \R} 
\le C_2|a-\grave{a}|
\end{equation}
holds.  Now set
$$
\omega_{\ep}^a := \omega_{\ep} + t_{\ep}^a.  
$$
We are ready to prove the rest of Theorem \ref{periodic main theorem}.

\begin{itemize}[leftmargin=*]
\item 
{\bf{Proof of (i).}}
Undoing the fixed point set-up of Appendix \ref{fixed point conversion appendix}, we find that
$$
\thetab(X) 
:=a\varphib_{\ep}^a(X)
:=a\nub(\omega_{\ep}^aX)+a\psib_{\ep}^a(\omega_{\ep}^aX)
$$
solves \eqref{lw eqns}.

\item
{\bf{Proof of (ii) and (iii).}}
Recalling the definitions of the components of $\Psib_{\ep}$ in \eqref{fp1}, \eqref{fp2}, and \eqref{fp3}, we compute
$$
(\psi_{1,\ep}^0, \psi_{2,\ep}^0,t_{\ep}^0) = \Psib_{\ep}(\psib_{\ep}^0,t_{\ep}^0,0) = (0,0,0),
$$
and so $\psib_{\ep}^0 = 0$ and $t_{\ep}^0 = 0$, hence $\omega_{\ep}^0 = \omega_{\ep} + t_{\ep}^0 = \omega_{\ep}$.  

For the $\O(\ep)$ bound on $\omega_{\ep}$, we know from \eqref{omegac intermediate bound} that with $\cep = \sqrt{c_{\kappa}^2+\ep^2}$ and $\omega_{\ep}  = \Omega_{\cep}/\ep$, we have
$$
\frac{1}{\ep}\left(\frac{\sqrt{2\kappa}}{\sqrt{c_{\kappa}^2+\ep_{\per}^2}}\right)
\le \frac{1}{\ep}\frac{\sqrt{2\kappa}}{\sqrt{c_{\kappa}^2+\ep^2}}
\le \omega_{\ep} 
\le \frac{1}{\ep}\frac{\sqrt{2 + 2\kappa}}{\sqrt{c_{\kappa}^2+\ep^2}}
\le \frac{1}{\ep}\left(\frac{\sqrt{2+2\kappa}}{c_{\kappa}}\right).
$$

\item 
{\bf{Proof of (iv).}}
Since $\psib_{\ep}^a = \F_{\ep}(\psib_{\ep}^a,t_{\ep}^a,a)$, we use \eqref{f bound} to estimate
$$
\norm{\psib_{\ep}^a}_r \le \M_{\map,r}^{\F}[\norm{\psib_{\ep}^a}_{r-1}].
$$
When $r = 3$, we know that $\norm{\psib_{\ep}^a}_2 \le r_0$, so the continuity of $\M_{\map,2}^{\F}$ implies 
$$
\norm{\psib_{\ep}^a}_3 \le a_{\per}\max_{0 \le t \le r_0} \M_{\map,2}^{\F}[t].
$$
Induction on $r$ then furnishes a constant $C_r > 0$ such that 
$$
\norm{\psib_{\ep}^a}_r \le C_r, \ r \ge 3.
$$
Next,
\begin{align*}
\norm{\psib_{\ep}^a-\psib_{\ep}^{\grave{a}}}_r
&= \norm{\F_{\ep}(\psib_{\ep}^a,t_{\ep}^a,a)-\F_{\ep}(\psib_{\ep}^{\grave{a}},t_{\ep}^{\grave{a}},\grave{a})}_r \\
\\
&\le \bunderbrace{\norm{\F_{\ep}(\psib_{\ep}^a,t_{\ep}^a,a)-\F_{\ep}(\psib_{\ep}^{\grave{a}},t_{\ep}^{\grave{a}},a)}_r}{\Delta_{1,r}}
+ \bunderbrace{\norm{\F_{\ep}(\psib_{\ep}^{\grave{a}},t_{\ep}^{\grave{a}},a)-\F_{\ep}(\psib_{\ep}^{\grave{a}},t_{\ep}^{\grave{a}},\grave{a})}_r}{\Delta_{2,r}}
\end{align*}
We bound $\Delta_{1,r}$ using \eqref{f lip}:
$$
\Delta_{1,r} \le \M_{\lip,r}^{\F}[\norm{\psib_{\ep}^a}_r,\norm{\psib_{\ep}^{\grave{a}}}_r]\left(\norm{\psib_{\ep}^a-\psib_{\ep}^{\grave{a}}}_{r-1}+|t_{\ep}^a-t_{\ep}^{\grave{a}}|\right)
$$
and $\Delta_{2,r}$ using \eqref{f a lip}:
$$
\Delta_{2,r} \le \M_{\max,r}^{\F}[\norm{\psib_{\ep}^{\grave{a}}}_r]|a-\grave{a}|.
$$
The uniform bounds on $\norm{\psib_{\ep}^a}_r$ and $\norm{\psib_{\ep}^a}_r$ and the continuity of $\M_{\lip,r}^{\F}$ and $\M_{\max,r}^{\F}$ then allow us to bound
$$
\Delta_{1,r} + \Delta_{2,r} \le C_r\left(\norm{\psib_{\ep}^a - \psib_{\ep}^{\grave{a}}}_{r-1} + |t_{\ep}^a-t_{\ep}^{\grave{a}}| + |a-\grave{a}|\right).
$$
All that remains is to induct on $r$ using the base case \eqref{periodic lipschitz conclusion}.
\end{itemize}

%% file: appendix_nanopteron_estimates_master.tex
\section{The Nanopteron Estimates} 
\label{all the pretty estimates appendix}

In this appendix we prove a long series of estimates that undergird the proof of Proposition \ref{main workhorse estimates}.  

\input{appendix_proof_of_main_workhorse_estimates}
\input{appendix_mapping_estimates}
\input{appendix_lipschitz_estimates}

\input{appendix_bootstrap_estimates}
\input{appendix_operator_estimates}

%% file: appendix_proof_of_main_workhorse_estimates.tex
\subsection{Proof of Proposition \ref{main workhorse estimates}} The following proposition contains the estimates needed to prove Proposition \ref{main workhorse estimates}.

\begin{proposition}
Let $\epbar \in (0,\ep_{\per})$ be as in Proposition \ref{faver wright nanopteron}.
\begin{enumerate}[label={\bf(\roman*)}]
\item There exists an increasing function $\M_{\map} \colon \R_+ \to \R_+$ such that 
\begin{equation}\label{M mapping est}
\norm{\nsystb^{\ep}(\etab,a)}_{\X^1} \le \M_{\map}[\norm{\etab}_{1,\qstar}]\left(\ep + \ep\norm{(\etab,a)}_{\X^1} + \norm{(\etab,a)}_{\X^1}^2\right)
\end{equation}
for all $\etab \in E_{\qstar}^1 \times E_{\qstar}^1$, $|a| \le a_{\per}$, and $0 < \ep < \epbar$.
\item There exists a radially increasing function $\M_{\lip} \colon \R_+^2 \to \R_+$ such that  
\begin{equation}\label{M lip est}
\norm{\nsystb^{\ep}(\etab,a)-\nsystb^{\ep}(\grave{\etab},\grave{a})}_{\X^0} \le \M_{\lip}[\norm{\etab}_{1,\qstar},\norm{\grave{\etab}}_{1,\qstar}]\left(\ep + \norm{(\etab,a)}_{\X^1} + \norm{(\grave{\etab},\grave{a})}_{\X^1}\right)\norm{(\etab,a)-(\grave{\etab},\grave{a})}_{\X^0}
\end{equation}
for all $\etab,\grave{\etab} \in E_{\qstar}^1 \times E_{\qstar}^1$, $|a| \le a_{\per}$, and $0 < \ep < \epbar$.
\item For all integers $r \ge 1$ there exists an increasing function $\M_{\boot,r} \colon \R_+ \to \R_+$ such that
\begin{equation}\label{M boot est 1}
\norm{\nsystb^{\ep}(\etab,a)}_{r+1,\qstar}
\le \M_{\boot,r}[\norm{\etab}_{r,\qstar}]
\left(\ep + \norm{\etab}_{r,\qstar}+|a|\ep^{-r}\norm{\etab}_{r,\qstar}+|a|\ep^{1-r}+a^2\ep^{1-2r}+|a|^3\ep^{1-3r}\right)
\end{equation}
and
\begin{equation}\label{M boot est 2}
|\nsyst_3^{\ep}(\etab,a)| 
\le \M_{\boot,r}[\norm{\etab}_{r,\qstar}]
\left(\ep^{r+1}+\ep^r\norm{\etab}_{r,\qstar}+|a|\norm{\etab}_{r,\qstar}+\ep|a|+a^2\ep^{1-r}+|a|^3\ep^{1-2r}\right)
\end{equation}
for all $\etab \in E_{\qstar}^r \times E_{\qstar}^r$, $|a| \le a_{\per}$, and $0 < \ep < \epbar$.
\end{enumerate}
\end{proposition}

The proof of \eqref{M mapping est} is developed in Appendix \ref{mapping estimates appendix}, of \eqref{M lip est} in Appendix \ref{lipschitz estimates appendix}, and of \eqref{M boot est 1} and \eqref{M boot est 2} in Appendix \ref{bootstrap estimates appendix}. Now we are ready to to prove Proposition \ref{main workhorse estimates}.

\begin{enumerate}[label={\bf{Proof of (\roman*).}}]
\item Let $M_{\star} > 0$ be such that 
$$
|\M_{\map}[\norm{\etab}_{1,\qstar}]| \le M_{\star} 
\quadword{and}
|\M_{\lip}[\norm{\etab}_{1,\qstar},\norm{\grave{\etab}}_{1,\qstar}]| \le M_{\star}
$$
whenever $\norm{\etab}_{1,\qstar} \le 1$ and $\norm{\grave{\etab}}_{1,\qstar} \le 1$.  Set $\tau = M_{\star} + 2$ and 
\begin{equation}\label{defn of ep-star}
\ep_{\star} = \min\left\{\frac{1}{M_{\star}(M_{\star}+2)}, \frac{1}{M_{\star}(M_{\star}+2)^2}, \frac{1}{2M_{\star}(M_{\star}+5)},\frac{1}{M_{\star}+2},1,\epbar\right\}
\end{equation}
Observe that if $(\etab,a) \in \U_{\ep,M_{\star}+2}^1$, then
\begin{equation}\label{<=1 obs}
\norm{\etab}_{1,\qstar}
\le
\norm{(\etab,a)}_{\X^1} 
\le \ep(M_{\star}+2)
\le \ep_{\star}(M_{\star}+2)
\le 1.
\end{equation}
\begin{itemize}[leftmargin=*]
\item For the mapping estimate, let $0 < \ep < \ep_{\star}$ and $(\etab,a) \in \U_{\ep,M_{\star}+2}^1$.  Then \eqref{<=1 obs} and \eqref{M mapping est} allow us to estimate
\begin{align*}
\norm{\nsystb^{\ep}(\etab,a)}_{\X^1} 
&\le M_{\star}\left(\ep + \ep^2(M_{\star}+2) + \ep^2(M_{\star}+2)^2\right) \\
\\
&= \ep\left(M_{\star} + M_{\star}(M_{\star}+2)\ep + M_{\star}(M_{\star}+2)^2\ep\right) \\
\\
&\le \ep\left(M_{\star} + M_{\star}(M_{\star}+2)\ep_{\star} + M_{\star}(M_{\star}+2)^2\ep_{\star}\right) \\
\\
&\le \ep\left(M_{\star} + 1 + 1\right) \\
\\
&= (M_{\star}+2)\ep.
\end{align*}
Hence $\nsystb^{\ep}(\etab,a) \in \U_{\ep,M_{\star}+2}^1$.

\item For the Lipschitz estimate, take $0 < \ep < \ep_{\star}$ and $(\etab,a)$, $(\grave{\etab},\grave{a}) \in \U_{\ep,M_{\star}+2}^1$.  Then \eqref{<=1 obs} and \eqref{M lip est} imply
\begin{align*}
\norm{\nsystb^{\ep}(\etab,a)-\nsystb^{\ep}(\grave{\etab},\grave{a})}_{\X^0}
&\le M_{\star}\left(\ep + \norm{(\etab,a)}_{\X^1} + \norm{(\grave{\etab},\grave{a})}_{\X^1}\right)\norm{(\etab,a)-(\grave{\etab},\grave{a})}_{\X^0} \\
\\
&\le M_{\star}\left(\ep + 2\ep(M_{\star}+2)\right)\norm{(\etab,a)-(\grave{\etab},\grave{a})}_{\X^0} \\
\\
&= \ep{M}_{\star}(5+M_{\star})\norm{(\etab,a)-(\grave{\etab},\grave{a})}_{\X^0} \\
\\
&\le \ep_{\star}{M}_{\star}(5+M_{\star})\norm{(\etab,a)-(\grave{\etab},\grave{a})}_{\X^0} \\
\\
&\le \frac{1}{2}\norm{(\etab,a)-(\grave{\etab},\grave{a})}_{\X^0}.
\end{align*}
\end{itemize}

\item Now we prove the bootstrap estimates. That $\nsystb^{\ep}(\etab,a) \in \X^{r+1}$ if $(\etab,a) \in \X^r$ follows from the smoothing properties of $\P_{\ep}$ and $\varpi^{\ep}$.  Given $\tau > 0$, let
$$
M_{r,\tau} = \max_{\norm{\etab}_{r,\qstar} \le \tau\ep_{\star}} \M_{\boot,r}[\etab].
$$
Let $0 < \ep < \ep_{\star}$ and take $(\etab,a) \in \U_{\ep,\tau}^r$, so that $\norm{\etab}_{r,\qstar} \le \tau\ep \le \tau\ep_{\star}$ and $|a| \le \tau\ep^r$.  Then \eqref{M boot est 1} implies
\begin{align*}
\norm{\nsystb^{\ep}(\etab,a)}_{\X^{r+1}}
&\le M_{r,\tau}\left(\ep + \tau\ep + \tau\ep^r(\ep^{1-r})+ (\tau\ep)^2\ep^{1-2r}+(\tau\ep^r)\ep^{-r}(\tau\ep)+(\tau\ep^r)^3\ep^{1-3r}\right) \\
\\
&= M_{r,\tau}\left(\ep+2\tau\ep+2\tau^2\ep+\tau^3\ep\right) \\
\\
&= M_{r,\tau}\left(1+2\tau+2\tau^3+\tau^3\right)\ep.
\end{align*}
In particular, we find
$$
\norm{(\nsyst_1^{\ep}(\etab,a),\nsyst_2^{\ep}(\etab,a))}_{r+1,\qstar} 
\le \norm{\nsystb^{\ep}(\etab,a)}_{\X^{r+1}}
\le M_{r,\tau}\left(1+2\tau+2\tau^3+\tau^3\right)\ep.
$$
We need to refine this estimate, however, for $\nsyst_3^{\ep}(\etab,a)$. Using \eqref{M boot est 2}, we find
\begin{align*}
|\nsyst_3^{\ep}(\etab,a)| 
&\le M_{\tau,r}\left(\ep^{r+1}+\ep^r(\tau\ep)+(\tau\ep^r)(\tau\ep)+\ep(\tau\ep^r)+(\tau\ep^r)^2\ep^{1-r}+(\tau\ep)^3\ep^{1-2r}\right) \\
\\
&= M_{\tau,r}\left(\ep^{r+1}+2\tau\ep^{r+1}+2\tau^2\ep^{r+1}+\tau^3\ep^{r+1}\right) \\
\\
&= M_{\tau,r}\left(1+2\tau+2\tau^2+\tau^3\right)\ep^{r+1}.
\end{align*}
Then with 
$$
\overline{\tau} := M_{\tau,r}\left(1+2\tau+2\tau^2+\tau^3\right),
$$
we find that if $(\etab,a) \in \U_{\ep,\tau}^r$, then $\nsystb^{\ep}(\etab,a) \in \U_{\ep,\overline{\tau}}^{r+1}.$
\end{enumerate}

%% file: appendix_mapping_estimates.tex
\subsection{The Mapping Estimates}\label{mapping estimates appendix}
Let
\begin{equation}\label{rhs map defn}
\rhs_{\map}^{\ep}(\etab,a) 
= \ep + \ep\norm{\etab}_{1,q_{\star}} + \ep|a| + \norm{\etab}_{1,q_{\star}}^2 + |a|^2.
\end{equation}
Observe that 
$$
\rhs_{\map}^{\ep}(\etab,a) \le \ep + \ep\norm{(\etab,a)}_{1,q_{\star}} + \norm{(\etab,a)}_{1,q_{\star}}^2,
$$
so to obtain \eqref{M mapping est} it suffices to prove the existence of an increasing function $\M_{\map} \colon \R_+ \to \R_+$ such that 
$$
\norm{\nsystb^{\ep}(\etab,a)}_{\X^1} \le \M_{\map}[\norm{\etab}_{1,q_{\star}}]\rhs_{\map}^{\ep}(\etab,a).
$$

\subsubsection{General strategy}
The bounds on $\A^{-1}$, $\K_2$, $\P_{\ep}$, $\iota_{\ep}$, and $\upsilon_{\ep}$ from Proposition \ref{faver wright nanopteron} let us estimate
\begin{align*}
\norm{\nanorhs_2^{\ep,\mod}(\etab,a)} 
&\le \sum_{k=1}^5 \bignorm{\l_{k1}^{\ep,\mod}(\etab,a)}_{1,q_{\star}}
+ \sum_{k=1}^5 \norm{\l_{k,2}^{\ep}(\etab,a)}_{1,q_{\star}} 
+ \norm{\l_6^{\ep}(\etab,a)}_{1,q_{\star}}, \\
\\
\norm{\nsyst_2^{\ep}(\etab,a)}_{1,q_{\star}}
&= \ep^2\bignorm{\P_{\ep}\nanorhs_2^{\ep,\mod}(\etab,a)}_{1,q_{\star}} \\
\\
&\le C\sum_{k=1}^5 \ep\bignorm{\l_{k1}^{\ep,\mod}(\etab,a)}_{1,q_{\star}}
+C\sum_{k=1}^5 \ep\norm{\l_{k,2}^{\ep}(\etab,a)}_{1,q_{\star}} 
+ C\ep\norm{\l_6^{\ep}(\etab,a)}_{1,q_{\star}} \\
\\
|\nsyst_3^{\ep}(\etab,a)|
&= \frac{1}{2\upsilon_{\ep}}\bignorm{\iota_{\ep}\big[\nanorhs_2^{\ep,\mod}(\etab,a)\big]}_{1,\qstar} \\
\\
&\le C\sum_{k=1}^5 \ep\norm{\l_{k1}^{\ep,\mod}(\etab,a)}_{1,q_{\star}}
+C\sum_{k=1}^5 \ep\bignorm{\l_{k,2}^{\ep}(\etab,a)}_{1,q_{\star}} 
+ C\ep\norm{\l_6^{\ep}(\etab,a)}_{1,q_{\star}}, 
\end{align*}
and
\begin{align*}
\norm{\nsyst_1^{\ep}(\etab,a)}_{1,q_{\star}}
&\le \bignorm{\A^{-1}\nanorhs_1^{\ep,\mod}(\etab,a)}_{1,q_{\star}} + \norm{\A^{-1}\K_2[\nsyst_2^{\ep}(\etab,a)}_{1,q_{\star}} \\
\\
&\le C\sum_{k=1}^5 \left(\bignorm{\j_{k1}^{\ep,\mod}(\etab,a)}_{1,\qstar} + \ep\norm{\l_{k1}^{\ep,\mod}(\etab,a)}_{1,q_{\star}}\right) \\
\\
&+ C\sum_{k=1}^5 \left(\bignorm{\j_{k2}^{\ep}(\etab,a)}_{1,\qstar} + \ep\norm{\l_{k,2}^{\ep}(\etab,a)}_{1,q_{\star}}\right) \\
\\
&+ C\left(\norm{\j_6^{\ep}(\etab,a)}_{1,\qstar} + \ep\norm{\l_6^{\ep}(\etab,a)}_{1,q_{\star}}\right).
\end{align*}
We saw in \cite{faver-wright} that the terms
$$
\norm{\j_{1k}^{\ep,\mod}(\etab,a)}_{1,q_{\star}}
\quadword{and}
\ep\norm{\l_{1k}^{\ep,\mod}(\etab,a)}_{1,q_{\star}}, \ k = 1,\ldots,5
$$
are bounded above, up to a constant, by $\rhs_{\map}^{\ep}(\etab,a)$.  (To be fair, of course the mass dimer versions of $\j_{1k}^{\ep,\mod}$ and $\l_{1k}^{\ep,\mod}$ were entirely different functions, but the structure of the necessary estimates is exactly the same.)  Now we show that the terms
$$
\norm{\j_{2k}^{\ep}(\etab,a)}_{1,q_{\star}}, 
\qquad
\ep\norm{\l_{2k}^{\ep}(\etab,a)}_{1,q_{\star}}, \ k = 1,\ldots,5,
\qquad
\norm{\j_6^{\ep}(\etab,a)}_{1,q_{\star}},
\quadword{and}
\ep\norm{\l_6^{\ep}(\etab,a)}_{1,q_{\star}}
$$
are bounded above by $\M_{\map}[\norm{\etab}_{1,q_{\star}}]\rhs_{\map}^{\ep}(\etab,a)$.

\begin{remark}
Here and elsewhere in these appendices, we will abuse notation and write
$$
\norm{\varphib}_{W^{r,\infty}}
:= \norm{\varphi_1}_{W^{r,\infty}} + \norm{\varphi_2}_{W^{r,\infty}}
$$
for a function $\varphib = (\varphi_1,\varphi_2) \in W^{r,\infty} \times W^{r,\infty}$.
\end{remark}

\subsubsection{Mapping estimates for $\j_{12}$ and $\l_{12}$} \label{J12 L12 mapping estimates}
We present this first series of estimates in detail to show the general techniques and reliance on Proposition \ref{mapping} that will permeate the subsequent mapping estimates.  We have
\begin{align*}
\norm{\j_{12}^{\ep}(\etab,a)}_{1,q_{\star}}
&= \norm{\varpi^{\ep}\cb_2^{\ep}(\sigmab,\sigmab,\an_{\ep}(\etab,a))}_{1,q_{\star}} \\
\\
&\le C\norm{\cb_2^{\ep}(\sigmab,\sigmab,\an_{\ep}(\etab,a))}_{1,q_{\star}} \\
\\
&\le C\norm{\cb^{\ep}(\sigmab,\sigmab,\an_{\ep}(\etab,a))}_{1,q_{\star}} \\
\\
&= C\norm{J_1^{\ep}M_{1/\kappa}\left((J^{\ep}\sigmab)^{.2}.\nl(\ep^2J^{\ep}\an_{\ep}(\etab,a))\right)}_{1,q_{\star}} \\
\\
&\le C\norm{(J^{\ep}\sigmab)^{.2}.\nl(\ep^2J^{\ep}\an_{\ep}(\etab,a))}_{1,q_{\star}}.
\end{align*}
We remark that 
$$
(J^{\ep}\sigmab)^{.2}.\nl(\ep^2J^{\ep}\an_{\ep}(\etab,a))
=
\bunderbrace{(J^{\ep}\sigmab)^{.2}}{H_{q_{\star}}^1}.\nl\bigg(\ep\big(\bunderbrace{\ep J^{\ep}(\sigmab+\etab)}{H_{q_{\star}}^1}+a\bunderbrace{(\ep{J}^{\ep}\varphib_{\ep}^a)}{W^{1,\infty}}\big)\bigg),
$$
and so we could control this quantity with Proposition \ref{mapping prop}.  However, that would be premature, as the resulting estimate would not have the conducive form of $\rhs_{\map}^{\ep}(\etab,a)$.  Instead, in order to factor out an all-important power of $\ep$ (a recurring theme in these estimates), we expand $\nl$ and find
\begin{align*}
C\norm{(J^{\ep}\sigmab)^{.2}.\nl(\ep^2J^{\ep}\an_{\ep}(\etab,a))}_{1,q_{\star}}
&= C\norm{(J^{\ep}\sigmab)^{.2}.(\ep^2J^{\ep}\an_{\ep}(\etab,a)).N(\ep^2J^{\ep}\an_{\ep}(\etab,a))}_{1,q_{\star}} \\
\\
&= C\ep\norm{(J^{\ep}\sigmab)^{.2}.(\ep{J}^{\ep}\an_{\ep}(\etab,a)).N(\ep^2J^{\ep}\an_{\ep}(\etab,a))}_{1,q_{\star}}.
\end{align*}
Next, we use the definition of $\an_{\ep}(\etab,a)$ and the triangle inequality to break this norm into two terms:
\begin{align*}
C\ep\norm{(J^{\ep}\sigmab)^{.2}.(\ep{J}^{\ep}\an_{\ep}(\etab,a)).N(\ep^2J^{\ep}\an_{\ep}(\etab,a))}_{1,q_{\star}} 
&\le \bunderbrace{C\ep\norm{(J^{\ep}\sigmab)^{.2}.(\ep{J}^{\ep}(\sigmab+\etab)).N(\ep^2J^{\ep}\an_{\ep}(\etab,a))}_{1,q_{\star}}}{\Pi_1} \\
\\
&+ \bunderbrace{C\ep\norm{(J^{\ep}\sigmab)^{.2}.(\ep{J}^{\ep}(a\varphib_{\ep}^a)).N(\ep^2J^{\ep}\an_{\ep}(\etab,a))}_{1,q_{\star}}}{\Pi_2}.
\end{align*}
Observe that the product in $\Pi_1$ really has the form of the factors in the estimate in Proposition \ref{mapping}:
$$
(J^{\ep}\sigmab)^{.2}.(\ep{J}^{\ep}(\sigmab+\etab)).N(\ep^2J^{\ep}\an_{\ep}(\etab,a)) 
= \bunderbrace{(J^{\ep}\sigmab)^{.2}.(\ep{J}^{\ep}(\sigmab+\etab))}{H_q^1}.N\bigg(\ep\big(\bunderbrace{J^{\ep}(\sigmab+\etab)}{H_q^1} + a(\bunderbrace{\ep{J}^{\ep}\varphib_{\ep}^a}{W^{1,\infty}})\big)\bigg)
$$
That proposition implies
$$
\Pi_1 
\le C\ep\M[\norm{J^{\ep}(\sigmab+\etab)}_{1,q_{\star}}]\left(1+|a|\ep^{1-1}\right)\norm{(J^{\ep}\sigmab)^{.2}.(\ep{J}^{\ep}(\sigmab+\etab))}_{1,q_{\star}},
$$
for some increasing function $\M \colon \R_+ \to \R_+$.  The triangle inequality and the uniform bound on $J^{\ep}$ from Proposition \ref{faver wright nanopteron} give constants $C_1$, $C_2 > 0$ such that 
$$
\M[\norm{J^{\ep}(\sigmab+\etab)}_{1,q_{\star}}]\norm{(J^{\ep}\sigmab)^{.2}.(\ep{J}^{\ep}(\sigmab+\etab))}_{1,q_{\star}},
\le \M[C_1\norm{\etab}_{1,q_{\star}}]C_2(1+\norm{\etab}_{1,q_{\star}})
=: \M_{\l_{12}}[\norm{\etab}_{1,q_{\star}}].
$$
That is, $\M_{\l_{12}}$ is increasing with  
$$
\Pi_1
\le \M_{\l_{12}}[\norm{\etab}_{1,q_{\star}}]\ep,
$$
and this bound has the desired form of $\rhs_{\map}^{\ep}(\etab,a)$ from \eqref{rhs map defn}.

The estimate on $\Pi_2$ proceeds just as the one for $\Pi_1$, except first we factor
$$
\Pi_2 
\le C\ep\norm{\ep{J}^{\ep}(a\varphib_{\ep}^a)}_{W^{1,\infty}}\norm{(J^{\ep}\sigmab)^{.2}.N(\ep^2J^{\ep}\an_{\ep}(\etab,a))}_{1,q_{\star}} 
\le C\ep\norm{(J^{\ep}\sigmab)^{.2}.N(\ep^2J^{\ep}\an_{\ep}(\etab,a))}_{1,q_{\star}}
$$
and then apply Proposition \ref{mapping prop} to the term on the right.  The whole estimate for $\l_{12}$ follows in an identical way, so we omit it.

\subsubsection{Mapping estimates for $\j_{22}$ and $\l_{22}$} These estimates are essentially the same as the ones for $\j_{12}$ and $\l_{12}$, except wherever we had a factor of $(J^{\ep}\sigmab)^{.2}$ in the previous section, now we have a factor of $(J^{\ep}\sigmab).(J^{\ep}\etab)$.  We omit the details.

\subsubsection{Mapping estimates for $\j_{32}$ and $\l_{32}$} To handle the new presence of $\varphib_{\ep}^a$, which costs a factor of $\ep$ each time we estimate its $W^{1,\infty}$-norm, we need to expose an additional factor of $\ep$ in the estimates.  We achieve this via the smoothing property of $\varpi^{\ep}$ for $\j_{32}$ and the extra $\ep$ that naturally comes along with $\l_{32}$.  We have
\begin{align*}
\norm{\j_{32}^{\ep}(\etab,a)}_{1,q_{\star}} 
&2|a|\norm{\varpi^{\ep}\cb_2^{\ep}(\sigmab,\varphib_{\ep}^a,\an_{\ep}(\etab,a))}_{1,q_{\star}} \\
\\
&\le C|a|\norm{\cb_2^{\ep}(\sigmab,\varphib_{\ep}^a,\an_{\ep}(\etab,a))}_{0,q} \\
\\
&\le C|a|\norm{\cb^{\ep}(\sigmab,\varphib_{\ep}^a,\an_{\ep}(\etab,a))}_{0,q} \\
\\
&\le C|a|\norm{(J^{\ep}\sigmab).(J^{\ep}\varphib_{\ep}^a).\nl(\ep^2J^{\ep}\an_{\ep}(\etab,a))}_{0,q} \\
\\
&\le C\ep|a|\norm{J^{\ep}\varphib_{\ep}^a}_{W^{0,\infty}}\norm{(J^{\ep}\sigmab).(\ep{J}^{\ep}\an_{\ep}(\etab,a)).N(\ep^2J^{\ep}\an_{\ep}(\etab,a))}_{0,q} \\
\\
&\le C\ep|a|\norm{(J^{\ep}\sigmab).(\ep{J}^{\ep}\an_{\ep}(\etab,a)).N(\ep^2J^{\ep}\an_{\ep}(\etab,a))}_{0,q}.
\end{align*}
We estimate the factor in the $H_q^0$-norm above with Proposition \ref{mapping} and bound that as in the case of $\j_{12}$ to find, ultimately,
$$
\norm{\j_{32}^{\ep}(\etab,a)}_{1,q_{\star}} \le \M_{\j_{32}}[\norm{\etab}_{1,q_{\star}}]\ep|a|.
$$

For the $\l_{32}$ estimate, we stay in the $H_q^1$ norm but use our extra factor of $\ep$ to counterbalance the $J^{\ep}\varphib_{\ep}^a$ factor.  Specifically,
\begin{align*}
\ep\norm{\l_{32}^{\ep}(\etab,a)}_{1,q_{\star}} 
&\le C\ep|a|\norm{(J^{\ep}\sigmab).(J^{\ep}\varphib_{\ep}^a).\nl(\ep^2J^{\ep}\an_{\ep}(\etab,a))}_{1,q_{\star}} \\
\\
&\le C\ep|a|\norm{J^{\ep}\varphib_{\ep}^a}_{W^{1,\infty}}\norm{(J^{\ep}\sigmab).\nl(\ep^2J^{\ep}\an_{\ep}(\etab,a))}_{1,q_{\star}} \\
\\
&\le C\ep|a|\norm{(J^{\ep}\sigmab).(\ep J^{\ep}\an_{\ep}(\etab,a)).N(\ep^2J^{\ep}\an_{\ep}(\etab,a))}_{1,q_{\star}}.
\end{align*}
We finish by applying Proposition \ref{mapping} to the last term above in the $H_q^1$ norm and find
$$
\norm{\l_{32}^{\ep}(\etab,a)}_{1,q_{\star}} \le \M_{\l_{32}}[\norm{\etab}_{1,q_{\star}}]\ep|a|.
$$

\subsubsection{Mapping estimates for $\j_{42}$ and $\l_{42}$} These estimates are the same as those for $\j_{32}$ and $\l_{32}$, except all factors of $J^{\ep}\sigmab$ are replaced by $J^{\ep}\etab$.  

\subsubsection{Mapping estimates for $\j_{52}$ and $\l_{52}$} These estimates are the same as those for $\j_{12}$ and $\l_{12}$ with the factor of $(J^{\ep}\sigmab)^{.2}$ replaced by $(J^{\ep}\etab)^{.2}$.

\subsubsection{Mapping estimates for $\j_6$ and $\l_6$} It is for these terms that we designed the estimate in Proposition \ref{general hrq lipschitz prop}.  We begin with a straightforward estimate on $\j_6$:
\begin{align*}
\norm{\j_6^{\ep}(\etab,a)}_{1,q_{\star}}
&\le Ca^2\norm{\cb^{\ep}(\varphib_{\ep}^a,\varphib_{\ep}^a,\an_{\ep}(\etab,a))-\cb^{\ep}(\varphib_{\ep}^a,\varphib_{\ep}^a,a\varphib_{\ep}^a)}_{1,q_{\star}} \\
\\
&= Ca^2\norm{(J^{\ep}\varphib_{\ep}^a)^{.2}.(\nl(\ep^2J^{\ep}\an_{\ep}(\etab,a))-\nl(\ep^2J^{\ep}(a\varphib_{\ep}^a)))}_{1,q_{\star}} \\
\\
&\le Ca^2\norm{J^{\ep}\varphib_{\ep}^a}_{W^{1,\infty}}^2\norm{\nl(\ep^2J^{\ep}\an_{\ep}(\etab,a))-\nl(\ep^2J^{\ep}(a\varphib_{\ep}^a))}_{1,q_{\star}}.
\end{align*}
Since
$$
\nl(\ep^2J^{\ep}\an_{\ep}(\etab,a))-\nl(\ep^2J^{\ep}(a\varphib_{\ep}^a))
=
\nl(\ep^2J^{\ep}(\sigmab+\etab) + \ep^2J^{\ep}(a\varphib_{\ep}^a)) - \nl(0 + \ep^2J^{\ep}(a\varphib_{\ep}^a)),
$$
Proposition \ref{general hrq lipschitz prop} applies to give
$$
\norm{\nl(\ep^2J^{\ep}\an_{\ep}(\etab,a))-\nl(\ep^2J^{\ep}(a\varphib_{\ep}^a))}_{1,q_{\star}}
\le \M[\norm{\ep{J}^{\ep}(\sigmab+\etab)}_{1,q_{\star}}]\left(1+|a|\ep^{1-1}\right)\norm{\ep^2J^{\ep}(\sigmab+\etab)-0}_{1,q_{\star}}. 
$$
Then estimating $\norm{J^{\ep}\varphib_{\ep}^a}_{W^{1,\infty}} \le C\ep^{-2}$, we find
\begin{align*}
\norm{\j_6^{\ep}(\etab,a)}_{1,q_{\star}}
&\le Ca^2\ep^{-2}\M[\norm{\ep{J}^{\ep}(\sigmab+\etab)}_{1,q_{\star}}]\left(1+|a|\ep^{1-1}\right)\norm{\ep^2J^{\ep}(\sigmab+\etab)-0}_{1,q_{\star}} \\
\\
&\le Ca^2\M[\norm{\ep{J}^{\ep}(\sigmab+\etab)}_{1,q_{\star}}]\norm{J^{\ep}(\sigmab+\etab)}_{1,q_{\star}}
\end{align*}
Taking the supremum over $0 < \ep < \epbar$ gives a bound of the form
$$
\norm{\j_6^{\ep}(\etab,a)}_{1,q_{\star}} 
\le \M_{\j_6}[\norm{\etab}_{1,q_{\star}}]a^2.
$$
The estimate for $\l_6$ follows in the same way.

%% file: appendix_lipschitz_estimates.tex
\subsection{The Lipschitz estimates}\label{lipschitz estimates appendix}
Let
\begin{equation}\label{rhs lip defn}
\rhs_{\lip}^{\ep}(\etab,\grave{\etab},a,\grave{a}) 
= \left(\ep + \norm{(\etab,a)}_{\X^1} + \norm{(\grave{\etab},\grave{a})}_{\X^1}\right)\norm{(\etab,a)-(\grave{\etab},\grave{a})}_{\X^0}.
\end{equation}
We prove the existence of an increasing function $\M_{\lip} \colon \R_+^2 \to \R_+$ such that 
\begin{equation}\label{bootstrap statement appendix}
\norm{\nsystb^{\ep}(\etab,a)-\nsystb^{\ep}(\grave{\etab},\grave{a})}_{\X^0} 
\le \M_{\lip}[\norm{\etab}_{1,q_{\star}},\norm{\grave{\etab}}_{1,q_{\star}}]\rhs_{\lip}^{\ep}(\etab,\grave{\etab},a,\grave{a})
\end{equation}
for all $\etab,\grave{\etab} \in E_{q_{\star}}^1 \times E_{q_{\star}}^1$, $|a| \le a_{\per}$, and $0 < \ep < \overline{\ep}$.

\subsubsection{General strategy} A first pass using the estimates in Proposition \ref{faver wright nanopteron} gives
\begin{align*}
\norm{\nsyst_1^{\ep}(\etab,a) - \nsyst_1^{\ep}(\grave{\etab},\grave{a})}_{1,\qstar/2}
&\le C\sum_{k=1}^5 \bignorm{\j_{k1}^{\ep,\mod}(\etab,a)-\j_{k1}^{\ep,\mod}(\grave{\etab},\grave{a})}_{1,\qstar/2} \\
\\
&+ C\sum_{k=1}^5\ep\bignorm{\l_{k1}^{\ep,\mod}(\etab,a)-\l_{k1}^{\ep,\mod}(\grave{\etab},\grave{a})}_{1,\qstar/2} \\
\\
&+ C\sum_{k=1}^5 \norm{\j_{k2}^{\ep}(\etab,a)-\j_{k2}^{\ep}(\grave{\etab},\grave{a})}_{1,\qstar/2} \\
\\
&+ C\sum_{k=1}^5 \ep\norm{\l_{k2}^{\ep}(\etab,a)-\l_{k2}^{\ep}(\grave{\etab},\grave{a})}_{1,\qstar/2} \\
\\
&+ C\norm{\j_6^{\ep}(\etab,a)-\j_6^{\ep}(\grave{\etab},\grave{a})}_{1,q_{\star}/2} \\
\\
&+ C\ep\norm{\l_6^{\ep}(\etab,a) - \l_6^{\ep}(\grave{\etab},\grave{a})}_{1,q_{\star}/2}, \\
\\
\norm{\nsyst_2^{\ep}(\etab,a)-\nsyst_1^{\ep}(\grave{\etab},\grave{a})}_{1,q_{\star}/2}
&\le C \sum_{k=1}^5 \ep\bignorm{\l_{k1}^{\ep,\mod}(\etab,a)-\l_{k1}^{\ep,\mod}(\grave{\etab},\grave{a})}_{1,q_{\star}/2} \\
\\
&+ C\sum_{k=1}^5 \ep\norm{\l_{k2}^{\ep}(\etab,a)-\l_{k2}^{\ep}(\grave{\etab},\grave{a})}_{1,q_{\star}/2} \\
\\ 
&+ C\ep\norm{\l_6^{\ep}(\etab,a)-\l_6^{\ep}(\grave{\etab},\grave{a})}_{1,q_{\star}/2}, 
\end{align*}
and
\begin{align*}
|\nsyst_3^{\ep}(\etab,a)-\nsyst_3^{\ep}(\grave{\etab},\grave{a})|
&\le C \sum_{k=1}^5 \ep\bignorm{\l_{k1}^{\ep,\mod}(\etab,a)-\l_{k1}^{\ep,\mod}(\grave{\etab},\grave{a})}_{1,q_{\star}/2} \\
\\
&+ C\sum_{k=1}^5 \ep\norm{\l_{k2}^{\ep}(\etab,a)-\l_{k2}^{\ep}(\grave{\etab},\grave{a})}_{1,q_{\star}/2} \\
\\ 
&+ C\ep\norm{\l_6^{\ep}(\etab,a)-\l_6^{\ep}(\grave{\etab},\grave{a})}_{1,q_{\star}/2}.
\end{align*}
In \cite{faver-wright} we bounded the differences
$$
\bignorm{\j_{k1}^{\ep,\mod}(\etab,a)-\j_{k1}^{\ep,\mod}(\grave{\etab},\grave{a})}_{1,q_{\star}/2},
\quadword{and}
\bignorm{\l_{k1}^{\ep,\mod}(\etab,a)-\l_{k1}^{\ep,\mod}(\grave{\etab},\grave{a})}_{1,q_{\star}/2}, 
\ k = 1,\ldots,5
$$
by $\rhs_{\lip}^{\ep}(\etab,\grave{\etab},a,\grave{a})$.  We just need to show that 
$$
\norm{\j_{k2}^{\ep}(\etab,a)-\j_{k2}^{\ep}(\grave{\etab},\grave{a})}_{1,q_{\star}/2}, 
\qquad
\ep\norm{\l_{k2}^{\ep}(\etab,a)-\j_{k2}^{\ep}(\grave{\etab},\grave{a})}_{1,q_{\star}/2}, \ k = 1,\ldots,5,
$$
$$
\norm{\j_6^{\ep}(\etab,a) - \j_6^{\ep}(\grave{\etab},\grave{a})}_{1,q_{\star}/2},
\quadword{and}
\ep\norm{\l_6^{\ep}(\etab,a) - \l_6^{\ep}(\grave{\etab},\grave{a})}_{1,q_{\star}/2}
$$
are all bounded by $\M_{\lip}[\norm{\etab}_{1,q_{\star}},\norm{\grave{\etab}}_{1,q_{\star}}]\rhs_{\lip}^{\ep}(\etab,\grave{\etab},a,\grave{a})$.


We will use a particular consequence of the essential estimate \eqref{varphib lipschitz estimate} often enough that it is worthwhile to single it out here.

\begin{lemma}\label{Jphi w1infty lipschitz}
For each $r \ge 0$ there is $C_r > 0$ such that for all $\ep \in (0,\overline{\ep})$ and $|a|$, $|\grave{a}| \le a_{\per}$, we have
$$
\norm{\sech\left(\frac{q_{\star}}{2}\cdot\right)J^{\ep}(a\varphib_{\ep}^a-\grave{a}\varphib_{\ep}^{\grave{a}})}_{W^{r,\infty}}
\le C_r\ep^{-r}|a-\grave{a}|.
$$
\end{lemma}

\subsubsection{Lipschitz estimates for $\j_{12}$ and $\l_{12}$} As with the mapping estimates, we spell out this first estimate in detail to show our reliance on the general Lipschitz estimates of Appendices \ref{Hrq lipschitz section}, \ref{H1q lipschitz section}, and \ref{W1infty lipschitz section}.  We begin with
\begin{align*}
\norm{\j_{12}^{\ep}(\etab,a) - \j_{12}^{\ep}(\grave{\etab},\grave{a})}_{1,q_{\star}/2}
&\le C\norm{\cb^{\ep}(\sigmab,\sigmab,\an_{\ep}(\etab,a))-\cb^{\ep}(\sigmab,\sigmab,\an_{\ep}(\grave{\etab},\grave{a}))}_{1,q_{\star}/2} \\
\\
&\le C\norm{(J^{\ep}\sigmab)^{.2}.(\nl(\ep^2J^{\ep}\an_{\ep}(\etab,a))-\nl(\ep^2J^{\ep}\an_{\ep}(\grave{\etab},\grave{a})))}_{1,q_{\star}/2} \\
\\
&\le \bunderbrace{C\norm{(J^{\ep}\sigmab)^{.2}.(\nl(\ep^2J^{\ep}\an_{\ep}(\etab,a))-\nl(\ep^2J^{\ep}(\grave{\etab},a)))}_{1,q_{\star}/2}}{\Delta_1} \\
\\
&+ \bunderbrace{C\norm{(J^{\ep}\sigmab)^{.2}.(\nl(\ep^2J^{\ep}\an_{\ep}(\grave{\etab},a))-\nl(\ep^2J^{\ep}\an_{\ep}(\grave{\etab},\grave{a})))}_{1,q_{\star}/2}}{\Delta_2}
\end{align*}

Since
$$
\nl(\ep^2J^{\ep}\an_{\ep}(\etab,a))-\nl(\ep^2J^{\ep}(\grave{\etab},a))
= \nl\bigg(\ep\big(\bunderbrace{\ep{J}^{\ep}(\sigmab+\grave{\etab})}{H_{q_{\star}}^1} + a(\bunderbrace{\ep{J}^{\ep}\varphib_{\ep}^a}{W^{1,\infty}})\big)\bigg)
- \nl\bigg(\ep\big(\bunderbrace{\ep{J}^{\ep}(\sigmab+\grave{\etab})}{H_{q_{\star}}^1}+a(\bunderbrace{\ep{J}^{\ep}\varphib_{\ep}^a)}{W^{1,\infty}}\big)\bigg),
$$
we can use Proposition \ref{general hrq lipschitz prop} to bound
\begin{equation}\label{Delta1 J21 L21 lipschitz}
\Delta_1 
\le \M_1[\norm{\ep^2J^{\ep}(\sigmab+\etab)}_{1,q_{\star}},\norm{\ep^2J^{\ep}(\sigmab+\grave{\etab})}_{1,q_{\star}}]\left(1+|a|\ep^{1-1}\right)\norm{\ep^2J^{\ep}(\sigmab+\etab)-\ep^2J^{\ep}(\sigmab+\grave{\etab})}_{1,q_{\star}/2},
\end{equation}
where $\M_1 \colon \R_+^2 \to \R_+$ is radially increasing.  We have
$$
\norm{\ep^2J^{\ep}(\sigmab+\etab)}_{1,q_{\star}} 
\le C(1+\norm{\etab})_{1,q_{\star}},
$$
and so
$$
\M_1[\norm{\ep^2J^{\ep}(\sigmab+\etab)}_{1,q_{\star}},\norm{\ep^2J^{\ep}(\sigmab+\grave{\etab})}_{1,q_{\star}}](1+|a|)
\le \bunderbrace{2\M_1[C(1+\norm{\etab})_{1,q_{\star}},C(1+\norm{\grave{\etab}})_{1,q_{\star}}]}{\M_2[\norm{\etab}_{1,q_{\star}},\norm{\grave{\etab}}_{1,q_{\star}}]}.
$$
We see that $\M_2$ is also radially increasing.  Then
$$
\Delta_1 
\le \M_2[\norm{\etab}_{1,q_{\star}},\norm{\grave{\etab}}_{1,q_{\star}}]\ep^2\norm{\etab-\grave{\etab}}_{1,q_{\star}/2}.
$$

For $\Delta_2$, we first note
\begin{multline*}
(J^{\ep}\sigmab)^{.2}.(\nl(\ep^2J^{\ep}\an_{\ep}(\grave{\etab},a))-\nl(\ep^2J^{\ep}\an_{\ep}(\grave{\etab},\grave{a}))) \\
= \bunderbrace{(J^{\ep}\sigmab)^{.2}}{H_{q_{\star}}^1}.\bigg(\nl\big(\bunderbrace{\ep^2J^{\ep}(\sigmab+\grave{\etab})}{H_{q_{\star}}^1}+\bunderbrace{\ep^2J^{\ep}(a\varphib_{\ep}^a}{W^{1,\infty}})\big)-\nl\big(\bunderbrace{\ep^2J^{\ep}(\sigmab+\grave{\etab})}{H_{q_{\star}}^1}+\bunderbrace{\ep^2J^{\ep}(\grave{a}\varphib_{\ep}^{\grave{a}})}{W^{1,\infty}}\big)\bigg).
\end{multline*}
Then we can use the estimate of Proposition \ref{lipschitz on phi}:
\begin{multline*}
\Delta_2 \le \M_3[\norm{\ep^2J^{\ep}(\sigmab+\grave{\etab})}_{1,q_{\star}},\norm{\ep^2J^{\ep}(a\varphib_{\ep}^a)}_{W^{1,\infty}},\norm{\ep^2J^{\ep}(\grave{a}\varphib_{\ep}^{\grave{a}})}_{W^{1,\infty}}] \\
\times \ep^2\norm{\sech\left(\frac{q_{\star}}{2}\cdot\right)(aJ^{\ep}\varphib_{\ep}^a-\grave{a}J^{\ep}\varphib_{\ep}^{\grave{a}})}_{W^{1,\infty}}\norm{(J^{\ep}\sigmab)^{.2}}_{1,q_{\star}},
\end{multline*}
where $\M_3 \colon \R_+^3 \to \R_+$ is radially increasing.  Taking the supremum over $\ep$, $a$, and $\grave{a}$ and using the estimate in Lemma \ref{Jphi w1infty lipschitz} on the $W^{1,\infty}$-factor, we find
$$
\Delta_2 \le \M_{4}[\norm{\grave{\etab}}_{1,q_{\star}}]\ep|a-\grave{a}|.
$$
for an increasing function $\M_4$.  We conclude that the sum $\Delta_1 + \Delta_2$ has an upper bound of the form $\rhs_{\lip}^{\ep}(\etab,\grave{\etab},a,\grave{a})$ from \eqref{rhs lip defn}, and the estimate for $\l_{12}$ is the same.

\subsubsection{Lipschitz estimates for $\j_{22}$ and $\l_{22}$} These estimates are mostly the same as the ones for $\j_{12}$ and $\l_{12}$ with a few small changes that are worth pointing out.  We estimate $\j_{22}$ to illustrate them:
\begin{align*}
\norm{\j_{22}^{\ep}(\etab,a) - \j_{22}^{\ep}(\grave{\etab},\grave{a})}_{1,q_{\star}/2}
&\le C\norm{\cb^{\ep}(\sigmab,\etab,\an_{\ep}(\etab,a)) - \cb^{\ep}(\sigmab,\grave{\etab},\an_{\ep}(\grave{\etab},\grave{a}))}_{1,q_{\star}/2} \\
\\
&\le \bunderbrace{\norm{\cb^{\ep}(\sigmab,\etab,\an_{\ep}(\etab,a))-\cb^{\ep}(\sigmab,\grave{\etab},\an_{\ep}(\etab,a))}_{1,q_{\star}/2}}{\Delta_1} \\
\\
&+ \bunderbrace{\norm{\cb^{\ep}(\sigmab,\grave{\etab},\an_{\ep}(\etab,a))-\cb^{\ep}(\sigmab,\grave{\etab},\an_{\ep}(\grave{\etab},a))}_{1,q_{\star}/2}}{\Delta_2} \\
\\
&+ \bunderbrace{\norm{\cb^{\ep}(\sigmab,\grave{\etab},\an_{\ep}(\grave{\etab},a))-\cb^{\ep}(\sigmab,\grave{\etab},\an_{\ep}(\grave{\etab},\grave{a}))}_{1,q_{\star}/2}}{\Delta_3}
\end{align*}
We have
$$
\Delta_1 \le C\norm{(J^{\ep}\sigmab).(J^{\ep}(\etab-\grave{\etab})).(\ep^2J^{\ep}\an_{\ep}(\etab,a)).N(\ep^2J^{\ep}\an_{\ep}(\etab,a))}_{1,q_{\star}/2},
$$
and using the algebra property of $H_q^1 \times H_q^1$, this factors as
$$
\Delta_1 \le C\ep\norm{J^{\ep}(\etab-\grave{\etab})}_{1,q_{\star}/2}\bunderbrace{\norm{(J^{\ep}\sigmab).(\ep{J}^{\ep}\an_{\ep}(\etab,a)).N(\ep^2J^{\ep}\an_{\ep}(\etab,a))}_{1,q_{\star}/2}}{\Pi}.
$$
We can bound the factor $\Pi$ using Proposition \ref{mapping prop} effectively as we did in Appendix \ref{J12 L12 mapping estimates} for the mapping estimates on $\j_{12}$ and $\l_{12}$.  Then
$$
\Delta_1 \le \M[\norm{\etab}_{1,q_{\star}}]\ep\norm{\etab-\grave{\etab}}_{1,q_{\star}/2}.
$$
The estimate for $\Delta_2$ again uses the algebra property of $H_q^1 \times H_q^1$ and then Proposition \ref{general hrq lipschitz prop} as we did in \eqref{Delta1 J21 L21 lipschitz} for the Lipschitz estimates on $\j_{21}$ and $\l_{21}$.  Finally, the estimate for $\Delta_3$ uses Proposition \ref{lipschitz on phi}. 

The estimate for $\l_{22}$ is identical.

\subsubsection{Lipschitz estimates for $\j_{32}$ and $\l_{33}$}\label{J32 L32 lipschitz}
As with the mapping estimates, we need to exploit the smoothing operator $\varpi^{\ep}$ on $\j_{32}$ and the extra $\ep$ on $\l_{33}$ to manage the presence of $\varphib_{\ep}^a$ and $\varphib_{\ep}^{\grave{a}}$.  We have
\begin{align*}
\norm{\j_{32}^{\ep}(\etab,a)-\j_{32}^{\ep}(\grave{\etab},\grave{a})}_{1,q_{\star}/2} 
&= 2\norm{\varpi^{\ep}(a\cb_2^{\ep}(\sigmab,\varphib_{\ep}^a,\an_{\ep}(\etab,a)) - \grave{a}\cb_2^{\ep}(\sigmab,\varphib_{\ep}^{\grave{a}},\an_{\ep}(\grave{\etab},\grave{a})))}_{1,q_{\star}/2} \\
\\
&\le C\norm{a\cb^{\ep}(\sigmab,\varphib_{\ep}^a,\an_{\ep}(\etab,a)) - \grave{a}\cb^{\ep}(\sigmab,\varphib_{\ep}^{\grave{a}},\an_{\ep}(\grave{\etab},\grave{a}))}_{0,q_{\star}/2} \\
\\
&\le \bunderbrace{|a-\grave{a}|\norm{\cb^{\ep}(\sigmab,\varphib_{\ep}^a,\an_{\ep}(\etab,a)}_{0,q_{\star}/2}}{\Delta_1} \\
\\
&+ \bunderbrace{|\grave{a}|\norm{\cb^{\ep}(\sigmab,\varphib_{\ep}^a,\an_{\ep}(\etab,a)-\cb^{\ep}(\sigmab,\varphib_{\ep}^{\grave{a}},\an_{\ep}(\etab,a))}_{0,q_{\star}/2}}{\Delta_2} \\
\\
&+ \bunderbrace{|\grave{a}|\norm{\cb^{\ep}(\sigmab,\varphib_{\ep}^{\grave{a}},\an_{\ep}(\etab,a))-\cb^{\ep}(\sigmab,\varphib_{\ep}^{\grave{a}},\an_{\ep}(\grave{\etab},a))}_{0,q_{\star}/2}}{\Delta_3} \\
\\
&+ \bunderbrace{|\grave{a}|\norm{\cb^{\ep}(\sigmab,\varphib_{\ep}^{\grave{a}},\an_{\ep}(\grave{\etab},a))-\cb^{\ep}(\sigmab,\varphib_{\ep}^{\grave{a}},\an_{\ep}(\grave{\etab},\grave{a}))}_{0,q_{\star}/2}}{\Delta_4}.
\end{align*}
We handle $\Delta_1$ easily using the mapping estimate of Proposition \ref{mapping prop}:
\begin{align*}
\Delta_1 &\le C\ep|a-\grave{a}|\norm{(J^{\ep}\sigmab).(J^{\ep}\varphib_{\ep}^a).(\ep{J}^{\ep}\an_{\ep}(\etab,a)).N(\ep^2J^{\ep}\an_{\ep}(\etab,a))}_{0,q_{\star}/2} \\
\\
&\le C\ep|a-\grave{a}|\norm{J^{\ep}\varphib_{\ep}^a}_{W^{0,\infty}}\norm{(J^{\ep}\sigmab).(\ep{J}^{\ep}\an_{\ep}(\etab,a)).N(\ep^2J^{\ep}\an_{\ep}(\etab,a))}_{0,q_{\star}/2} \\
\\
&\le \M_{\Delta_1}[\norm{\etab}_{1,q_{\star}}]\ep|a-\grave{a}|.
\end{align*}
For $\Delta_2$ we need both Proposition \ref{mapping prop} and Lemma \ref{Jphi w1infty lipschitz}:
\begin{align*}
\Delta_2
&\le C|\grave{a}|\norm{(J^{\ep}\sigmab).(J^{\ep}(\varphib_{\ep}^a-\varphib_{\ep}^{\grave{a}})).\nl(\ep^2\an_{\ep}(\etab,a))}_{0,q_{\star}/2} \\
\\
&\le C\ep|\grave{a}|\norm{J^{\ep}(\varphib_{\ep}^a-\varphib_{\ep}^{\grave{a}})}_{W^{0,\infty}}\norm{(J^{\ep}\sigmab).(\ep{J}^{\ep}\an_{\ep}(\etab,a)).N(\ep^2J^{\ep}\an_{\ep}(\etab,a))}_{0,q_{\star}/2} \\
\\
&\le \M_{\Delta_2}[\norm{\etab}_{1,q_{\star}}]\ep|\grave{a}||a-\grave{a}|.
\end{align*}
We factor $\Delta_3$ and then employ Proposition \ref{general hrq lipschitz prop}:
\begin{align*}
\Delta_3 
&\le C|\grave{a}|\norm{(J^{\ep}\sigmab).(J^{\ep}\varphib_{\ep}^{\grave{a}}).(\nl(\ep^2J^{\ep}\an_{\ep}(\etab,a))-\nl(\ep^2J^{\ep}\an_{\ep}(\grave{\etab},a)))}_{0,q_{\star}/2} \\
\\
&\le C|\grave{a}|\norm{J^{\ep}\varphib_{\ep}^a}_{W^{0,\infty}}\norm{J^{\ep}\sigmab}_{W^{0,\infty}}\norm{\nl(\ep^2J^{\ep}\an_{\ep}(\etab,a))-\nl(\ep^2J^{\ep}\an_{\ep}(\grave{\etab},a))}_{0,q_{\star}/2} \\
\\
&\le \M_{\Delta_3}[\norm{\etab}_{1,q_{\star}},\norm{\grave{\etab}}_{1,q_{\star}}]\ep^2|\grave{a}|\norm{\etab-\grave{\etab}}_{0,q_{\star}/2}.
\end{align*}
Finally, $\Delta_4$ uses Proposition \ref{lipschitz on phi}:
\begin{align*}
\Delta_4 
&\le C|\grave{a}|\norm{(J^{\ep}\sigmab).(J^{\ep}\varphib_{\ep}^a).(\nl(\ep^2J^{\ep}\an_{\ep}(\grave{\etab},a))-\nl(\ep^2J^{\ep}\an_{\ep}(\grave{\etab},\grave{a})))}_{0,q_{\star}/2} \\
\\
&\le C|\grave{a}|\norm{J^{\ep}\varphib_{\ep}^a}_{W^{0,\infty}}\norm{(J^{\ep}\sigmab).(\nl(\ep^2J^{\ep}\an_{\ep}(\grave{\etab},a))-\nl(\ep^2J^{\ep}\an_{\ep}(\grave{\etab},\grave{a})))}_{0,q_{\star}/2} \\
\\
&\le \M[\norm{\grave{\etab}}_{1,q_{\star}}]\ep^2|\grave{a}|\norm{\sech\left(\frac{q_{\star}}{2}\cdot\right)J^{\ep}(a\varphib_{\ep}^a-\grave{a}\varphib_{\ep}^{\grave{a}})}_{W^{0,\infty}}\norm{J^{\ep}\sigmab}_{1,q_{\star}} \\
\\
&\le \M_{\Delta_4}[\norm{\grave{\etab}}_{1,q_{\star}}]\ep^2|\grave{a}||a-\grave{a}|.
\end{align*}
\subsubsection{Lipschitz estimates for $\j_{42}$ and $\l_{42}$} Again, we smooth with $\varpi^{\ep}$ on $\j_{42}$ and use the extra $\ep$ on $\l_{42}$ to our advantage; the mechanics are the same as the estimates for $\j_{32}$ and $\l_{32}$ with one exception, which we highlight below.  We bound
\begin{align*}
\norm{\j_{42}^{\ep}(\etab,a)-\j_{42}^{\ep}(\grave{\etab},\grave{a})}_{1,q_{\star}/2}
&\le C\norm{a\cb^{\ep}(\etab,\varphib_{\ep}^a,\an_{\ep}(\etab,a))-\grave{a}\cb^{\ep}(\grave{\etab},\varphib_{\ep}^{\grave{a}},\an_{\ep}(\grave{\etab},\grave{a}))}_{0,q_{\star}/2} \\
\\
&\le |a-\grave{a}|\norm{\cb^{\ep}(\etab,\varphib_{\ep}^a,\an_{\ep}(\etab,a))}_{0,q_{\star}/2} \\
\\
&\bunderbrace{+|\grave{a}|\norm{\cb^{\ep}(\etab,\varphib_{\ep}^a,\an_{\ep}(\etab,a))-\cb^{\ep}(\grave{\etab},\varphib_{\ep}^a,\an_{\ep}(\etab,a))}_{0,q_{\star}/2}}{\Delta} \\
\\
&+ |\grave{a}|\norm{\cb^{\ep}(\grave{\etab},\varphib_{\ep}^a,\an_{\ep}(\etab,a))-\cb^{\ep}(\grave{\etab},\varphib_{\ep}^{\grave{a}},\an_{\ep}(\etab,a))}_{0,q_{\star}/2} \\
\\
&+ |\grave{a}|\norm{\cb^{\ep}(\grave{\etab},\varphib_{\ep}^{\grave{a}},\an_{\ep}(\etab,a))-\cb^{\ep}(\grave{\etab},\varphib_{\ep}^{\grave{a}},\an_{\ep}(\grave{\etab},a))}_{0,q_{\star}/2} \\
\\
&+ |\grave{a}|\norm{\cb^{\ep}(\grave{\etab},\varphib_{\ep}^{\grave{a}},\an_{\ep}(\grave{\etab},a))-\cb^{\ep}(\grave{\etab},\varphib_{\ep}^{\grave{a}},\an_{\ep}(\grave{\etab},\grave{a}))}_{0,q_{\star}/2}
\end{align*}
The first term above can be bounded with our standard mapping estimate from Proposition \ref{mapping prop}.  The last three terms above are analogous to $\Delta_2$, $\Delta_3$, and $\Delta_4$ in Appendix \ref{J32 L32 lipschitz}.  The exception is the term that we have labeled $\Delta$ above:
\begin{align*}
\Delta
&\le C|\grave{a}|\norm{(J^{\ep}(\etab-\grave{\etab})).(J^{\ep}\varphib_{\ep}^a).\nl(\ep^2J^{\ep}\an_{\ep}(\etab,a))}_{0,q_{\star}/2} \\
\\
&\le C\ep|\grave{a}|\norm{J^{\ep}\varphib_{\ep}^a}_{W^{0,\infty}}\norm{\ep{J}^{\ep}\an_{\ep}(\etab,a)}_{W^{0,\infty}}\norm{(J^{\ep}(\etab-\grave{\etab})).N(\ep^2J^{\ep}\an_{\ep}(\etab,a))}_{0,q_{\star}/2} \\
\\
&\le C\ep|\grave{a}|\norm{(J^{\ep}(\etab-\grave{\etab})).N(\ep^2J^{\ep}\an_{\ep}(\etab,a))}_{0,q_{\star}/2} \\
\\
&\le C\ep|\grave{a}|\M[\norm{\etab}_{1,q_{\star}}]\norm{\etab-\grave{\etab}}_{0,q_{\star}/2}.
\end{align*}
For the last inequality we have used Proposition \ref{easy peasy lipschitz}.  The estimate for $\l_{42}$ goes through in the same way.

\subsubsection{Lipschitz estimates for $\j_{52}$ and $\l_{52}$} As we have no factor of $\varphib_{\ep}^a$ here, we do not need to use smoothing or an extra factor of $\ep$ to avoid problems. We bound
\begin{align*}
\norm{\j_{52}^{\ep}(\etab,a)-\j_{52}^{\ep}(\grave{\etab},\grave{a})}_{1,q_{\star}/2} 
&\le C\norm{\cb^{\ep}(\etab,\etab,\an_{\ep}(\etab,a))-\cb^{\ep}(\grave{\etab},\grave{\etab},\an_{\ep}(\grave{\etab},\grave{a}))}_{1,q_{\star}/2} \\
\\
&\le \bunderbrace{C\norm{\cb^{\ep}(\etab,\etab,\an_{\ep}(\etab,a))-\cb^{\ep}(\grave{\etab},\grave{\etab},\an_{\ep}(\etab,a))}_{1,q_{\star}/2}}{\Delta} \\
\\
&+ C\norm{\cb^{\ep}(\grave{\etab},\grave{\etab},\an_{\ep}(\etab,a))-\cb^{\ep}(\grave{\etab},\grave{\etab},\an_{\ep}(\grave{\etab},a))}_{1,q_{\star}/2} \\
\\
&+ C\norm{\cb^{\ep}(\grave{\etab},\grave{\etab},\an_{\ep}(\grave{\etab},a))-\cb^{\ep}(\grave{\etab},\grave{\etab},\an_{\ep}(\grave{\etab},\grave{a}))}_{1,q_{\star}/2}.
\end{align*} 
Of the three terms above, we know how to estimate the second and third using Propositions \ref{general hrq lipschitz prop} and \ref{lipschitz on phi}; we bound $\Delta$ by
\begin{align*}
\Delta &\le C\norm{((J^{\ep}\etab)^{.2}-(J^{\ep}\grave{\etab})^{.2}).\an_{\ep}(\etab,a)}_{1,q_{\star}/2} \\
\\
&\le C\ep\norm{J^{\ep}(\etab+\grave{\etab})}_{1,q_{\star}/2}\norm{(J^{\ep}(\etab-\grave{\etab})).(\ep{J}^{\ep}\an_{\ep}(\etab,a)).N(\ep^2J^{\ep}\an_{\ep}(\etab,a))}_{1,q_{\star}/2} \\
\\
&\le \M[\norm{\etab}_{1,q_{\star}},\norm{\grave{\etab}}_{1,q_{\star}}]\ep\norm{\etab-\grave{\etab}}_{1,q_{\star}/2}.
\end{align*}
For the second inequality, we factored the difference of squares and used the algebra property of $H_{q_{\star}/2}^1$ and for the third we used Proposition \ref{easy peasy lipschitz}. The estimate for $\l_{52}$ is the same.

\subsubsection{Lipschitz estimates for $\j_6$ and $\l_6$} We work on $\j_6$; the strategy for $\l_6$ is the same.  We have
$$
\norm{\j_6^{\ep}(\etab,a)-\j_6^{\ep}(\grave{\etab},\grave{a})}_{1,q_{\star}/2}
\le C\sum_{k=1}^5 \norm{\Delta_k}_{1,q_{\star}/2},
$$
where 
\begin{align*}
\Delta_1 &
= (a^2-\grave{a}^2)(\cb^{\ep}(\varphib_{\ep}^a,\varphib_{\ep}^a,\an_{\ep}(\etab,a))-\cb^{\ep}(\varphib_{\ep}^a,\varphib_{\ep}^a,a\varphib_{\ep}^a)) \\
\\
\Delta_2 
&= \grave{a}^2(\cb^{\ep}(\varphib_{\ep}^a,\varphib_{\ep}^a,\an_{\ep}(\etab,a))-\cb^{\ep}(\varphib_{\ep}^a,\varphib_{\ep}^a,a\varphib_{\ep}^a) 
- (\cb^{\ep}(\varphib_{\ep}^{\grave{a}},\varphib_{\ep}^a,\an_{\ep}(\etab,a))-\cb^{\ep}(\varphib_{\ep}^{\grave{a}},\varphib_{\ep}^a,a\varphib_{\ep}^a))) \\
\\
\Delta_3 
&= \grave{a}^2((\cb^{\ep}(\varphib_{\ep}^{\grave{a}},\varphib_{\ep}^a,\an_{\ep}(\etab,a))-\cb^{\ep}(\varphib_{\ep}^{\grave{a}},\varphib_{\ep}^a,a\varphib_{\ep}^a)))
-(\cb^{\ep}(\varphib_{\ep}^{\grave{a}},\varphib_{\ep}^{\grave{a}},\an_{\ep}(\etab,a))-\cb^{\ep}(\varphib_{\ep}^{\grave{a}},\varphib_{\ep}^{\grave{a}},a\varphib_{\ep}^a)) \\
\\
\Delta_4
&= \grave{a}^2((\cb^{\ep}(\varphib_{\ep}^{\grave{a}},\varphib_{\ep}^{\grave{a}},\an_{\ep}(\etab,a))-\cb^{\ep}(\varphib_{\ep}^{\grave{a}},\varphib_{\ep}^{\grave{a}},a\varphib_{\ep}^a))
-(\cb^{\ep}(\varphib_{\ep}^{\grave{a}},\varphib_{\ep}^{\grave{a}},\an_{\ep}(\grave{\etab},a))-\cb^{\ep}(\varphib_{\ep}^{\grave{a}},\varphib_{\ep}^{\grave{a}},a\varphib_{\ep}^a))) \\
\\
\Delta_5
&= \grave{a}^2((\cb^{\ep}(\varphib_{\ep}^{\grave{a}},\varphib_{\ep}^{\grave{a}},\an_{\ep}(\grave{\etab},a))-\cb^{\ep}(\varphib_{\ep}^{\grave{a}},\varphib_{\ep}^{\grave{a}},a\varphib_{\ep}^a))
- (\cb^{\ep}(\varphib_{\ep}^{\grave{a}},\varphib_{\ep}^{\grave{a}},\an_{\ep}(\grave{\etab},\grave{a}))-\cb^{\ep}(\varphib_{\ep}^{\grave{a}},\varphib_{\ep}^{\grave{a}},\grave{a}\varphib_{\ep}^{\grave{a}})))
\end{align*}

Proposition \ref{general hrq lipschitz prop} lets us bound $\Delta_1$ as
$$
\norm{\Delta_1}_{1,q_{\star}/2} \le \M[\norm{\etab}_{1,q_{\star}}]\ep|a-\grave{a}|.
$$
We use this proposition again on $\Delta_2$:
\begin{align*}
\norm{\Delta_2}_{1,q_{\star}/2} 
& \le C\grave{a}^2\norm{(J^{\ep}(\varphib_{\ep}^a-\varphib_{\ep}^{\grave{a}})).(J^{\ep}\varphib_{\ep}^a).(\nl(\ep^2J^{\ep}\an_{\ep}(\etab,a))-\nl(\ep^2J^{\ep}(a\varphib_{\ep}^a)))}_{1,q_{\star}/2} \\
\\
&\le C\grave{a}^2\norm{J^{\ep}(\varphib_{\ep}^a-\varphib_{\ep}^{\grave{a}})}_{W^{1,\infty}}\norm{J^{\ep}\varphib_{\ep}^a}_{W^{1,\infty}}\norm{\nl(\ep^2J^{\ep}\an_{\ep}(\etab,a))-\nl(\ep^2J^{\ep}(a\varphib_{\ep}^a)))}_{1,q_{\star}/2} \\
\\
&\le \grave{a}^2\ep^{-2}|a-\grave{a}|\M[\norm{\etab}_{1,q_{\star}}]\ep^2 \\
\\
&= \M[\norm{\etab}_{1,q_{\star}}]|\grave{a}|^2|a-\grave{a}| \\
\\
&\le \M[\norm{\etab}_{1,q_{\star}}]|a-\grave{a}|
\end{align*}
since $|\grave{a}| \le a_{\per} < 1$.

The estimate for $\Delta_3$ is exactly the same as the one for $\Delta_2$, while for $\Delta_4$ the $\pm\cb^{\ep}(\varphib_{\ep}^{\grave{a}},\varphib_{\ep}^{\grave{a}},a\varphib_{\ep}^a)$ terms nicely cancel to give us
\begin{align*}
\norm{\Delta_4}_{1,q_{\star}/2}
&\le C\grave{a}^2\norm{(J^{\ep}\varphib_{\ep}^{\grave{a}})^{.2}.(\nl(\ep^2J^{\ep}\an_{\ep}(\etab,a))-\nl(\ep^2J^{\ep}\an_{\ep}(\grave{\etab},a)))}_{1,q_{\star}/2} \\
\\
&\le \grave{a}^2\ep^{-2}\M[\norm{\etab}_{1,q_{\star}},\norm{\grave{\etab}}_{1,q_{\star}}]\ep^2\norm{\etab-\grave{\etab}}_{1,q_{\star}/2} \\
\\
&\le \M[\norm{\etab}_{1,q_{\star}},\norm{\grave{\etab}}_{1,q_{\star}}]|\grave{a}|\norm{\etab-\grave{\etab}}_{1,q_{\star}/2}.
\end{align*}

After factoring out $(J^{\ep}\varphib_{\ep}^{\grave{a}})^{.2}$, we rewrite $\Delta_5$ in a form amenable to Proposition \ref{use the ftc, luke}:
\begin{align*}
\norm{\Delta_5}_{1,q_{\star}/2}
&\le C\grave{a}^2\ep^{-2}\norm{(\nl(\ep^2J^{\ep}\an_{\ep}(\grave{\etab},a))-\nl(\ep^2J^{\ep}(a\varphib_{\ep}^a)))-(\nl(\ep^2J^{\ep}\an_{\ep}(\grave{\etab},\grave{a}))-\nl(\ep^2J^{\ep}(\grave{a}\varphib_{\ep}^{\grave{a}})))}_{1,q_{\star}/2} \\
\\
&\le \M[\norm{\grave{\etab}}_{1,q_{\star}}]|\grave{a}|\ep^{-2}\norm{\sech\left(\frac{q_{\star}}{2}\cdot\right)\ep^2(aJ^{\ep}\varphib_{\ep}^a-\grave{a}J^{\ep}\varphib_{\ep}^{\grave{a}})}_{W^{1,\infty}}\norm{\ep^2J^{\ep}(\sigmab+\grave{\etab})}_{1,q_{\star}} \\
\\
&\le \M[\norm{\grave{\etab}}_{1,q_{\star}}]|\grave{a}|\ep|a-\grave{a}|.
\end{align*}

%% file: appendix_bootstrap_estimates.tex
\subsection{The bootstrap estimates}\label{bootstrap estimates appendix}
Let
\begin{equation}\label{rhs boot defn}
\rhs_{\boot,r}^{\ep}(\etab,a)
:= \ep 
+ |a|\ep^{1-r}
+ a^2\ep^{1-2r}
+ |a|^3\ep^{1-3r}
+ \norm{\etab}_{r,\star}
+ |a|\norm{\etab}_{r,\qstar}\ep^{-r}.
\end{equation}
We prove the existence of increasing functions $\M_{\boot,r} \colon \R_+ \to \R_+$ such that 
\begin{equation}\label{bootstrap statement appendix 1}
\norm{(\nsyst_1^{\ep}(\etab,a),\nsyst_2^{\ep}(\etab,a))}_{r+1,\qstar}
\le \M_{\boot,r}[\norm{\etab}_{r,\qstar}]\rhs_{\boot,r}^{\ep}(\etab,a)
\end{equation}
and
\begin{equation}\label{bootstrap statement appendix 2}
|\nsyst_3^{\ep}(\etab,a)| 
\le \M_{\boot,r}[\norm{\etab}_{r,\qstar}]\ep^r\rhs_{\boot,r}^{\ep}(\etab,a)
\end{equation}
for all $\etab \in \X^r$ $|a| \le a_{\per}$, and $0 < \ep < \epbar$.

\subsubsection{General strategy} 
It is for the sake of these bootstrap estimates that we proved Propositions \ref{mapping prop} and \ref{general hrq lipschitz prop} for arbitrary $r$ (thereby complicating the proofs considerably as opposed to doing them for just $r=1$).  In these sections $r$ will always be a positive integer.

The smoothing property \eqref{Pep estimate} of $\P_{\ep}$ gives 
$$
\norm{\nsyst_2^{\ep}(\etab,a)}_{r+1,\qstar}
\le C_r\sum_{k=1}^5 \bignorm{\l_{k1}^{\ep,\mod}(\etab,a)}_{r+1,\qstar}
+ C_r\sum_{k=1}^5\norm{\l_{k2}^{\ep}(\etab,a)}_{r+1,\qstar} 
+ \norm{\l_6^{\ep}(\etab,a)}_{r+1,\qstar},
$$
and this along with the boundedness of $\A^{-1}$ implies
\begin{align*}
\norm{\nsyst_1^{\ep}(\etab,a)}_{r+1,\qstar}
&\le C_r\sum_{k=1}^5 \left(\bignorm{\j_{k1}^{\ep,\mod}(\etab,a)}_{r+1,\qstar}
+ \bignorm{\l_{k1}^{\ep,\mod}(\etab,a)}_{r,\qstar}\right) \\
\\
&+ C_r\sum_{k=1}^5 \left(\norm{\j_{k2}^{\ep}(\etab,a)}_{r+1,\qstar}
+ \norm{\l_{k2}^{\ep}(\etab,a)}_{r,\qstar}\right) \\
\\
&+C_r\left(\norm{\j_6^{\ep}(\etab,a)}_{r+1,\qstar} + \norm{\l_6^{\ep}(\etab,a)}_{r,\qstar}\right).
\end{align*}
Last, the estimate \eqref{iota ep estimate} on $\iota_{\ep}$ and the boundedness of $\upsilon_{\ep}$ from \eqref{upsilon ep estimate} give
$$
|\nsyst_3^{\ep}(\etab,a)| 
\le C\ep^r\sum_{k=1}^5\left(\bignorm{\l_{k1}^{\ep,\mod}(\etab,a)}_{r,\qstar} 
+ \norm{\l_{k2}^{\ep}(\etab,a)}_{r,\qstar}\right) + C\ep^r\norm{\l_6^{\ep}(\etab,a)}_{r,\qstar}. 
$$
In \cite{faver-wright} we saw that the terms 
$$
\bignorm{\j_{1k}^{\ep,\mod}(\etab,a)}_{r+1,\qstar}
\quadword{and}
\bignorm{\l_{1k}^{\ep,\mod}(\etab,a)}_{r,\qstar},
\ k = 1,\ldots,5
$$
are all bounded by $\rhs_{\boot,r}^{\ep}(\etab,a)$.  Now we bound the remaining terms
$$
\norm{\j_{2k}^{\ep}(\etab,a)}_{r+1,\qstar},
\qquad
\norm{\l_{2k}^{\ep}(\etab,a)}_{r,\qstar},
\ k = 1,\ldots,5,
\qquad
\norm{\j_6^{\ep}(\etab,a)}_{r+1,\qstar},
\quadword{and} 
\norm{\l_6^{\ep}(\etab,a)}_{r,\qstar}
$$
by $\M_{\boot,r}[\norm{\etab}_{r,\qstar}]\rhs_{\boot,r}^{\ep}(\etab,a)$.  

We will only show the estimates for the $\j$ terms, as once we have smoothed by $\varpi^{\ep}$ on a $\j$ term, the resulting upper bound is a constant multiple of the upper bound for the corresponding $\l$ term.

\subsubsection{Bootstrap estimates for $\j_{12}$} As with the mapping and Lipschitz estimates, we write this section in particular detail to illustrate our techniques. Smoothing by $\varpi^{\ep}$ per \eqref{varpi ep estimate}, we have
\begin{align*}
\norm{\j_{12}^{\ep}(\etab,a)}_{r+1,\qstar}
&\le C\norm{\cb^{\ep}(\sigmab,\sigmab,\an_{\ep}(\etab,a))}_{r,\qstar} \\
\\
&\le C\norm{(J^{\ep}\sigmab)^{.2}.\nl(\ep^2J^{\ep}\an_{\ep}(\etab,a))}_{r,\qstar} \\
\\
&\le \bunderbrace{C\ep^2\norm{(J^{\ep}\sigmab)^{.2}.(J^{\ep}(\sigmab+\etab)).N(\ep^2J^{\ep}\an_{\ep}(\etab,a))}_{r,\qstar}}{\Pi_1} \\
\\
&+ \bunderbrace{C\ep^2\norm{(J^{\ep}\sigmab)^{.2}.(aJ^{\ep}\varphib_{\ep}^a).N(\ep^2J^{\ep}\an_{\ep}(\etab,a))}_{r,\qstar}}{\Pi_2}.
\end{align*}
Then Proposition \ref{mapping prop} with $C_{\star} = \max\{C_1,\ldots,C_r\}$ and $C_1,\ldots,C_r$ satisfying the estimate \eqref{varphib Wrinfty estimate} implies 
$$
\Pi_1 
\le \ep^2\M_{1,r}[\ep^2\norm{J^{\ep}(\sigmab+\etab)}_{r,\qstar}]\left(1+|a|\ep^{1-r}\right)\norm{(J^{\ep}\sigmab)^{.2}.(J^{\ep}(\sigmab+\etab))}_{r,\qstar}.
$$
After factoring, the same proposition gives
\begin{align*}
\Pi_2
&\le \ep^2|a|\norm{J^{\ep}\varphib_{\ep}^a}_{W^{r,\infty}}\norm{(J^{\ep}\sigmab)^{.2}.N(\ep^2J^{\ep}\an_{\ep}(\etab,a))}_{r,\qstar} \\
\\
&\le C_r\ep^{2-r}|a|\M_{2,r}[\norm{\ep^2J^{\ep}(\sigmab+\etab)}_{r,\qstar}]\left(1+|a|\ep^{1-r}\right)\norm{(J^{\ep}\sigmab)^{.2}}_{r,\qstar}.
\end{align*}
Since we can assume $\M_{1,r}$ and $\M_{2,r}$ are increasing, we take the supremum over $0 < \ep < \epbar < 1$ and find
\begin{equation}\label{pi1 bootstrap}
\Pi_1
\le \tilde{\M}_{1,r}[\norm{\etab}_{r,\qstar}]\ep\left(1+|a|\ep^{1-r}\right)
\end{equation}
and
\begin{equation}\label{pi2 bootstrap}
\Pi_2
\le \tilde{\M}_{2,r}[\norm{\etab}_{r,\qstar}]\ep^{2-r}|a|\left(1+|a|\ep^{1-r}\right)
\le \tilde{\M}_{2,r}[\norm{\etab}_{r,\qstar}]|a|\ep^{1-r}+a^2\ep^{1-2r}
\end{equation}
for increasing functions $\tilde{\M}_{1,r}$ and $\tilde{\M}_{2,r}$.  
All together, these estimates give an upper bound on $\norm{\j_{12}^{\ep}(\etab,a)}_{r+1,\qstar}$ of the form $\M[\norm{\etab}_{r,\qstar}]\rhs_{\boot,r}^{\ep}(\etab,a)$ given in \eqref{rhs boot defn}.

\subsubsection{Bootstrap estimates for $\j_{22}$} We have
\begin{align*}
\norm{\j_{22}^{\ep}(\etab,a)}_{r+1,\qstar}
&\le C\norm{\cb^{\ep}(\sigmab,\etab,\an_{\ep}(\etab,a))}_{r,\qstar} \\
\\
&\le C\norm{(J^{\ep}\sigmab).(J^{\ep}\etab).\nl(\ep^2J^{\ep}\an_{\ep}(\etab,a))}_{r,\qstar} \\
\\
&\le C\norm{J^{\ep}\etab}_{r,\qstar}\norm{(J^{\ep}\sigmab).\nl(\ep^2J^{\ep}\an_{\ep}(\etab,a))}_{r,\qstar} \\
\\
&\le \M_r[\norm{\etab}_{r,\qstar}]\left(1+|a|\ep^{1-r}\right)\norm{\etab}_{r,\qstar},
\end{align*}
using Proposition \ref{mapping prop} for the last inequality.

\subsubsection{Bootstrap estimates for $\j_{32}$} We begin with
\begin{align*}
\norm{\j_{32}^{\ep}(\etab,a)}_{r+1,\qstar}
&\le C|a|\norm{\cb^{\ep}(\sigmab,\varphib_{\ep}^a,\an_{\ep}(\etab,a))}_{r,\qstar} \\
\\
&\le C|a|\norm{(J^{\ep}\sigmab).(J^{\ep}\varphib_{\ep}^a).\nl(\ep^2J^{\ep}\an_{\ep}(\etab,a))}_{r,\qstar} \\
\\
&\le C|a|\ep^{-r}\norm{(J^{\ep}\sigmab).\nl(\ep^2J^{\ep}\an_{\ep}(\etab,a))}_{r,\qstar}.
\end{align*}
Expanding $\nl$ into its product form, we find
\begin{align*}
\norm{\j_{32}^{\ep}(\etab,a)}_{r+1,\qstar}
&\le C|a|\ep^{2-r}\norm{(J^{\ep}\sigmab).(J^{\ep}(\sigmab+\etab)).N(\ep^2J^{\ep}\an_{\ep}(\etab,a))}_{r,\qstar} \\
\\
&+ C|a|\ep^{2-r}\norm{(J^{\ep}\sigmab).(aJ^{\ep}\varphib_{\ep}^a).N(\ep^2J^{\ep}\an_{\ep}(\etab,a))}_{r,\qstar}.
\end{align*}
We estimate these terms as we did $\Pi_1$ and $\Pi_2$ in \eqref{pi1 bootstrap} and \eqref{pi2 bootstrap} and find
$$
\norm{\j_{32}^{\ep}(\etab,a)}_{r+1,\qstar}
\le \M_r[\norm{\etab}_{r,\qstar}]\left(|a|\ep^{1-r}+a^2\ep^{1-2r}\right).
$$

\subsubsection{Bootstrap estimates for $\j_{42}$} Routine estimates give
\begin{align*}
\norm{\j_{42}^{\ep}(\etab,a)}_{r+1,\qstar}
&\le C|a|\norm{\cb^{\ep}(\etab,\varphib_{\ep}^a,\an_{\ep}(\etab,a))}_{r,\qstar} \\
\\
&\le C|a|\norm{(J^{\ep}\etab).(J^{\ep}\varphib_{\ep}^a).\nl(\ep^2J^{\ep}\an_{\ep}(\etab,a))}_{r,\qstar} \\
\\
&\le C|a|\ep^{-r}\norm{(J^{\ep}\etab).\nl(\ep^2J^{\ep}\an_{\ep}(\etab,a))}_{r,\qstar}.
\end{align*}
From here we follow the $\j_{32}$ estimates with $\sigmab$ replaced by $\etab$.

\subsubsection{Bootstrap estimates for $\j_{52}$} Straightforward estimates and one invocation of Proposition \ref{mapping prop} give
\begin{align*}
\norm{\j_{52}^{\ep}(\etab,a)}_{r+1,\qstar}
&\le C\norm{\cb^{\ep}(\etab,\etab,\an_{\ep}(\etab,a))}_{r,\qstar} \\
\\
&\le C\norm{(J^{\ep}\etab)^{.2}.\nl(\ep^2J^{\ep}\an_{\ep}(\etab,a))}_{r,\qstar} \\
\\
&\le \M_r[\ep^2\norm{J^{\ep}(\sigmab+\etab)}_{r,\qstar}]\left(1+|a|\ep^{1-r}\right)\norm{(J^{\ep}\etab)^{.2}}_{r,\qstar} \\
\\
&\le \M_r[\norm{\etab}_{r,\qstar}]\left(1+|a|\ep^{1-r}\right)\norm{\etab}_{r,\qstar}^2.
\end{align*}

\subsubsection{Bootstrap estimates for $\j_{6}$} We rely on Proposition \ref{general hrq lipschitz prop}:
\begin{align*}
\norm{\j_6^{\ep}(\etab,a)}_{r+1,\qstar}
&\le Ca^2\norm{\cb^{\ep}(\varphib_{\ep}^a,\varphib_{\ep}^a,\an_{\ep}(\etab,a))-\cb^{\ep}(\varphib_{\ep}^a,\varphib_{\ep}^a,a\varphib_{\ep}^a)}_{r,\qstar} \\
\\
&\le Ca^2\norm{(J^{\ep}\varphib_{\ep}^a)^{.2}.\left(\nl(\ep^2J^{\ep}\an_{\ep}(\etab,a))-\nl(\ep^2J^{\ep}(a\varphib_{\ep}^a))\right)}_{r,\qstar} \\
\\
&\le Ca^2\ep^{-2r}\norm{\nl(\ep^2J^{\ep}\an_{\ep}(\etab,a))-\nl(\ep^2J^{\ep}(a\varphib_{\ep}^a))}_{r,\qstar} \\
\\
&\le Ca^2\ep^{2-2r}\M_r[\norm{J^{\ep}(\sigmab+\etab)}_{r,\qstar}]\left(1+|a|\ep^{1-r}\right)\norm{J^{\ep}(\sigmab+\etab)}_{r,\qstar}.
\end{align*}
Taking the supremum over $\ep \in (0,\epbar)$ in the $\M_r$ factor, we conclude
$$
\norm{\j_6^{\ep}(\etab,a)}_{r+1,\qstar}
\le \M_r[\norm{\etab}_{r,\qstar}]a^2\ep^{2-2r}\left(1+|a|\ep^{1-r}\right)
\le \M_r[\norm{\etab}_{r,\qstar}]\left(a^2\ep^{1-2r} + |a|^3\ep^{1-3r}\right).
$$

%% file: appendix_operator_estimates.tex
\section{Operator Estimates}

In this appendix we present two long collections of results that we proved earlier in \cite{faver-wright}. The set-up there is slightly different, of course; for example, the eigenvector operators $J$ and $J_1$ are not the same, and the ``symmetry'' in the mass dimer problem analogous to the ``even $\times$ even'' symmetry in Lemma \ref{even even symmetry lemma} is ``even $\times$ odd.''  But the proofs, happily, do not depend on these superficial differences, and so we just provide the long litany of estimates and other properties below.

\subsection{Estimates for the periodic problem} 

\begin{proposition}\label{master periodic lemma}
There exists $\ep_{\per} \in (0,1)$ with the following properties.

\begin{enumerate}[label={\bf(\roman*)}]
\item For $\ep \in (0,\ep_{\per})$ and $t \in \R$, there is a function $\rhs_{\ep}$ such that
\begin{equation}\label{xi symbol}
\txi_{\cep}(\ep\omega_{\ep}+t)
= \bunderbrace{\txi_{\cep}'(\ep\omega_{\ep})}{\Upsilon_{\ep}}t + t^2\rhs_{\ep}(t).
\end{equation}
The constant $b_0 > 0$ from \eqref{omegac derivative bound} also satisfies the estimate
\begin{equation}\label{upsilon bound}
|\Upsilon_{\ep}| \ge b_0
\end{equation}
for all $\ep \in (0,\ep_{\per})$.

\item For $|t| \le 1$ and $\ep \in (0,\ep_{\per})$, the multiplier $\xi^{\ep,t}$ maps $E_{\per}^{r+2}$ bijectively onto the space 
$$
\set{\psi \in E_{\per}^r}{\hat{\psi}(1) = 0}. 
$$
In particular, if $\Pi_2$ is the multiplier defined in Appendix \ref{fixed point conversion appendix} with symbol $\tilde{\Pi}_2(k) = 1-\delta_{|k|,1}$, then $\xi^{\ep,t}$ is invertible on the range of $\Pi_2$.  

\item 
There exists $C_{\map} > 0$ such that 
\begin{equation}\label{L1 mapping estimate}
\sup_{\substack{0 < \ep < \ep_{\per} \\
|t| \le 1 \\
r \in \R}}
\norm{\L_1^{\ep}(t)}_{\b(\W^r,\W^{r+2})} \le C_{\map},
\end{equation}
\begin{equation}\label{L2 L3 L4 mapping estimate}
\sup_{\substack{0 < \ep < \ep_{\per} \\
|t| \le 1 \\
r \in \R \\
j = 2,3}}
\norm{\L_j^{\ep}(t)}_{\b(\W^r,\W^r)} \le C_{\map},
\end{equation}
and
\begin{equation}\label{txi mapping estimate}
\sup_{\substack{0 < \ep < \ep_{\per} \\ t \in \R}} |\rhs_{\ep}(t)| \le C_{\map}.
\end{equation}

\item 
There exists $C_{\lip} > 0$ such that 
\begin{equation}\label{L1 lipschitz estimate}
\sup_{\substack{0 < \ep < \ep_{\per} \\
r \in \R}}
\norm{\L_1^{\ep}(t) - \L_1^{\ep}(\grave{t})}_{\b(\W^r,\W^r)}
\le C_{\lip}|t-\grave{t}|,
\end{equation}
\begin{equation}\label{L2 L3 L4 Lipschitz estimate}
\sup_{\substack{0 < \ep < \ep_{\per} \\
r \in \R \\
j = 2,3,4}}
\norm{\L_j^{\ep}(t)-\L_j^{\ep}(\grave{t})}_{\b(\W^r,\W^{r-1})} \le C_{\lip}|t-\grave{t}|,
\end{equation}
and
\begin{equation}\label{txi lipschitz estimate}
\sup_{0 < \ep < \ep_{\per}} |\rhs_{\ep}(t) - \rhs_{\ep}(\grave{t})| \le C_{\lip}|t-\grave{t}|.
\end{equation}
for any $|t|$, $|\grave{t}| \le 1$. 
\end{enumerate}
\end{proposition}

\subsection{Estimates for the nanopteron equations} \label{nanopteron estimates from FW}

\begin{proposition}\label{faver wright nanopteron}
There exists $q_{\star} \in (0,1/2\sqrt{\alpha_{\kappa}})$ and $\epbar \in (0,\ep_{\per})$ such that the following hold.

\begin{enumerate}[label={\bf(\roman*)}]
\item For all $r \ge 0$, there exists $C_r > 0$ such that for all $\ep \in (0,\epbar)$ and all $a$, $\grave{a} \in [-a_{\per},a_{\per}]$, the periodic solutions $\varphib_{\ep}^a$ defined in Theorem \ref{periodic main theorem} satisfy
\begin{equation}\label{varphib Wrinfty estimate}
\norm{\varphib_{\ep}^a}_{W^{r,\infty}} + \norm{J^{\ep}\varphib_{\ep}^a}_{W^{r,\infty}} \le C_r\ep^{-r}
\end{equation}
and
\begin{equation}\label{varphib lipschitz estimate}
|\partial_X^r[J^{\ep}(\varphib_{\ep}^a-\varphib_{\ep}^{\grave{a}})](X)|
\le C_r\ep^{-r}|a-\grave{a}|(1+|X|), \ X \in \R.
\end{equation}

\item There exists $C > 0$ such that for all $r \ge 0$, $q > 0$, $\ep \in (0,\epbar)$, and and $f \in H_q^r$, the operator $\iota_{\ep}$ defined in \eqref{iota ep defn} satisfies
\begin{equation}\label{iota ep estimate}
|\iota_{\ep}[f]|
\le \frac{C\ep^r}{\sqrt{q}}\norm{f}_{r,q}.
\end{equation}

\item There exists $C > 0$ such that the quantity $\upsilon_{\ep}$, defined in \eqref{upsilon ep defn}, satisfies
\begin{equation}\label{upsilon ep estimate}
|\upsilon_{\ep}| \ge C
\end{equation}
for all $\ep \in (0,\epbar)$.

\item The operator $\T_{\ep}$ defined in \eqref{T ep defn} has the following properties.

\begin{itemize}[leftmargin=*]
\item Let $r \ge 0$ and $q \in [0,q_{\star}]$.  Given $g \in H_q^r$, there exists $f \in H_q^{r+2}$ such that $\T_{\ep}f = g$ if and only if $\hat{g}(\pm\omega_{\ep}) = 0$, in which case $f$ is unique;
\item Let $r \ge 0$, $q \in [0,q_{\star}]$, and $\ep \in (0,\epbar)$.  For all $g \in E_q^r$, there exists a unique $f \in E_q^{r+2}$ such that 
\begin{equation}\label{Pep defn}
\T_{\ep}f 
= g - \frac{1}{\upsilon_{\ep}}\iota_{\ep}[g]\chi_{\ep}.
\end{equation} 
As in \eqref{Pep defn orig}, we set $f := \P_{\ep}g$.  Equivalently,
$$
\P_{\ep}g := \T_{\ep}^{-1}\left(g - \frac{1}{\upsilon_{\ep}}\iota_{\ep}[g]\chi_{\ep}\right).
$$
\item For each $q \in [0,q_{\star}]$, there exists $C_q > 0$ such that
\begin{equation}\label{Pep estimate}
\norm{\P_{\ep}g}_{r+j,q} \le \frac{C_q}{\ep^{j+1}}\norm{g}_{r,q}, \ j = 0 ,1,2
\end{equation}
for all $r \ge 0$ and $\ep \in (0,\epbar)$.  
\end{itemize}

\item There exists $C> 0$ such that for all $r \ge 0$, $\ep \in (0,\epbar)$, $q \in (0,q_{\star}]$, and $\hb \in H_q^r \times H_q^r$
\begin{equation}\label{Jep0 estimate}
\norm{(J^{\ep}-J^0)\hb}_{r,q} \le C\ep\norm{\hb}_{r+1,q}.
\end{equation}

\item There exists $C > 0$ such that for all $\ep \in (0,\epbar)$, $q \in (0,\qstar]$, $r \ge 0$, and $f \in H_q^r$, the Fourier multiplier $\varpi^{\ep}$, whose symbol is given in \eqref{tvarpi ep defn}, satisfies
\begin{equation}\label{varpi ep estimate}
\norm{\varpi^{\ep}f}_{r+2,q} \le C\norm{f}_{r,q}.
\end{equation}

\item For each $r \ge 0$ and $q \in [0,q_{\star}]$ the operator $\A$ defined in \eqref{A defn} is bijective from $E_q^r$ to $E_q^r$.  

\item For each $r$, $q \ge 0$, the operator $\K_2$ defined in \eqref{death star defns} is bounded from $E_q^r$ to $E_q^r$.

\item Let $\mu$ be one of the operators $J$, $J_1$, $\lambda_{\pm}$.  Then with $\mu^{\ep}$ defined in Appendix \ref{fourier multipliers appendix}, the operator norms $\norm{\mu^{\ep}}$ satisfy
$$
\sup_{0 < \ep < \epbar} \norm{\mu^{\ep}} < \infty
$$
for each $q$, $r \ge 0$.
\end{enumerate}

\end{proposition}